\documentclass[reqno,12pt]{amsart}
\usepackage[dvips]{color}
\usepackage[hypertex]{hyperref}

\usepackage[utf8]{inputenc}

\usepackage{amsfonts}
\usepackage{amsthm}
\usepackage{amssymb}
\usepackage{texdraw}

\usepackage[all]{xy}
\allowdisplaybreaks

\numberwithin{equation}{section}

\newtheorem{Theorem}{Theorem}[section]
\newtheorem{Proposition}[Theorem]{Proposition}
\newtheorem{Lemma}[Theorem]{Lemma}
\newtheorem{Definition}[Theorem]{Definition}
\newtheorem{Corollary}[Theorem]{Corollary}

\theoremstyle{definition}
\newtheorem{Remark}[Theorem]{Remark}
\newtheorem{Example}[Theorem]{Example}

\makeatletter
\@namedef{subjclassname@2020}{\textup{2020} Mathematics Subject Classification}
\makeatother

\begin{document}
\vspace*{-2cm}
\begin{center}
\noindent{\bf This paper will appear in {\em Annales de l'Institut Fourier}}
\vspace*{\baselineskip}
\end{center}

\vspace*{2cm}

\title{Jacobian curve of  singular foliations}
\author{Nuria Corral}
%\date{\today}
\address{Nuria Corral, Departamento de Matemáticas, Estadística y Computación, Universidad de Cantabria, Avda. de los Castros s/n, 39005 -- Santander, SPAIN}
\email{nuria.corral@unican.es}
\dedicatory{Dedicated to Felipe Cano, with admiration and gratitude}
\thanks{The author is supported by the Spanish research project  PID2019-105621GB-I00}
\maketitle
\begin{abstract} Topological properties of the jacobian curve  ${\mathcal J}_{\mathcal{F},\mathcal{G}}$ of two foliations $\mathcal{F}$ and $\mathcal{G}$ are des\-cribed in terms of invariants associated to the foliations. The main result gives a decomposition of the jacobian curve ${\mathcal J}_{\mathcal{F},\mathcal{G}}$ which depends on how similar are the foliations $\mathcal{F}$ and $\mathcal{G}$. The similarity between foliations is codified in terms of the Camacho-Sad indices of the foliations with the notion of collinear point or divisor. Our approach allows to recover the results concerning the factorization of the jacobian curve of two plane curves and of the polar curve of a curve or a foliation.
\end{abstract}

\section{Introduction}
Given two germs of holomorphic functions $f,g \in {\mathbb C}\{x,y\}$,  the Jacobian determinant
$$J(f,g)=f_x g_y-f_yg_x$$
defines a curve called the {\em jacobian curve} of $f$ and $g$ (see \cite{Mau,Cas-07} for instance). The analytic type of the jacobian curve is an invariant of the analytic type of the pair of curves $f=0$ and $g=0$ but its topological type is not a topological invariant of the pair of curves (see \cite{Mau99}). Properties of the jacobian curve  have been studied by several authors in terms of properties of the curves $f=0$ and $g=0$ (see for instance \cite{Kuo-P-2004,Cas-07},  and  \cite{Gar-G} when $g$ is a characteristic approximated root of $f$).

This notion can be studied in the more general context given by the theory of singular foliations: given two germs of foliations $\mathcal F$ and $\mathcal G$
 in $({\mathbb C}^2,0)$,
defined by the 1-forms $\omega=0$ and $\eta=0$, the {\em jacobian curve} ${\mathcal J}_{{\mathcal F},{\mathcal G}}$   of $\mathcal F$ and $\mathcal G$ is the curve given by
$$\omega \wedge \eta=0.$$
Note that this is the curve of tangency between both foliations. It is easy to show that the branches of ${\mathcal J}_{{\mathcal F},{\mathcal G}}$ are not separatrices of $\mathcal F$ or $\mathcal G$ provided that the foliations $\mathcal F$ and $\mathcal G$ do not have  common separatrices.

If the foliation $\mathcal G$ is non-singular, the jacobian curve ${\mathcal J}_{{\mathcal F},{\mathcal G}}$  coincides with the polar curve of the foliation $\mathcal F$. Properties of the equisingularity type of  polar curves of foliations have been studied in \cite{Rou,Cor-2003,Cor-2009,Cor-2009-BullBraz}. Moreover, if the foliation $\mathcal F$ is given by $df=0$ with $f \in {\mathbb C}\{x,y\}$, we recover the notion of polar curve of a plane curve.  The local study of these curves has also been widely treated by many authors (see for instance \cite{Mer,Cas-00,Gar,Le-M-W} or the recent works \cite{Alb-G,Hef-H-I}).

Moreover, the use of polar curves of foliations allowed to describe properties of foliations. In \cite{Can-C-M}, the  study of intersection properties of polar curves of foliations permitted to characterize generalized curve foliations as well as second type foliations; an expression of the GSV-index can also be given in terms of these invariants (see also \cite{Gen-M} for the dicritical case). There are also some recent works that show the interest of  jacobian curves or polar curves of foliations in the study of analytic invariants of curves  (see for instance \cite{Gom}) or singular foliations (see \cite{Ort-R-V}).

\medskip
The aim of  this paper is to describe properties of the equisingularity type of ${\mathcal J}_{{\mathcal F},{\mathcal G}}$ in terms of invariants associated to the foliations $\mathcal F$ and $\mathcal G$. Note that, in general, the locus ${\mathcal J}_{{\mathcal F},{\mathcal G}}$ cannot be described from the data of $\mathcal{F}$ and $\mathcal{G}$. It is enough to consider the non-singular foliations $\mathcal{F}$ given by $dx=0$ and $\mathcal{G}$ defined by $dx+h(x,y)dy=0$,  hence the jacobian curve ${\mathcal J}_{{\mathcal F},{\mathcal G}}$ is defined by $h(x,y)=0$.

  To illustrate the kind of conditions we are going to ask to the foliations and the type of results that we can obtain, let us explain the relationship between the multiplicity at the origin  of the jacobian curve $\nu_0({\mathcal J}_{{\mathcal F},{\mathcal G}})$ and the multiplicities of the foliations $\mathcal{F}$ and $\mathcal{G}$.
If the 1-forms $\omega$ and $\eta$ defining $\mathcal{F}$ and $\mathcal{G}$ are given by $\omega=A(x,y) dx + B(x,y) dy$ and $\eta=P(x,y) dx + Q(x,y) dy$ respectively, the jacobian curve ${\mathcal J}_{{\mathcal F},{\mathcal G}}$  is defined by $J(x,y)=0$ where
\begin{equation} \label{eq:jacobiano}
J(x,y)=\left| \begin{array}{cc}
A(x,y) & B(x,y) \\
P(x,y) & Q(x,y)
\end{array} \right|.
\end{equation}
Thus the multiplicity at the origin $\nu_0({\mathcal J}_{{\mathcal F},{\mathcal G}})$ of the jacobian curve satisfies
\begin{equation}\label{eq:multiplicidad}
\nu_0({\mathcal J}_{{\mathcal F},{\mathcal G}}) \geq \nu_0({\mathcal F}) + \nu_0({\mathcal G})
\end{equation}
where $\nu_0( {\mathcal F} ), \ \nu_0({\mathcal G})$ denote the multiplicity at the origin of the foliations $\mathcal F$ and $\mathcal G$ respectively.
One of the first results  describing  the properties of the jacobian curve shows that equality in equation~\eqref{eq:multiplicidad} holds, that is,
$$
\nu_0({\mathcal J}_{{\mathcal F},{\mathcal G}}) = \nu_0({\mathcal F}) + \nu_0({\mathcal G})
$$
  provided that the foliations $\mathcal{F}$ and $\mathcal{G}$ have different Camacho-Sad index at  any singular point  in the exceptional divisor $E^1$ obtained after one blow-up (see Lemma~\ref{lema:multiplicidad-E1}).

The main result in this paper, Theorem \ref{th:desc-general}, gives a factorization of the jacobian curve ${\mathcal J}_{\mathcal{F},\mathcal{G}}$ of two generalized curve foliations $\mathcal F$ and $\mathcal G$ in terms of invariants given by the dual graph of the common minimal reduction of singularities of $\mathcal F$ and $\mathcal G$. This result gives a decomposition of the jacobian curve ${\mathcal J}_{{\mathcal F},{\mathcal G}}$ in two classes of components: one for which we can control some properties of the topology from the data of $\mathcal{F}$ and $\mathcal{G}$ and another one for which such a control is impossible.
The properties of the components in that decomposition depend on how ``similar'' are the foliations $\mathcal F$ and $\mathcal G$ in terms of its singularities and Camacho-Sad indices at the common singularities. We introduce the notion of collinear point and collinear divisor to measure this similarity between the foliations (see section~\ref{sec:collinear} where properties of collinear and non-collinear divisors are given).

The strategy used to prove the decomposition result is to study first the case when the foliations $\mathcal{F}$ and $\mathcal{G}$ have separatrices with non-singular irreducible components (section~\ref{sec:jacobian}). In this case, thanks to the hypothesis over the separatrices, we can compute ``by hand'' the infinitely near points of the jacobian curve under certain hypothesis over the foliations related with the notion of collinearity (a key point is Lemma~\ref{lema:parte-inicial-jac} relating the weighted initial part of the 1-forms defining the foliations $\mathcal{F}$ and $\mathcal{G}$ and the one of the equation of the jacobian curve).
%that, for non-collinear divisors, the weighted initial part of the equation defining the jacobian curve can be obtained form the ones of the 1-forms defining the foliations (see Lemma~\ref{lema:parte-inicial-jac}).
In these computations, we use the existence of logarithmic models for generalized curve foliations (proved in \cite{Cor-2003}) and the properties shared by a foliation and its logarithmic model. Moreover, we describe the relationship between the jacobian curve of two foliations and the one of its logarithmic models (see Lemmas~\ref{lema:parte-incial-modelo-log} and \ref{lema:puntos-jac-log-nolog}). These results will allow us to do some computations for the jacobian curve of two logarithmic foliations (see Theorem~\ref{th-multiplicidad-jac}) and thus we get it for the jacobian curve of any non-dicritical generalized curve foliations.

Then we use a ramification $\rho:(\mathbb{C}^2,0) \to (\mathbb{C}^2,0)$ to reduce the general case to the previous one (section~\ref{sec:caso-general}). This strategy works since we can prove that the curves $\rho^{-1} \mathcal{J}_{\mathcal{F},\mathcal{G}}$ and $\mathcal{J}_{\rho^*\mathcal{F},\rho^*\mathcal{G}}$ ``share'' the same infinitely near points in the common reduction of singularities of $\rho^*\mathcal{F}$ and $\rho^*\mathcal{G}$ (see Lemma~\ref{lemma:ptos-inf-prox-ram}). Thus the results obtained in section~\ref{sec:jacobian} can be used to describe properties of $\rho^{-1} \mathcal{J}_{\mathcal{F},\mathcal{G}}$ and hence, recover properties of the curve $\mathcal{J}_{\mathcal{F},\mathcal{G}}$ since the equisingularity data of a curve can be recovered from the one obtained after ramification (see~\cite{Cor-2009}). We include an appendix (appendix~\ref{ap:ramificacion}) devoted to explain all the details concerning the ramification process.

Section \ref{sec:local-invariants} is devoted to introduce notations and local invariants of curves and foliations which will be used throughout the paper. In section~\ref{sec:logarithmic} we recall the notion of logarithmic model (introduced in \cite{Cor-2003}) and some properties of logarithmic foliations.

In the last part of the article (section~\ref{sec:polares}) we show the role that the Camacho-Sad indices play  to explain some behaviours of jacobian curves of plane curves. In particular, we show  how our results  imply the results of T.-C. Kuo and A. Parusi\'nski concerning jacobian curves of plane curves \cite{Kuo-P-2004}, the results of E. García Barroso and J. Gwo\'{z}dziewicz about the jacobian curve of a plane curve and its approximate roots \cite{Gar-G} and also previous results about polar curves of foliations (given in \cite{Rou,Cor-2003}). All these results can be consider as particular cases of the results in this paper.

The article finishes with two appendices. The first one contains results concerning ramification. The second one is devoted to prove some formulas which describe the multiplicity of intersection of the jacobian curve with the separatrices of the foliations $\mathcal F$ and $\mathcal G$ in terms of the local invariants associated to $\mathcal F$ and $\mathcal G$. These formulas generalize some properties of polar curves of a foliation given in \cite{Cor-2003,Can-C-M} which were key in the proof of the characterization of generalized curve foliations and second type foliations given in \cite{Can-C-M}.
\medskip

\noindent{\bf Acknowledgements.} The author thanks the referee for his carefully reading of the paper and for his contribution to a good presentation of it.

\section{Local invariants }\label{sec:local-invariants}
\subsection{Foliations}\label{subsec:inv-foliations}
Let $\mathbb F$ be the space of singular foliations in $({\mathbb C}^2,0)$. An element ${\mathcal F} \in {\mathbb F}$ is defined by a 1-form $\omega=0$, with $\omega=A(x,y)dx+B(x,y)dy$, or by the vector field ${\bf v} = -B(x,y) \partial/\partial x + A(x,y) \partial/\partial y$ where
 $A, B \in {\mathbb C}\{x,y\}$ are relatively prime. The origin is a singular point if $A(0)=B(0)=0$.
The {\em multiplicity\/} $\nu_0({\mathcal F})$ of $\mathcal F$ at the origin is the minimum of the orders $\nu_0(A)$, $\nu_0(B)$ at the origin. Thus, the origin is a {\em singular point\/} of $\mathcal F$ if $\nu_0({\mathcal F}) \geq 1$.

Consider a germ of irreducible analytic curve $S$ at $({\mathbb C}^2,0)$. We say that $S$ is a {\em separatrix\/} of $\mathcal F$ at the origin if $S$ is an invariant curve of the foliation $\mathcal F$. Therefore, if $f=0$ is a reduced equation of $S$, we have that $f$ divides $\omega \wedge df$.

Let us now recall the desingularization process of a foliation. We say that the origin is a {\em simple singularity} of $\mathcal F$ if there are local coordinates $(x,y)$ in $({\mathbb C}^2,0)$ such that $\mathcal F$ is given by a $1$-form of the type
$$\lambda y dx - \mu x dy + \text{ h.o.t}$$
with $\mu \neq 0$ and $\lambda/\mu \not\in {\mathbb Q}_{>0}$. If $\lambda=0$, the singularity is called a {\em saddle-node\/}.
There are two formal invariant curves $\Gamma_x$ and $\Gamma_y$ which are tangent to $x=0$ and $y=0$ respectively, and such that they are both convergent in the case that $\lambda \mu \neq 0$. In the saddle-node situation with $\lambda=0$ and $\mu\neq 0$, we say  that the saddle-node is {\em well oriented\/} with respect to the curve $\Gamma_y$.

 Let $\pi_1: X_1 \to ({\mathbb C}^2,0)$ be the blow-up of the origin with $E_1=\pi_1^{-1}(0)$ the exceptional divisor. We say that the blow-up $\pi_1$ (or the exceptional divisor $E_1$) is {\em non-dicritical\/} if $E_1$ is invariant by the strict transform $\pi_1^*{\mathcal F}$ of $\mathcal F$; otherwise, the exceptional divisor $E_1$ is generically transversal to $\pi_1^* {\mathcal F}$ and we say that the blow-up $\pi_1$ (or the divisor $E_1$) is {\em dicritical}.

A {\em reduction of singularities\/} of $\mathcal F$ is a morphism $\pi: X \to ({\mathbb C}^2,0)$, composition of a finite number of punctual blow-ups, such that the strict transform $\pi^* {\mathcal F}$ of $\mathcal F$ verifies that
\begin{itemize}
  \item each irreducible component of the exceptional divisor $\pi^{-1}(0)$ is either invariant by $\pi^*{\mathcal F}$ or transversal to $\pi^*{\mathcal F}$;
  \item all the singular points of $\pi^*{\mathcal F}$ are simple and do not belong to a dicritical component of the exceptional divisor.
\end{itemize}
There exists a reduction of singularities as a consequence of Seidenberg's Desingularization Theorem \cite{Sei}. Moreover, there is a minimal morphism $\pi$ such that any other reduction of singularities of $\mathcal F$ factorizes through the minimal one. The centers of the blow-ups of a reduction of singularities of $\mathcal F$ are called {\em infinitely near points\/} of $\mathcal F$.
If all the irreducible components of the exceptional divisor are invariant by $\pi^*{\mathcal F}$ we say that the foliation $\mathcal F$ is {\em non-dicritical}; otherwise, $\mathcal F$ is called a {\em dicritical foliation}.

A non-dicritical foliation $\mathcal F$ is called a {\em generalized curve foliation\/} if there are not saddle-node singularities in the reduction of singularities (see \cite{Cam-S-LN}).
We will denote $\mathbb G$ the space of non-dicritical generalized curve foliations in $({\mathbb C}^2,0)$. The foliation $\mathcal F$ is of {\em second type} if all saddle-nodes of $\pi^*{\mathcal F}$ are well oriented with respect to the exceptional divisor $E=\pi^{-1}(0)$ (see \cite{Mat-S}).

In order to describe properties of generalized curve foliations and second type foliations, let us recall some   local invariants used in the local study of foliations in dimension two (see for instance \cite{Can-C-D}). The {\em Milnor number} $\mu_0({\mathcal F})$ is given by
$$\mu_0({\mathcal F}) = \text{dim}_{\mathbb C}\frac{{\mathbb C}\{x,y\}}{(A,B)}= (A,B)_0,$$
where $(A,B)_0$ stands for the intersection multiplicity.  Note that, if the foliation is defined by $df=0$, the Milnor number of the foliation coincides with the one of the curve given by $f=0$. Given an irreducible curve $S$ and a primitive parametrization $\gamma: ({\mathbb C},0) \to ({\mathbb C}^2,0)$ of $S$ with $\gamma(t)=(x(t),y(t))$, we have that $S$ is a separatrix of $\mathcal F$ if and only if $\gamma^* \omega=0$. In this case,  the {\em Milnor number} $\mu_0({\mathcal F},S)$ of $\mathcal F$ {\em along\/} $S$ is given by
%$$\mu_0({\mathcal F},S)= \text{ord}_t {\bf w}(t),
%$$
%where ${\bf w}(t)$ is the unique vector field at $({\mathbb C},0)$ such that $\gamma_* {\bf w}(t)={\bf v} \circ \gamma(t)$. We have that
\begin{equation}
\label{multiplicidaderelativa1}
 \mu_0({\mathcal F},S) =
\begin{cases}
 {\rm ord}_{t}(B(\gamma(t))) -  {\rm ord}_{t}(x(t)) + 1 \ \ \ \mbox{if $x(t) \neq 0$}   \\
 {\rm ord}_{t}(A(\gamma(t))) -  {\rm ord}_{t}(y(t)) + 1 \ \ \ \mbox{if $y(t) \neq 0$}
\end{cases}
\end{equation}
(this number is also called {\em multiplicity of $\mathbf{v}$ along $S$}, see~\cite{Cam-S-LN}, p. 152-153). %, definition in page 152 and proposition 3 in page 153).
If $S$ is not a separatrix, we define the {\em tangency order} $\tau_0({\mathcal F},S)$ by
\begin{equation}\label{def-tangencia}
\tau_0({\mathcal F},S)= \text{ord}_t (\alpha(t))
\end{equation}
where $\gamma^* \omega= \alpha(t) dt$.
If $S=(y=0)$ is a non-singular invariant curve of the foliation $\mathcal F$,  the {\em Camacho-Sad index of $\mathcal F$ relative to $S$ at the origin} is given by
\begin{equation}\label{eq:indice-C-S}
  \mathcal{I}_0({\mathcal F},S) = -\text{Res}_0 \frac{a(x,0)}{b(x,0)}
\end{equation}
where the 1-form defining $\mathcal F$ is written as $y a(x,y) dx + b(x,y) dy$ (see \cite{Cam-S}).

Next result summarizes some of the properties of second type and generalized curve foliations that we will use throughout the text:
\begin{Theorem}\cite{Cam-S-LN,Mat-S,Can-C-M} \label{th-curva-gen}
Let $\mathcal F$ be a non-dicritical foliation and consider ${\mathcal G}_f$ the foliation defined by $df=0$ where $f$ is  a reduced equation of the curve $S_{\mathcal F}$ of separatrices of $\mathcal F$. Let $\pi: X \to ({\mathbb C}^2,0)$ be the minimal reduction of singularities of $\mathcal F$.
\begin{itemize}
 \item[(i)] $\pi$ is a reduction of singularities of $S_{\mathcal F}$. Moreover, $\pi$ is the minimal reduction of singularities of the curve $S_{\mathcal F}$ if and only if $\mathcal F$ is of second type;
 \item[(ii)]  $\nu_0({\mathcal F}) \geq \nu_0({\mathcal G}_f)$ and the equality holds if and only if $\mathcal F$ is of second type;
  \item[(iii)] $\mu_0({\mathcal F}) \geq \mu_0({\mathcal G}_f)$ and the equality holds if and only if $\mathcal F$ is a generalized curve;
  \item[(iv)] if $S$ is an irreducible curve which is not a separatrix of $\mathcal F$, then $\tau_0({\mathcal F},S) \geq \tau_0({\mathcal G}_f,S) $ and the equality holds if and only if $\mathcal F$ is of second type.
\end{itemize}
\end{Theorem}
Recall that for the hamiltonian foliation ${\mathcal G}_f$ we have that $\nu_0({\mathcal G}_f)= \nu_0(S_{\mathcal F})-1$,  $\mu_0({\mathcal G}_f) =\mu_0(S_{\mathcal F})$ and $\tau_0({\mathcal G}_f,S)=(S_{\mathcal F},S)_0 -1$ where $\nu_0(S_{\mathcal F})$ is the multiplicity of the curve $S_{\mathcal F}$ at the origin, $\mu_0(S_{\mathcal F})$ is the Milnor number of the curve $S_\mathcal{F}$ and $(S_{\mathcal F},S)_0$ denotes the intersection multiplicity of the curves $S_{\mathcal F}$ and $S$ at the origin.

\medskip
\paragraph{\bf Notation.} Given a plane curve $C$ in $({\mathbb C}^2,0)$, we denote by ${\mathbb F}_C$ the sub-space of $\mathbb F$ composed by the foliations having $C$ as curve of separatrices and ${\mathbb G}_C$ the foliations of ${\mathbb F}_C$ which are generalized curve foliations.

\medskip
Moreover, for generalized curve foliations we have that
\begin{Lemma}\cite{Cor-2009-AnnAcBras}
Assume that $\mathcal F$ is a non-dicritical generalized curve foliation. Let $\pi: (X,P) \to ({\mathbb C}^2,0)$ be a morphism composition of a finite number of punctual blow-ups and take  an irreducible component $E$ of the exceptional divisor $\pi^{-1}(0)$ with $P \in E$. Then, the strict transforms $\pi^*{\mathcal F}$ and $\pi^* {\mathcal G}_f$ satisfy that
\begin{itemize}
  \item[1.] $\nu_P(\pi^*{\mathcal F})=\nu_P(\pi^* {\mathcal G}_f)$;
  \item[2.] $\mu_P(\pi^*{\mathcal F},E)=\mu_P(\pi^*{\mathcal G}_f,E)$.
\end{itemize}
where $f=0$ is a reduced equation of the curve $S_{\mathcal F}$ of separatrices  of $\mathcal F$.
\end{Lemma}

\subsection{Weighted initial forms and Jacobian curves.} \label{sec:parte-inicial} Fix  coordinates  $(x,y)$   in $({\mathbb C}^2,0)$. Given   a 1-form $\omega$, we can write $\omega=\sum_{i,j} \omega_{ij}$ where $\omega_{ij}=A_{ij}x^{i-1}y^{j} dx + B_{ij} x^i y^{j-1} dy$. We denote $\Delta(\omega)=\Delta(\omega;x,y)=\{(i,j) \ : \ \omega_{ij} \neq 0\}$ and    the {\em Newton polygon\/} ${\mathcal N}({\mathcal F};x,y)={\mathcal N}({\mathcal F})={\mathcal N}(\omega)$ is given by the convex envelop of $\Delta(\omega) + ({\mathbb R}_{\geq 0})^2$.

Given a rational number $\alpha \in {\mathbb Q}$, we define the {\em initial form\/} of $\omega$ {\em with weight $\alpha$}
$$\text{In}_{\alpha} (\omega;x,y) = \sum_{i+\alpha j =k} \omega_{ij}$$
where $i+\alpha j =k$ is the equation of the first line of slope $-1/\alpha$ which intersects the Newton polygon ${\mathcal N}({\mathcal F})$ in the coordinates $(x,y)$. Note that $k=\nu_{(1,\alpha)}(\omega)$ where  $\nu_{(1,\alpha)}(\omega)=\min \{i+\alpha j \ : \ \omega_{ij} \neq 0\}$ is  the $(1,\alpha)$-degree of $\omega$. Hence, we have that the multiplicity $\nu_0(\mathcal{F})=\nu_{(1,1)}(\omega)$.

In a similar way, given any function $f =\sum_{ij} f_{ij} x^i y^j \in {\mathbb C}\{x,y\}$, we denote $\Delta(f)=\Delta(f;x,y)=\{(i,j) \ : \ f_{ij} \neq 0\}$ and the Newton polygon ${\mathcal N}(C;x,y)={\mathcal N}(C)$ of the curve $C=(f=0)$ is the convex envelop of $\Delta(f) +({\mathbb R}_{\geq 0})^2$. Note that ${\mathcal N}(C;x,y)={\mathcal N}(df;x,y)$ and ${\mathcal N}({\mathcal F})={\mathcal N}(C)$ if ${\mathcal F} \in {\mathbb G}_{C}$.
Thus, we can define the {\em initial form\/} $\text{In}_{\alpha} (f;x,y)=\sum_{(i,j) \in L} f_{ij}x^i y^j$ where $L$ is the first line of slope $-1/\alpha$ which intersects ${\mathcal N}(C)$.
Note that, if $f=0$ is an equation of the curve $C$, then $\text{In}_{1} (f;x,y)$ gives an equation of the tangent cone of $C$, and hence $\text{In}_{1} (f;x,y)=\sum_{i+j=\nu_0(C)} f_{ij}x^i y^j$.

With these notations we can state the first result which illustrates the type of conditions we are going to ask the foliations $\mathcal{F}$ and $\mathcal{G}$ in order to be able to describe properties of the jacobian curve $\mathcal{J}_{\mathcal{F},\mathcal{G}}$.
\begin{Lemma}\label{lema:multiplicidad-E1}
Let $\mathcal{F}$ and $\mathcal{G}$ be two foliations in $(\mathbb{C}^2,0)$ and consider $\mathcal{J}_{\mathcal{F},\mathcal{G}}$ its jacobian curve. Let $\pi_1: X_1 \to (\mathbb{C}^2,0)$ be the blow-up of the origin and $E^1=\pi^{-1}_1(0)$ be the exceptional divisor. If there is a point $R \in E^1$ such that the Camacho-Sad indices verify that $\mathcal{I}_R(\pi_1^*\mathcal{F},E^1) \neq \mathcal{I}_R(\pi_1^*\mathcal{G},E^1)$, then
$$\nu_0(\mathcal{J}_{\mathcal{F},\mathcal{G}})=\nu_0(\mathcal{F}) + \nu_0(\mathcal{G}).$$
\end{Lemma}
\begin{proof}
Take $(x,y)$ coordinates such that $x=0$ is not tangent to the foliations $\mathcal{F}$ and $\mathcal{G}$ and let $(x_1,y_1)$ be coordinates in the first chart of the blow-up such that $\pi_1(x_1,y_1)=(x_1,x_1y_1)$ and $E^1=(x_1=0)$. Assume that $\nu_0(\mathcal{J}_{\mathcal{F},\mathcal{G}})>\nu_0(\mathcal{F}) + \nu_0(\mathcal{G})$, then $\text{In}_1(\omega) \wedge \text{In}_1(\eta)\equiv 0$. Thus, if  write
\begin{align*}
  \text{In}_1(\omega) & = A_{\nu_0(\mathcal{F})}(x,y) dx + B_{\nu_0(\mathcal{F})}(x,y) dy \\
   \text{In}_1(\eta) & = P_{\nu_0(\mathcal{G})}(x,y) dx + Q_{\nu_0(\mathcal{G})}(x,y) dy
\end{align*}
then we have that
\begin{equation}\label{eq:cono-tg}
 A_{\nu_0(\mathcal{F})}(x,y)   Q_{\nu_0(\mathcal{G})}(x,y)  - B_{\nu_0(\mathcal{F})}(x,y) P_{\nu_0(\mathcal{G})}(x,y) \equiv 0.
\end{equation}
  The computation of the Camacho-Sad index at a point $R$ given by $(0,c)$ in the coordinates $(x_1,y_1)$ gives
\begin{align*}
  \mathcal{I}_{R}(\pi_1^*\mathcal{F},E^1) & = -\text{Res}_{y=c} \frac{B_{\nu_0(\mathcal{F})}(1,y)}{ A_{\nu_0(\mathcal{F})}(1,y)+yB_{\nu_0(\mathcal{F})}(1,y)} \\
  \mathcal{I}_{R}(\pi_1^*\mathcal{G},E^1) & = -\text{Res}_{y=c} \frac{Q_{\nu_0(\mathcal{G})}(1,y)}{ P_{\nu_0(\mathcal{G})}(1,y)+yQ_{\nu_0(\mathcal{G})}(1,y)}
\end{align*}
The equality in \eqref{eq:cono-tg}, implies
$$ \mathcal{I}_{R}(\pi_1^*\mathcal{F},E^1) = \mathcal{I}_{R}(\pi_1^*\mathcal{G},E^1)$$
for  any point $R$ in $E^1$. This gives a contradiction with the hypothesis over the foliations $\mathcal{F}$ and $\mathcal{G}$.
\end{proof}

The condition over the Camacho-Sad indices of the foliations in the previous lemma is related with the notion of collinear point that will be introduced in section~\ref{sec:collinear}.

\subsection{Equisingularity data and dual graph of plane curves}\label{Apendice:grafo-dual}
In this subsection we will fix some notations concerning the equisingularity data of a plane curve $C=\cup_{i=1}^r C_i$ in $({\mathbb C}^2,0)$ (see appendix~\ref{ap:equisingularidad} for more details). Given an irreducible component $C_i$ of $C$, we will denote $n^i=\nu_0(C_i)$ the multiplicity of $C_i$ at the origin, $\{\beta_0^i,\beta_1^i,\ldots,\beta_{g_i}^i\}$ the characteristic exponents of $C_i$ and $\{(m_l^i,n_l^i)\}_{l=1}^{g_i}$ the Puiseux pairs of $C_i$.

Let us denote $\pi_C : X_C \to ({\mathbb C}^2,0)$ the minimal reduction of singularities of the curve $C$. The {\em dual
graph\/} $G(C)$ is constructed as follows:  each irreducible
component $E$ of the exceptional divisor $\pi_{C}^{-1}(0)$ is represented by a vertex which
we also call $E$ (we identify a divisor and its  associated vertex
in the dual graph). Two vertices are joined by an edge if and only
if the associated divisors intersect. Each irreducible component of
$C$ is represented by an arrow joined to the only divisor which
meets the strict transform of $C$ by $\pi_C$. We can give a weight to
each vertex $E$ of $G(C)$ equal to the  self-intersection of the
divisor $E \subset X_C$ and this weighted dual graph is equivalent to the
equisingularity data of $C$.

If we  denote by $E^1$ the irreducible component of $\pi_{C}^{-1}(0)$ corresponding to the divisor
obtained by the blow-up of the origin, we can give an orientation
 to the graph $G(C)$ beginning from the first divisor $E^1$. The {\em geodesic\/}
of a divisor $E$ is the path which joins the first divisor $E^1$
with the divisor $E$. The geodesic of a curve is the geodesic of the
divisor that meets the strict  transform of the curve. Thus,
there is a partial order in the set of vertices of $G(C)$ given by
$E < E'$ if and only if the geodesic of $E'$ goes through $E$. A {\em maximal divisor} in $G(C)$ will be a maximal element in the set of vertices of $G(C)$ with this partial order.
 Given a divisor $E$ of $G(C)$, we denote by $I_E$ the set of indices $i \in \{1,2,\ldots, r\}$ such that $E$ belongs to the geodesic of the curve $C_i$.

Given a vertex $E$ of $G(C)$, we define the number $b_E^C$ in the following way: $b_E^C +1$ is the valence of $E$ if $E \neq E^1$ and $b_{E^1}^C$ is the valence of $E^1$  in $G(C)$ (recall that the valence of a divisor $E$ in $G(C)$ is the number of arrows and edges attached to $E$ in $G(C)$). Given a divisor $E$ of $G(C)$, we say that $E$ is a {\em bifurcation divisor\/} of $G(C)$ if $b_E^C \geq 2$ and a {\em terminal divisor\/} of $G(C)$ if $b_E^C=0$. A {\em dead arc\/} is a path which joins a bifurcation divisor with a terminal one without going through other bifurcation divisor. We  denote by $B(C)$ the set of bifurcation divisors of $G(C)$. If there is no confusion with the curve $C$ we are working with, we will denote $b_E=b_E^C$ for any divisor $E$ in $G(C)$.

Given   an irreducible component $E$ of $\pi^{-1}_C(0)$, we denote by $\pi_E: X_E \to ({\mathbb C}^2,0)$ the {\em morphism reduction\/} of $\pi_C$ to $E$ (see \cite{Cor-2009}), that is, the morphism which verifies that
\begin{itemize}
  \item the morphism $\pi_C$ factorizes as $\pi_C= \pi_E \circ \, \pi_E'$ where $\pi_E$ and $\pi_E'$ are composition of punctual blow-ups;
  \item the divisor $E$ is the strict transform by $\pi_E'$ of an irreducible component $E_{red}$ of $\pi_E^{-1}(0)$ and $E_{red} \subset X_E$ is the only component of $\pi_E^{-1}(0)$ with self-intersection equal to $-1$.
\end{itemize}
We will denote by $\pi_E^* C$ the strict transform of $C$ by the morphism $\pi_E$. The points $\pi_E^* C \cap E_{red}$ are called {\em infinitely near points\/} of $C$ in $E$.

\begin{Remark}
If $C$ is a curve with only non-singular irreducible components and $E$ is an irreducible component of $\pi_{C}^{-1}(0)$, then the number of infinitely near points of $\pi_E^*C$ in $E_{red}$ is equal to $b_E$. That is, the cardinal of the set $\pi_E^*C \cap E_{red}$ coincides with $b_E$.
\end{Remark}

%%%%%%%%%%%%%%%%%%%%%%%%%%%%%%%%%%%%%%%%%%%%%%%%%%%%%%%%%%%%%
\section{Logarithmic foliations}\label{sec:logarithmic}
Consider a germ of plane curve $C=\cup_{i=1}^r C_i$ in $({\mathbb C}^2,0)$.  Take $f \in {\mathbb C}\{x,y\}$ such that $C=(f=0)$ and let us write $f= f_1 \cdots f_r$ with $f_i \in {\mathbb C}\{x,y\}$ irreducible. Given $\lambda= (\lambda_1, \ldots, \lambda_r) \in {\mathbb C}^r$, we can consider the logarithmic foliation ${\mathcal L}_{\lambda}^C$ defined by
\begin{equation}\label{eq:1-forma-log}
  f_1 \cdots f_r \sum_{i=1}^{r} \lambda_i \frac{df_i}{f_i}=0.
\end{equation}
The logarithmic foliation ${\mathcal L}_{\lambda}^C$ belongs to  ${\mathbb G}_C$ provided that $\lambda$ avoids certain rational resonances. Each generalized curve foliation ${\mathcal F} \in {\mathbb G}_C$ has a {\em logarithmic model\/} ${\mathcal L}_\lambda^C$, that is, a logarithmic foliation such that the Camacho-Sad indices of ${\mathcal F}$ and ${\mathcal L}_\lambda^C$ coincide along the reduction of singularities (see \cite{Cor-2003}); note that $\mathcal F$ and ${\mathcal L}_\lambda^C$ have the same separatrices and the same minimal  reduction of singularities. Moreover, the logarithmic model of $\mathcal F$ is unique once a reduced equation of the separatrices is fixed. Thus, for each foliation ${\mathcal F} \in {\mathbb G}_C$, we denote by $\lambda({\mathcal F})$ the exponent vector of the logarithmic model of $\mathcal F$ and we denote ${\mathbb G}_{C,\lambda}$ the set of foliations ${\mathcal F} \in {\mathbb G}_C$ such that $\lambda({\mathcal F})=\lambda$.
A particular case of logarithmic foliation is the ``hamiltonian" foliation defined by $df=0$ which corresponds to $\lambda=(1,\ldots,1)$; this foliation coincides with the foliation
${\mathcal G}_f$ used in section~\ref{sec:local-invariants}.

Let us fix some notations concerning logarithmic foliations that will be used in the sequel.  Assume that the curve $C=\cup_{i=1}^r C_i$   has only non-singular irreducible components and consider   a non-dicritical logarithmic foliation  ${\mathcal L}_\lambda^C$ given by~\eqref{eq:1-forma-log}. Let $\pi_C : X_C \to ({\mathbb C}^2,0)$ be the minimal reduction of singularities of $C$, take $E$  an irreducible component of $\pi_C^{-1}(0)$ and consider $\pi_E: X_E \to ({\mathbb C}^2,0)$ the morphism reduction of $\pi_C$ to $E$ (see section~\ref{Apendice:grafo-dual}).
Given an irreducible component $C_j$ of $C$ and a divisor $F$ of $\pi_C^{-1}(0)$, we denote $\varepsilon_{F}^{C_j}=1$ if the geodesic of $C_j$ contains the divisor $F$ and  $\varepsilon_{F}^{C_j}=0$ otherwise, that is,
$$
\varepsilon_{F}^{C_j}=
\left\{
  \begin{array}{ll}
    1, & \hbox{ if } j \in I_F \\
    0, & \hbox{ if } j \not \in I_F.
  \end{array}
\right.
$$
The {\em residue of the logarithmic foliation} ${\mathcal L}_\lambda^C$ {\em along the divisor} $E$ is given by
\begin{equation}\label{eq:peso-divisor-log}
\kappa_E({\mathcal L}_\lambda^C)=\sum_{j=1}^{r} \lambda_j \sum_{E' \leq E} \varepsilon_{E'}^{C_j},
\end{equation}
where $E' \leq E$ means  all divisors in $G(C)$ which are in the geodesic of $E$ (including $E$ itself).
Note that $\kappa_E({\mathcal L}_\lambda^C)= \sum_{j=1}^{r} \lambda_j m_E^{C_j}$ where $m_E^{C_j}$ is the multiplicity of $f_j \circ \pi_E$ along the divisor $E$ (see \cite{Pau-89,Pau-95}).

Let $\{R_1^E, R_2^E, \cdots, R_{b_E}^E\}$ be the set of points  $\pi_E^* C \cap E_{red}$ where we denote $b_E=b_E^C$ and
put $I_{R_l^E}^C=\{ i \in \{1,\ldots, r\} \ : \ \pi_E^*C_i \cap E_{red}= \{R_l^E\}\}$ for $l=1,2,\ldots, b_E$, that is, $i \in I_{R_l^E}^C$ if $E$ belongs to the geodesic of the curve $C_i$ in $G(C)$ and $R_i^E$ is an infinitely near point of $C_i$. With the notations introduced in subsection~\ref{Apendice:grafo-dual}, we have that
$I_E= \cup_{l=1}^{b_E} I_{R_l^E}^C$.

 An easy computation shows that the Camacho-Sad index of $\pi_E^* {\mathcal L}_\lambda^C$ relative to $E_{red}$ at a point $R_l^E$ is given by
\begin{equation}\label{eq:camacho-sad-log}
\mathcal{I}_{R_l^E}(\pi_E^* {\mathcal L}_\lambda^C,E_{red})= -\frac{\ \  \sum_{i \in I_{R_l^E}^C} \lambda_i \ \ }{\kappa_E({\mathcal L}_\lambda^C)}
\end{equation}
(see \cite{Cor-2003}, section 4).

In appendix~\ref{ap:ramificacion} we include a subsection where we explain the behaviour of the above invariants associated to logarithmic foliations after ramification (see subsection~\ref{subsec:ram-log}).

\section{Collinear and non-collinear points and divisors}\label{sec:collinear}

In this section we will introduce some notations and definitions in order to describe properties of the jacobian curve that will be given in section~\ref{sec:jacobian}.

Let $C$ and $D$ be two plane curves  in $({\mathbb C}^2,0)$ without common branches and assume that the curve $Z=C \cup D$ has only non-singular irreducible components.
Take now ${\mathcal F} \in {\mathbb G}_C$ and ${\mathcal G} \in {\mathbb G}_D$ and let % be two singular foliations in $({\mathbb C}^2,0)$
 $\pi_Z: X_Z \to ({\mathbb C}^2,0)$ be the minimal  reduction of singularities of $Z$ which gives a common reduction of singularities of $\mathcal F$ and $\mathcal G$.
Recall that, given   an irreducible component $E$ of $\pi_Z^{-1}(0)$, we denote by $\pi_E: X_E \to ({\mathbb C}^2,0)$ the {\em morphism reduction\/} of $\pi_Z$ to $E$ (see section~\ref{Apendice:grafo-dual}).

\begin{Remark}\label{rk:coordenadas-adaptadas}
With the above assumptions about $Z$, if $v(E)=p$ (see appendix~\ref{ap:equisingularidad}), then the morphism $\pi_E$ is a composition of $p$ punctual blow-ups
$$
({\mathbb C}^2,0) \overset{\sigma_1}{\longleftarrow} (X_{1},P_1)  \overset{\sigma_2}{\longleftarrow}  \cdots \overset{\sigma_{p-1}}{\longleftarrow} (X_{p-1},P_{p-1})  \overset{\sigma_{p}}{\longleftarrow} X_p=X_{E}.
$$
Moreover, if $(x,y)$ are coordinates in $({\mathbb C}^2,0)$ there is a change
of coordinates $(x,y)=(\bar{x}, \bar{y}+
\varepsilon_E(\bar{x}))$, with $\varepsilon_E(x)= a_1 x + \cdots +
a_{p-1} x^{p-1}$,  such that the blow-up $\sigma_j$ is given by
$x_{j-1}=x_j$, $y_{j-1}=x_j y_j$, for $j=1,2, \ldots, p$, where
$(x_j,y_j)$ are coordinates centered at $P_j$ and
$(x_0,y_0)=(\bar{x},\bar{y})$. We say that
$(\bar{x},\bar{y})$ are {\em coordinates in\/} $({\mathbb
C}^2,0)$ {\em adapted to\/} $E$.
Note that in these coordinates, if $(x_p,y_p)$ are coordinates in the first chart of $E_{red}$ we have that $\pi_E(x_p,y_p)=(x_p,x_p^p y_p)$ and $E_{red}=(x_p=0)$.
\end{Remark}

In this section, we will denote $b_E=b_E^Z$.
Let $\{R_1^E, R_2^E, \ldots, R_{b_E}^E\}$  be  the infinitely near points of $Z$ in $E$, that is, $\pi_E^* Z \cap E_{red} = \{R_1^E, R_2^E, \ldots, R_{b_E}^E\}$. Note that these points are the union of the singular points of $\pi_E^* {\mathcal F}$ and $\pi_E^* {\mathcal G}$ in the first chart of  $E_{red}$ (the singular points of the foliations which do not correspond to a corner of the divisor).
%, that is, $\pi_E^* Z \cap E_{red} = \{R_1^E, R_2^E, \ldots, R_{b_E}^E\}$ are the infinitely near points of $Z$ in $E$.
We denote
$$\Delta_E^{{\mathcal F},{\mathcal G}}(R_i^E)= \mathcal{I}_{R_i^E}(\pi_E^* {\mathcal G},E_{red}) - \mathcal{I}_{R_i^E}(\pi_E^* {\mathcal F},E_{red}) $$
where $\mathcal{I}_{R_i^E}(\pi_E^* {\mathcal F},E_{red})$ is the Camacho-Sad index of $\pi_E^* {\mathcal F}$ relative to $E_{red}$ at the point $R_i^E$ (see definition given in  \eqref{eq:indice-C-S}). We will denote $\Delta_E(R_i^E)=\Delta_E^{{\mathcal F},{\mathcal G}}(R_i^E)$ if it is clear the foliations $\mathcal F$ and $\mathcal G$ we are working with.

In view of  the notations given in \cite{Kuo-P-2002,Kuo-P-2004} for curves, we introduce the following definitions for foliations:

\begin{Definition}\label{def:collinear}
We say that an infinitely near point $R_l^E$ of $Z$ is a {\em collinear point\/} for the foliations $\mathcal F$ and $\mathcal G$ in $E$ if $\Delta_E(R_l^E)=0$; otherwise we say that $R_l^E$ is a {\em non-collinear point\/}.

We say that a divisor $E$ is {\em collinear (for the foliations $\mathcal F$ and $\mathcal G$)\/} if $\Delta_E(R_l^E)=0$ for all $l=1,\ldots, b_E$; otherwise $E$ is called a {\em non-collinear divisor\/}. A divisor $E$ is called {\em purely non-collinear\/} if $\Delta_E(R_l^E) \neq 0$  for each $l=1,\ldots, b_E$.
\end{Definition}
We denote by $\text{Col}(E)$ the set of collinear points of $E$ and by $\operatorname{NCol}(E)$ the set of non-collinear points (for the foliations $\mathcal{F}$ and $\mathcal{G}$). It is clear that $\text{Col}(E) \cup \operatorname{NCol}(E) = \{R_1^E, R_2^E, \ldots, R_{b_E}^E\}$.
\begin{Remark}\label{rk:max-div-non-colinear}
Note that if $E$ is a maximal bifurcation divisor (with the partial order in $G(C)$ given in subsection~\ref{Apendice:grafo-dual}), then $E$ is purely non-collinear. This follows from the fact that, if $E$ is a maximal bifurcation divisor, then each infinitely near point $R_l^E$ of $Z$ in $E$ is in the geodesic of only one irreducible component of $Z$ and hence it is a  singular point for only one of the foliations $\pi_E^* \mathcal{F}$ or $\pi_E^* \mathcal{G}$. Moreover, it is a simple singularity. In fact, we have that
$$\Delta_E(R_l^E)=\left\{
                    \begin{array}{ll}
                      -\mathcal{I}_{R_l^E}(\pi_E^* \mathcal{F},E_{red}), & \hbox{ if } R_l^E \in \pi_E^* C \cap E \\
                      \mathcal{I}_{R_l^E}(\pi_E^* \mathcal{G},E_{red}) , & \hbox{ if } R_l^E \in \pi_E^* D \cap E
                    \end{array}
                  \right.
$$
If  $R_l^E \in \pi_E^* C \cap E$, we have that $\mathcal{I}_{R_l^E}(\pi_E^* \mathcal{F},E_{red}) \neq 0$ since $\mathcal{F}$ is a generalized curve foliation and $R_l^E$ is a simple singularity of $\pi_E^* \mathcal{F}$ (similarly when $R_l^E \in \pi_E^* D \cap E$).
Consequently,  $\Delta_E(R_l^E) \neq 0$  for each $l=1,\ldots, b_E$.
\end{Remark}

Although the definition of $\Delta_E(R_l^E)$ seems different to the one given by Kuo and Parusi\'nski in \cite{Kuo-P-2002,Kuo-P-2004}, we will show in subsection~\ref{sec:curvas} that both definitions coincide in the case of curves.

\medskip
Take coordinates $(x_p,y_p)$ in the first chart of $E_{red}$ such that $\pi_E(x_p,y_p)=(x_p,x_p^p y_p)$, $E_{red}=(x_p=0)$ and assume that $R_{l}^E =(0,c_l^E)$, $l=1,2,\ldots, b_E$, in these coordinates.
We define the {\em rational function ${\mathcal M}_E (z)={\mathcal M}^{{\mathcal F},{\mathcal G}}_E(z)$ associated to the divisor\/} $E$ {\em for the foliations $\mathcal F$ and $\mathcal G$} by
\begin{equation}\label{eq:funcion-meromorfa}
{\mathcal M}_E (z)=\sum_{l=1}^{b_E} \frac{\Delta_E(R_l^E)}{z-c_l^E}.
\end{equation}
\begin{Remark}
Note that although  $\Delta_E^{\mathcal{F},\mathcal{G}}(R_l^E)$ and ${\mathcal M}_E^{\mathcal{F},\mathcal{G}}(z)$ depend on the foliations $\mathcal F$ and $\mathcal G$, we have that
$$
\Delta_E^{\mathcal{F},\mathcal{G}}(R_l^E) = \Delta_E^{\mathcal{L}_\lambda^C,\mathcal{L}_{\mu}^D}(R_l^E); \qquad {\mathcal M}^{{\mathcal F},{\mathcal G}}_E(z)= {\mathcal M}^{\mathcal{L}_\lambda^C,\mathcal{L}_{\mu}^D}_E(z)
$$
provided that ${\mathcal F} \in {\mathbb G}_{C,\lambda}$ and ${\mathcal G} \in {\mathbb G}_{D,\mu}$.
\end{Remark}
\begin{Remark}
Observe that if $E$ is a non-collinear divisor, then ${\mathcal M}_E(z) \not \equiv 0$.
\end{Remark}
Let $M(E)=\{Q_1^E, \ldots, Q_{s_E}^E\}$ be the set of points of $E_{red}$ given by $Q_l^{E}=(0,q_l)$ in coordinates $(x_p,y_p)$ where $\{q_1, \ldots, q_{s_E}\}$ is the set of zeros of ${\mathcal M}_E (z)$. We denote by $t_{Q_l^{E}}$ the multiplicity of $q_l$ as a zero of ${\mathcal M}_E (z)$ and $t(E)=\sum_{l=1}^{s_E} t_{Q_l^{E}}$ the degree of the numerator of the rational function $\mathcal{M}_E(z)$. We put $t_P=0$ for any $P \in E \smallsetminus M(E)$. Note that it can happen that $M(E)=\emptyset$ (see example in \cite{Kuo-P-2004}, p. 584).
\begin{Lemma}\label{lema:E-no-colineal}
If $\operatorname{NCol}(E) \neq \emptyset$, that is, $E$ is a non-collinear divisor, then we have that
%$$
\begin{equation}\label{eq:no-colineal}
\operatorname{NCol}(E) \cap M(E) = \emptyset \quad \text{ and } \quad \sharp \operatorname{NCol}(E) \geq 1 +\sum_{P \in M(E)} t_P = 1+ t(E).
\end{equation}
Moreover, if $\sum_{R_l^E \in \operatorname{NCol}(E)} \Delta_{E}(R_l^E) \neq 0$, then we have that $$\sharp \operatorname{NCol}(E) =1 +\sum_{P \in M(E)} t_P.$$
\end{Lemma}
\begin{proof}
With the notations above, we can write ${\mathcal M}_E (z)$ as follows
$$
{\mathcal M}_E (z)=\sum_{R_l^E \in \operatorname{NCol}(E)} \frac{\Delta_E(R_l^E)}{z-c_l^E}.
$$
Thus, the set of zeros of ${\mathcal M}_E(z)$ is given by the roots of the polynomial
\begin{equation}\label{eq:ceros-ME}
\sum_{R_l^E \in \operatorname{NCol}(E)} \Delta_{E}(R_l^E) \prod_{\substack{j \text{ with }  j \neq l  \\  R_j^E \in \operatorname{NCol}(E) }}  (z-c_j^E).
\end{equation}
Consequently, if $z=c_{l_0}^E$ is a zero of ${\mathcal M}_E(z)$, then $\Delta_{E}(R_{l_0}^E) \prod_{\substack{R_j^E  \in \operatorname{NCol}(E) \\ j \neq l_0}} (c_{l_0}^E-c_j^E)=0$ which implies that $\Delta_E(R_{l_0}^E)=0$ and hence $R_{l_0}^E \not \in \operatorname{NCol}(E)$. Moreover, the degree of the polynomial given in equation~\eqref{eq:ceros-ME} is $\leq \sharp \operatorname{NCol}(E) -1$; the equality is attained when $\sum_{R_l^E \in \operatorname{NCol}(E)} \Delta_{E}(R_l^E) \neq 0$.   Thus we have the statements of the lemma.

\end{proof}
\begin{Remark}
Observe that we can have that $\sum_{R_l^E \in \operatorname{NCol}(E)} \Delta_{E}(R_l^E) = 0$ even if $E$ is a purely non-collinear divisor. This can happen for instance when $E=E^1$ is a bifurcation divisor since in this situation $\sum_{i=1}^{b_E} \mathcal{I}_{R_i^E} (\pi_E^* \mathcal{F},E_{red})=\sum_{i=1}^{b_E} \mathcal{I}_{R_i^E} (\pi_E^* \mathcal{G},E_{red})=-1$ and this implies
 $\sum_{R_l^E \in \operatorname{NCol}(E)} \Delta_{E}(R_l^E) = 0$ (see also Remark~\ref{rmk:suma-delta} and Corollary~\ref{cor:colinear}).
\end{Remark}
\begin{Remark}\label{rmk:CE-ME}
  Note that it can happen that $\operatorname{Col}(E) \cap M(E) \neq \emptyset$. With the notations of section~\ref{sec:logarithmic}, consider the foliations ${\mathcal F}={\mathcal L}^C_{\lambda}$ and ${\mathcal G}={\mathcal L}^D_{\mu}$ where
 $$ C=(f=0), \quad f(x,y)=(y-x) (y+x^2) (y-x^2) (y+2x^2), \quad \lambda=(1,1,1,3)$$
 $$  D=(g=0), \quad g(x,y)=(y+x) (y+x^2+x^3) (y-x^2+x^3), \quad \mu=(3,3,1).$$
Take $Z=C \cup D$. Consider the morphism $\sigma= \pi_1 \circ \pi_2$ where $\pi_1: X_1 \to ({\mathbb C}^2,0)$ is the blow-up of the origin, $E_1= \pi_1^{-1}(0)$ and $\pi_2: X_2 \to (X_1,P_1)$ is the blow-up of the origin $P_1$ of the first chart of $E_1$ and $E_2=\pi_2^{-1}(P_1)$. Taking coordinates $(x_2,y_2)$ in the first chart of $E_2$ such that $\sigma(x_2,y_2)=(x_2,x_2^2y_2)$ we have that $\sigma^*Z \cap E_2 =\{R_1^{E_2},R_{2}^{E_2},R_{3}^{E_2}\}$ where $R_1^{E_2}=(0,-1)$, $R_{2}^{E_2}=(0,1)$ and $R_{3}^{E_2}=(0,-2)$.
A simple computation shows that
$$\Delta_{E}(R_1^{E_2})= -\frac{2}{11}; \quad \Delta_{E}(R_2^{E_2})= 0;\quad  \Delta_{E}(R_3^{E_2})= \frac{3}{11}.$$
Thus $\operatorname{Col}(E_2)=\{R_2^{E_2}\}$ and $\operatorname{NCol}(E_2)=\{R_1^{E_2}, R_3^{E_2}\}$. Moreover, we have that
$${\mathcal M}_{E_2}(z)= -\frac{2}{11} \frac{1}{(z+1)} + \frac{3}{11} \frac{1}{(z+2)}  = \frac{z-1}{11(z+2)(z+1)}$$
which implies that
$M(E_2)=\{R_2^{E_2}\}$.
\end{Remark}

\medskip
Given a non-collinear divisor $E$ and a point $P \in E_{red}$, we define

\begin{equation}\label{eq:tau_E}
\tau_E(P)=\left\{
              \begin{array}{ll}
                t_P, & \hbox{ if } P \in M(E) \\
                - 1, & \hbox{ if } P \in \operatorname{NCol}(E) \\
                0, & \hbox{ otherwise.}
              \end{array}
            \right.
\end{equation}
Note that $\sum_{P \in E_{red}} \tau_E(P)=t(E) - \sharp \operatorname{NCol}(E)$ which is a negative integer (it is the degree of the rational function $\mathcal{M}_E(z)$).

\subsection{Collinear and non-collinear infinitely near points} Let us explain now  the behaviour of collinear (resp. non-collinear) infinitely near points by blowing-up. Recall that here we denote $b_E=b_E^Z$ for any divisor $E$ in $G(Z)$.

\begin{Lemma}\label{lema:explosio-colineal}
Let $E$ and $E'$ be two consecutive divisors in $G(Z)$ with $E <E'$ and $b_{E'}=1$. We can write  $\pi_{E'}=  \pi_{E} \circ \sigma $ where   $\sigma: X_{E'} \to (X_E,P)$ is the blow-up with center at a point $P \in E_{red}$. Let $Q$ be the point $\pi_{E'}^* Z \cap E_{red}'$.

If $P$ is a collinear point (resp. a non-collinear point) for the foliations $\mathcal F$ and $\mathcal G$ in $E$, then $Q$ is a collinear point (resp. non-collinear point) for $E'$.
\end{Lemma}
\begin{proof}
Let us denote by $\widetilde{E}_{red}$ the strict transform of $E_{red}$ by $\sigma$ and $\widetilde{P}=\widetilde{E}_{red} \cap E_{red}'$. Given any singular foliation $\mathcal F$,
the Camacho-Sad indices  verify the following equalities (see \cite{Cam-S}) %({\color{red} notaciones para los transformados de los puntos y divisores})
$$
\mathcal{I}_{\widetilde{P}}(\pi^*_{E'} {\mathcal F},\widetilde{E}_{red})  = \mathcal{I}_P(\pi_{E}^*{\mathcal F},E_{red})-1 $$
$$ \mathcal{I}_{\widetilde{P}}(\pi^*_{E'}{\mathcal F},E'_{red})+\mathcal{I}_Q(\pi^*_{E'}{\mathcal F},E'_{red})=-1
$$
Moreover, since $\mathcal F$ is a generalized curve foliation, then $\widetilde{P}$ is a simple singularity for $\pi_{E'}^* {\mathcal F}$ and hence we have that
$$
 \mathcal{I}_{\widetilde{P}}(\pi^*_{E'} {\mathcal F},\widetilde{E}_{red}) \cdot \mathcal{I}_{\widetilde{P}}(\pi^*_{E'} {\mathcal F},E'_{red})=1.
$$
Then the index $\mathcal{I}_Q(\pi^*_{E'}{\mathcal F},E'_{red})$ can be computed as
$$\mathcal{I}_Q(\pi^*_{E'}{\mathcal F},E'_{red})= -1 -\frac{1}{\mathcal{I}_P(\pi_{E}^*{\mathcal F},E_{red})-1}= - \frac{\mathcal{I}_P(\pi_{E}^*{\mathcal F},E_{red})}{\mathcal{I}_P(\pi_{E}^*{\mathcal F},E_{red})-1}.$$
Thus, using the expression above for the foliations $\mathcal F$ and $\mathcal G$,  we have that
$$
  \Delta_{E'}(Q)  =\left| \begin{array}{cc}
  1 & \mathcal{I}_Q(\pi^*_{E'} {\mathcal F},E'_{red}) \\
  1 &  \mathcal{I}_Q(\pi^*_{E'} {\mathcal G},E'_{red})
  \end{array} \right|
     =\frac{\Delta_E(P)}{(\mathcal{I}_P(\pi_{E}^*{\mathcal F},E_{red})-1)(\mathcal{I}_P(\pi_{E}^*{\mathcal G},E_{red})-1)}
$$
and then $\Delta_{E'}(Q)=0$ if and only if $\Delta_E(P)=0$. This gives the result.
\end{proof}

Consider now $E$ and $E'$  two consecutive bifurcation divisors in $G(Z)$, that is, there is a chain of consecutive divisors
$$E_0=E < E_1 < \cdots < E_{k-1} < E_k=E'$$
with $b_{E_l}=1$ for $l=1, \ldots, k-1$ and the morphism $\pi_{E'}=  \pi_E \circ \sigma $ where $\sigma: X_{E'} \to (X_E,P)$ is a composition of $k$ punctual blow-ups
\begin{equation} \label{eq:bifurcacion-consecutivos}
(X_E, P) \overset{\sigma_1}{\longleftarrow} (X_{E_1},P_1)  \overset{\sigma_2}{\longleftarrow}  \cdots \overset{\sigma_{k-1}}{\longleftarrow} (X_{E_{k-1}},P_{k-1})  \overset{\sigma_{k}}{\longleftarrow} X_{E'}.
\end{equation}
If $E$ and $E'$ are two consecutive bifurcation divisors as above, we say that $E'$ {\em arises from $E$ at $P$} and we denote $E <_P E'$.

As a consequence of Lemma~\ref{lema:explosio-colineal}, we have that if $P$ is a collinear point (resp. a non-collinear point) for $\mathcal F$ and $\mathcal G$ relative to $E$, then $P_l$ is  a collinear point (resp. non-collinear point) relative to $E_l$ for $l=1,\ldots, k-1$. Moreover, we have that

\begin{Corollary}
Let $E$ be the first bifurcation divisor in $G(Z)$. We can write $\pi_E = \sigma_1 \circ \sigma _2 \circ \cdots \circ \sigma_k$ as a composition of $k$ punctual blow-ups
$$
({\mathbb C}^2, 0) \overset{\sigma_1}{\longleftarrow} (X_1,P_1)  \overset{\sigma_2}{\longleftarrow}  \cdots \overset{\sigma_{k-1}}{\longleftarrow} (X_{k-1},P_{k-1})  \overset{\sigma_{k}}{\longleftarrow} X_k=X_{E}.
$$
We denote $E_i=\sigma_i^{-1}(P_{i-1})$ with $P_0=0$.
Then all the divisors $E_i$, $1 \leq i \leq k-1$, are collinear.
\end{Corollary}
\begin{proof}
Since $b_{E_l}=1$ for $l=1,\ldots, k-1$, it is enough to prove that $P_1$ is a collinear point for $\mathcal F$ and $\mathcal G$ relative to $E_1=E^1$. But this is a consequence of the fact that $I_{P_1} (\sigma_1^{\ast}{\mathcal F},E_{1,red})=I_{P_1} (\sigma_1^{\ast}{\mathcal G},E_{1,red})=-1$.
Thus the result follows straightforward.
\end{proof}
\begin{Remark}\label{rmk:suma-delta}
Note that, if $E$ is the first bifurcation divisor, the properties of the Camacho-Sad index imply that
$$ \sum_{l=1}^{b_E} \Delta_E(R_l^E)=0$$
where $\pi_E^* Z \cap E_{red}= \{ R_1^E, \ldots, R_{b_E}^E\}$.
\end{Remark}
The above equality also holds in the following context:
\begin{Corollary}\label{cor:colinear}
Let $E$ and $E'$ be two consecutive bifurcation divisors in $G(Z)$ such that $E'$ arises from $E$ at $P$. If $P$ is a collinear point, then
$$\sum_{l=1}^{b_{E'}} \Delta_{E'}(R_{l}^{E'})=0$$
where $\pi_{E'}^* Z \cap E_{red}'=\{R_1^{E'},\ldots, R_{b_{E'}}^{E'}\}$.
\end{Corollary}
\begin{proof}
As we have explained before we have that $\pi_{E'}=\pi_E \circ \sigma$ where $\sigma: X_{E'} \to (X_E,P)$ is a composition of $k$ punctual blow-ups

\begin{equation*}
(X_E, P) \overset{\sigma_1}{\longleftarrow} (X_{E_1},P_1)  \overset{\sigma_2}{\longleftarrow}  \cdots \overset{\sigma_{k-1}}{\longleftarrow} (X_{E_{k-1}},P_{k-1})  \overset{\sigma_{k}}{\longleftarrow} X_{E'}.
\end{equation*}
We denote $E_i=\sigma_{i}^{-1}(P_{i-1})$ with $P_0=P$ and we have that $b_{E_i}=1$ for $i=1, \ldots, k-1$. Since $P$ is a collinear point, then $\mathcal{I}_P(\pi_E^* \mathcal{F},E_{red})= \mathcal{I}_P(\pi_E^* \mathcal{G},E_{red})$ and the properties of the Camacho-Sad indices imply that $\mathcal{I}_{P_i}(\pi_{E_i}^* \mathcal{F},F_{i,red})= \mathcal{I}_{P_i}(\pi_{E_i}^* \mathcal{G},E_{i,red})$ for $i=1, \ldots, k-1$. Consequently, we have that
$$
\sum_{l=1}^{b_{E'}} \mathcal{I}_{R_l^{E'}} (\pi_{E'}^* \mathcal{F},E_{red}') = \sum_{l=1}^{b_{E'}} \mathcal{I}_{R_l^{E'}} (\pi_{E'}^* \mathcal{G},E_{red}')
$$
which is equivalent to $\sum_{l=1}^{b_{E'}} \Delta_{E'}(R_{l}^{E'})=0$.
\end{proof}
\medskip
\subsection{Weighted initial forms and non-collinear divisors}
Let us introduce the following notation in order to describe the relationship between the Newton polygon and the infinitely near points of a curve (see subsection~\ref{sec:parte-inicial} and also \cite{Cor-2009}). From now on we will always assume that we choose coordinates $(x,y)$ such that $x=0$ is not tangent to the curve $Z=C \cup D$ union of the separatrices of $\mathcal F$ and $\mathcal G$. This will imply that the first side of the Newton polygons ${\mathcal N}({\mathcal F};x,y)$ and ${\mathcal N}({\mathcal G};x,y)$ has slope greater or equal to $-1$.

Assume that $Z$ is a curve with only non-singular irreducible components and consider $\pi_Z:X_Z \to ({\mathbb C}^2,0)$ its minimal reduction of singularities. Take any divisor $E$ of $\pi^{-1}_Z(0)$ with $v(E)=p$ and consider $\pi_E:X_E \to ({\mathbb C}^2,0)$  the morphism reduction of $\pi_Z$ to $E$. With the notations introduced in Remark~\ref{rk:coordenadas-adaptadas}, if $(x,y)$ are coordinates adapted to $E$, the points $\pi_E^* Z \cap E_{red}$ are determined by $\text{In}_{p} (h;x,y)$ where $h=0$ is a reduced equation of the curve $Z$. More precisely, if we take $(x_p,y_p)$ coordinates in the first chart of $E_{red}$ with $\pi_E(x_p,y_p)=(x_p,x_p^py_p)$ and $E_{red}=(x_p=0)$, thus the points of $\pi_E^*Z \cap E_{red}$ are given by $x_p=0$ and $\sum_{i+pj=k}h_{ij} y_p^j=0$ where $\text{In}_{p} (h;x,y)=\sum_{i+pj=k}h_{ij} x^i y^j$.

We are interested in determine the points $\pi_E^*{\mathcal J}_{\mathcal{F},\mathcal{G}} \cap E_{red}$, thus if $v(E)=p$ and $(x,y)$ are coordinates adapted to $E$, we would like to determine ${\rm In}_p (J;x,y)$ where $J(x,y)=0$ is an equation of the jacobian curve. Next result proves that the initial form  ${\rm In}_p (J;x,y)$    is determined by the initial forms  ${\rm In}_p(\omega), {\rm In}_p(\eta)$ of the 1-forms defining the foliations $\mathcal F$ and $\mathcal G$ provided that the divisor $E$ is non-collinear.

\begin{Lemma}\label{lema:parte-inicial-jac}
  Let $E$ be an irreducible component of $\pi_{Z}^{-1}(0)$ with $v(E)=p$ and take $(x,y)$ coordinates adapted to $E$. If $E$ is a non-collinear divisor, then
$${\rm In}_p(\omega) \wedge {\rm In}_p(\eta) \not \equiv 0,$$
where ${\rm In}_p(\omega)= {\rm In}_p(\omega;x,y)$ and ${\rm In}_p(\eta)={\rm In}_p(\eta;x,y)$,
and hence
$${\rm In}_p (J;x,y) = J_p(x,y)$$
with ${\rm In}_p(\omega) \wedge {\rm In}_p(\eta) = J_p (x,y) dx \wedge dy$.
\end{Lemma}
\begin{proof}
Take an irreducible component $E$ of $\pi_{Z}^{-1}(0)$ with $v(E)=p$ and let $(x,y)$ be coordinates adapted to $E$. Assume that ${\rm In}_p(\omega) \wedge {\rm In}_p(\eta) \equiv 0$, that is, if we write
\begin{align*}
  {\rm In}_p(\omega) & = A_I(x,y) dx + B_I (x,y) dy \\
  {\rm In}_p(\eta) & = P_I(x,y) dx + Q_I (x,y) dy
\end{align*}
then
\begin{equation}\label{eq:parte-inicial}
A_I(x,y) Q_I(x,y) - B_I(x,y) P_I(x,y)\equiv 0.
\end{equation}
Note that $(A_I,B_I)\not \equiv (0,0)$ and $(P_I,Q_I) \not \equiv (0,0)$.

Take coordinates $(x_p,y_p)$ in the first chart of $E_{red}$ such that $\pi_E(x_p,y_p)=(x_p,x_p^p y_p)$, $E_{red}=(x_p=0)$ and assume that $R_{l}^E =(0,c_l^E)$, for $l=1,\ldots, b_E$, in these coordinates where $\pi_E^*Z \cap E_{red} = \{ R_1^E, \ldots, R_{b_E}^E\}$. Let $\omega^E$ and $\eta^E$ be the strict transforms of $\omega$ and $\eta$ by $\pi_E$ with
\begin{align*}
\omega^E &= A^E(x_p,y_p) dx_p + x_p B^E(x_p,y_p) dy_p,  \\
\eta^E &= P^E(x_p,y_p) dx_p + x_p Q^E(x_p,y_p) dy_p.
\end{align*}
From the definition of the Camacho-Sad index we have that
\begin{align*}
  \mathcal{I}_{R_l^E}(\pi_E^* {\mathcal F},E_{red}) &  = - \text{Res}_{y=c_l^E} \frac{B^E(0,y)}{A^E(0,y)}, \\
  \mathcal{I}_{R_l^E}(\pi_E^* {\mathcal G},E_{red}) &  = - \text{Res}_{y=c_l^E} \frac{Q^E(0,y)}{P^E(0,y)}.
\end{align*}
Note that $A^E(0,y)=A_I(1,y) + py B_I(1,y)$, $B^E(0,y)=B_I(1,y)$, $P^E(0,y)=P_I(1,y) + p yQ_I(1,y)$ and $Q^E(0,y)=Q_I(1,y)$. Thus, the equality given in~\eqref{eq:parte-inicial} implies
$$
\frac{B^E(0,y)}{A^E(0,y)}=\frac{Q^E(0,y)}{P^E(0,y)}
$$
and hence $\mathcal{I}_{R_l^E}(\pi_E^* {\mathcal F},E_{red})= \mathcal{I}_{R_l^E}(\pi_E^* {\mathcal G},E_{red})$ for $l=1,\ldots, b_E$ in contradiction with the fact that the divisor $E$ is non-collinear.
\end{proof}
Observe that the result above is true for the first divisor $E^1$ although the curves $C$ and $D$ have singular irreducible components, and hence we have that Lemma~\ref{lema:multiplicidad-E1} can be obtained as a consequence of the previous result since the conditions over the foliations $\mathcal{F}$ and $\mathcal{G}$ in Lemma~\ref{lema:multiplicidad-E1} imply that $E^1$ is a non-collinear divisor.
\medskip
Moreover, next result shows that given ${\mathcal F}, \widetilde{\mathcal F} \in {\mathbb G}_{C,\lambda}$ and ${\mathcal G},\widetilde{\mathcal G} \in {\mathbb G}_{D,\mu}$, we have that
$$\text{In}_p(J_{{\mathcal F},\mathcal{G}};x,y)= \text{In}_p(J_{\widetilde{\mathcal F},\widetilde{\mathcal{G}}};x,y)$$
provided that $E$ is a non-collinear divisor with $v(E)=p$, where $(x,y)$ are coordinates adapted to $E$ and $J_{{\mathcal F},\mathcal{G}}(x,y)=0$ and $J_{\widetilde{\mathcal F},\widetilde{\mathcal{G}}}(x,y)=0$ are equations of the jacobian curves ${\mathcal J}_{\mathcal{F},\mathcal{G}}$ and ${\mathcal J}_{\widetilde{\mathcal{F}},\widetilde{\mathcal{G}}}$ respectively. Given a foliation $\mathcal F$, we will denote by $\omega_{\mathcal F}$ a $1$-form defining $\mathcal F$ and $\text{In}_p(\omega_{\mathcal F}) =\text{In}_p(\omega_{\mathcal F};x,y)$. Thus we have the following result

\begin{Lemma}\label{lema:parte-incial-modelo-log}
Let $E$ be an irreducible component of $\pi^{-1}(0)$ with $v(E)=p$ and assume that $(x,y)$ are coordinates adapted to $E$. Consider foliations  ${\mathcal F}, \widetilde{\mathcal F} \in {\mathbb G}_{C,\lambda}$ and ${\mathcal G},\widetilde{\mathcal G} \in {\mathbb G}_{D,\mu}$ then
$${\rm In}_p(\omega_{\mathcal F})= {\rm In}_p(\omega_{\widetilde{{\mathcal F}}}); \qquad {\rm In}_p(\omega_{\mathcal G})= {\rm In}_p(\omega_{\widetilde{{\mathcal G}}}).$$
Hence, if $E$ is a non-collinear divisor, we have that
$${\rm In}_p(J_{{\mathcal F},\mathcal{G}};x,y)= {\rm In}_p(J_{\widetilde{\mathcal F},\widetilde{\mathcal{G}}};x,y).$$
\end{Lemma}
\begin{proof}
Let us prove that ${\rm In}_p(\omega_{\mathcal F})= {\rm In}_p(\omega_{\widetilde{{\mathcal F}}})$.  We can write
\begin{align*}
  \text{In}_p(\omega_{\mathcal F}) & = A_I^{\mathcal F}(x,y) dx + B_I^{\mathcal F}(x,y) dy \\
  \text{In}_p(\omega_{\widetilde{{\mathcal F}}}) & = A_I^{\widetilde{\mathcal F}}(x,y) dx + B_I^{\widetilde{\mathcal F}}(x,y) dy.
\end{align*}
Take $(x_p,y_p)$ coordinates in the first chart of $E_{red}$ such that $\pi_E(x_p,y_p)=(x_p,x_p^p y_p)$ and $E_{red}=(x_p=0)$. Let $\omega^E_{\mathcal F}$ and $\omega^E_{\widetilde{\mathcal F}}$ be the strict transforms of $\omega_{\mathcal F}$ and $\omega_{\widetilde{{\mathcal F}}}$ by $\pi_E$ with
\begin{align*}
\omega^E_{\mathcal F} &= A^E_{\mathcal F}(x_p,y_p) dx_p + x_p B^E_{\mathcal F}(x_p,y_p) dy_p,  \\
\omega^E_{\widetilde{\mathcal F}} &= A^E_{\widetilde{\mathcal F}}(x_p,y_p) dx_p + x_p B^E_{\widetilde{\mathcal F}}(x_p,y_p) dy_p.
\end{align*}
Recall that we have that
\begin{align*}
  A^E_{\mathcal F}(0,y) & = A_I^{\mathcal F}(1,y) + pyB_I^{\mathcal F}(1,y); \quad & B^{E}_{\mathcal F}(0,y)& =B_I^{\mathcal F}(1,y) \\
   A^E_{\widetilde{\mathcal F}}(0,y)  & = A_I^{\widetilde{\mathcal F}}(1,y) + pyB_I^{\widetilde{\mathcal F}}(1,y); \quad & B^{E}_{\widetilde{\mathcal F}}(0,y)&=B_I^{\widetilde{\mathcal F}}(1,y)
\end{align*}
and that the Camacho-Sad indices coincide for $\mathcal F$ and $\widetilde{\mathcal F}$, that is,
$$
\mathcal{I}_{P_l^E}(\pi_E^* {\mathcal F},E_{red}) = \mathcal{I}_{P_l^E}(\pi_E^* \widetilde{\mathcal F},E_{red}), \qquad \text{ for } l=1,2,\ldots,k,
$$
where $\pi^*_E C \cap E_{red}=\{P_1^E,\ldots, P_k^E\}$.  Since $C$ has only non-singular irreducible components, if we write $P_l^E=(0,d_l^E)$ in coordinates $(x_p,y_p)$ and denote $m_{l}^C=\nu_{P_l^E}(\pi_E^*C)$, we have that
$$
A^E_{\mathcal F}(0,y)=A^E_{\widetilde{\mathcal F}}(0,y)=\prod_{l=1}^{k} (y-d_l^E)^{m_{l}^C}
$$
up to divide $\omega^E_{\mathcal F}$ and $\omega_{\widetilde{\mathcal F}}^E$ by a constant. Moreover, if we consider $(x_l,y_l)$ coordinates centered at $P_l^E$ with $x_l=x_p$ and $y_l=y_p-d_l^E$, the equality of the Newton polygons ${\mathcal N}(\pi_E^* {\mathcal F};x_l,y_l)$, ${\mathcal N}(\pi_E^* \widetilde{\mathcal F};x_l,y_l)$ and ${\mathcal N}(\pi_E^*C;x_l,y_l)$ implies
$$
\text{ord}_{y=d_l^E} (B^E_{\mathcal F}(0,y)) \geq m_{l}^C-1
\qquad \text{ord}_{y=d_l^E} (B^E_{\widetilde{\mathcal F}}(0,y)) \geq m_{l}^C-1
$$
(see \cite{Cor-2009-AnnAcBras}, Lemma 1) and we can write $B^E_{\mathcal F}(0,y)=\prod_{l=1}^{k} (y-d_l^E)^{m_{l}^C-1} \tilde{B}^E_{\mathcal F}(y)$, $B^E_{\widetilde{\mathcal F}} (0,y)=\prod_{l=1}^{k} (y-d_l^E)^{m_{l}^C-1} \tilde{B}^E_{\widetilde{\mathcal F}}(y)$. Thus,
from the definition of the Camacho-Sad index given in \eqref{eq:indice-C-S}, we have that
\begin{align*}
  \mathcal{I}_{P_l^E}(\pi_E^* {\mathcal F},E_{red}) &  = - \text{Res}_{y=d_l^E} \frac{B^E_{\mathcal F}(0,y)}{A^E_{\mathcal F}(0,y)}= - \text{Res}_{y=d_l^E} \frac{\tilde{B}^E_{\mathcal F}(y)}{\prod_{l=1}^k(y-d_l^E)}, \\
  \mathcal{I}_{P_l^E}(\pi_E^* \widetilde{\mathcal F},E_{red}) &  = - \text{Res}_{y=d_l^E} \frac{B^E_{\widetilde{\mathcal F}}(0,y)}{A^E_{\widetilde{\mathcal F}}(0,y)} = - \text{Res}_{y=d_l^E} \frac{\tilde{B}^E_{\widetilde{\mathcal F}}(y)}{\prod_{l=1}^k(y-d_l^E)}.
\end{align*}
The equality of the Camacho-Sad indices $\mathcal{I}_{P_l^E}(\pi_E^* {\mathcal F},E_{red})=\mathcal{I}_{P_l^E}(\pi_E^* \widetilde{\mathcal F},E_{red})$, for $l=1,\ldots, k$, implies $\tilde{B}^E_{\mathcal F}(y)=\tilde{B}^E_{\widetilde{\mathcal F}}(y)$ and hence
$\text{In}_p(\omega_{\mathcal F})=\text{In}_p(\omega_{\widetilde{\mathcal F}})$. %The same arguments show that $\text{In}_p(\eta)=\text{In}_p(\eta_\mu)$.

Finally, if $E$ is a non-collinear divisor,   the equality
$$\text{In}_p(J_{{\mathcal F},\mathcal{G}};x,y)= \text{In}_p(J_{\widetilde{\mathcal F},\widetilde{\mathcal{G}}};x,y)$$
is a direct consequence of Lemma~\ref{lema:parte-inicial-jac}.
\end{proof}

\begin{Remark}\label{rmk:x-cono-tg}
Note that $x=0$ can be a branch of ${\mathcal J}_{\mathcal{F},\mathcal{G}}$ although $x=0$ is not tangent to the curve $Z$. Let us consider the foliations $\mathcal F$ and $\mathcal G$ given by $\omega=0$ and $\eta=0$ with
\begin{align*}
  \omega & = (xy-6x^2)dx + (y^2-6xy+10x^2)dy \\
  \eta & = -6x^5 dx + 3y^2 dy.
\end{align*}
Thus ${\mathcal J}_{\mathcal{F},\mathcal{G}}$ is given by $J(x,y)=0$ with
$$J(x,y)=3x (y^3-6xy^2+2x^4y^2-12x^5y+20x^6).$$
In this example, if we consider the blow-up $\pi_1: X_1 \to ({\mathbb C}^2,0)$ of the origin, the first divisor $E^1$ is non-collinear. Thus the result of Lemma~\ref{lema:parte-inicial-jac} above holds: we have that $\text{In}_1(\omega)=\omega$, $\text{In}_1(\eta)=3y^2dy$ and hence $\text{In}_1(J)= 3xy^2(y-6x)$.

Note that the rational function ${\mathcal M}_{E^1}(z)$ is given by
$${\mathcal M}_{E^1}(z)=- \frac{z-6}{z(z-1)(z-2)(z-3)}$$
which determines the branch $J^{E^1}_{nc}$ whose tangent cone is given by $y-6x=0$ but we cannot determine the branch $x=0$ of ${\mathcal J}_{\mathcal{F},\mathcal{G}}$ (see statement of Theorem~\ref{th:descomposicion-no-sing}).
\end{Remark}

The strategy to prove the results about the jacobian curve is based on the properties that share a foliation and its logarithmic model. Next lemma will allow to follow this strategy.
\begin{Lemma}\label{lema:puntos-jac-log-nolog}
Consider foliations ${\mathcal F}, {\mathcal L}_\lambda^C \in {\mathbb G}_{C,\lambda}$ and ${\mathcal G}, {\mathcal L}_\mu^D \in {\mathbb G}_{D,\mu}$. Let ${\mathcal J}_{\mathcal{F},\mathcal{G}}$   be the jacobian curve of $\mathcal F$ and $\mathcal G$ and $\mathcal{J}_{\lambda,\mu}$ the    jacobian curve of $\mathcal{L}_\lambda^C$ and $\mathcal{L}_{\mu}^D$. Let $E$ be an irreducible component of $\pi_Z^{-1}(0)$. If $E$ is a non-collinear divisor, we have that
$$\pi_E^* {\mathcal J}_{\mathcal{F},\mathcal{G}} \cap E_{red} = \pi_E^* \mathcal{J}_{\lambda,\mu} \cap E_{red}$$
and the multiplicities satisfy that
$$\nu_P(\pi_E^* {\mathcal J}_{\mathcal{F},\mathcal{G}}) = \nu_P (\pi_E^* \mathcal{J}_{\lambda,\mu})$$
at each point $P \in \pi_E^* {\mathcal J}_{\mathcal{F},\mathcal{G}} \cap E_{red}$.
\end{Lemma}
\begin{proof}
Let $E$ be an irreducible component of $\pi_Z^{-1}(0)$ with $v(E)=p$ and take $(x,y)$ coordinates adapted to $E$. The result is a direct consequence of the equality
$$
{\rm In}_p(J_{\mathcal{F},\mathcal{G}};x,y)={\rm In}_p(J_{\lambda,\mu};x,y)
$$
given in Lemma~\ref{lema:parte-incial-modelo-log} provided that $E$ is a non-collinear divisor, where jacobian curves ${\mathcal J}_{\mathcal{F},\mathcal{G}}$, $\mathcal{J}_{\lambda,\mu}$ are given by $J_{\mathcal{F},\mathcal{G}}(x,y)=0$ and $J_{\lambda,\mu}(x,y)=0$ respectively.
\end{proof}

\section{Properties of the jacobian curve} \label{sec:jacobian} %\label{sec:pruebas}

Let us consider two singular foliations  $\mathcal F$ and $\mathcal G$    in $({\mathbb C}^2,0)$ defined by the 1-forms $\omega=0$ and $\eta=0$ with $\omega=A(x,y) dx + B(x,y) dy$ and $\eta=P(x,y) dx + Q(x,y) dy$.
Recall that the jacobian curve ${\mathcal J}_{{\mathcal F},{\mathcal G}}$ is defined by $J(x,y)=0$ where
$$J(x,y)=A(x,y) Q(x,y) - B(x,y) P(x,y).$$
Next remark shows that the jacobian curve behaves well by a change of coordinates.
\begin{Remark}\label{rk:cambio-coord}

If $F:({\mathbb C}^2,0) \to  ({\mathbb C}^2,0)$ is a change of coordinates with $F=(F_1,F_2)$, the Jacobian curve of $F^*{\mathcal F}$ and $F^* {\mathcal G}$ is given by
$$
\left| \begin{array}{cc}
A \circ F & B \circ F\\
P \circ F & Q \circ F
\end{array} \right| \cdot \left| \begin{array}{cc}
F_{1,x} & F_{1,y} \\
F_{2,x} & F_{2,y}
\end{array} \right|  =0.$$
Thus, the curve ${\mathcal J}_{F^*{\mathcal F},F^* {\mathcal G}}$ is defined by $J \circ F=0$. Hence, ${\mathcal J}_{F^*{\mathcal F},F^* {\mathcal G}}= F^{-1}({\mathcal J}_{{\mathcal F},{\mathcal G}})$.

In particular, we get that the analytic type of the jacobian curve of the foliations $\mathcal{F}$ and $\mathcal{G}$ is an invariant of the analytic type of the pair of foliations $\mathcal F$ and $\mathcal G$.
\end{Remark}

Assume that ${\mathcal F} \in {\mathbb G}_{C,\lambda}$ and ${\mathcal G} \in {\mathbb G}_{D,\mu}$ where
 $C=\cup_{i=1}^r C_i$ and $D=\cup_{i=1}^s D_i$ are two plane curves in $({\mathbb C}^2,0)$ without common irreducible components. In this section, we will assume that  the curve $Z=C \cup D$ has only non-singular irreducible components and we consider $\pi_Z: X_Z \to (\mathbb{C}^2,0)$ the minimal reduction of singularities of $Z$. Note that since $Z$ has only non-singular irreducible components,  the centers of the blow-ups to obtain $\pi_Z$ are all free infinitely near points of $Z$. The general case will be treated in section~\ref{sec:caso-general}.

In section~\ref{sec:collinear} we have introduced all the notions we need to state the results concerning the properties of the jacobian curve.
If we want to compute the infinitely near points of the jacobian curve of two foliations, Lemma~\ref{lema:puntos-jac-log-nolog} will allow to do computations for the jacobian curve of two logarithmic foliations and then get the result for the jacobian curve of two generalized curve foliations.
The first result gives the multiplicity of the jacobian curve at a point in the reduction of singularities in terms of the multiplicities of the curves $C$ and $D$. In particular, we obtain that all infinitely near points of the jacobian curve in the first chart of a divisor $E$ of $\pi_Z^{-1}(0)$ are either  infinitely near points of $Z$ or a point in $M(E)$. Note that $x=0$ can be tangent to the jacobian curve although it is not tangent to $Z$ and we cannot control this with the rational function ${\mathcal M}_E(z)$ (see Remark~\ref{rmk:x-cono-tg}). More precisely, given a divisor $E$, we can fix coordinates $(x,y)$ adapted to $E$ (see Remark~\ref{rk:coordenadas-adaptadas}) and we can denote by $E_{red}^*$ the points in the first chart of $E_{red}$. Then,  we have
\begin{Theorem}\label{th-multiplicidad-jac}
Let $E$ be an irreducible component  of $\pi^{-1}_Z(0)$ and assume that $E$ is a non-collinear divisor. Given any $P \in E_{red}^*$, we have that
$$
\nu_P(\pi_E^*{\mathcal J}_{\mathcal{F},\mathcal{G}})= \nu_P(\pi_E^* C) + \nu_P(\pi_E^*D) + \tau_E(P).
$$
In particular, if $P \in E_{red}^*$ with $\nu_P(\pi_E^*{\mathcal J}_{\mathcal{F},\mathcal{G}})>0$, then $P$ is an infinitely near point of $Z$ or a point in $M(E)$.
%$P \in C(E) \cup N(E) \cup M(E)$.
\end{Theorem}

\begin{proof}
We prove here the result for logarithmic foliations.
The general case, when $\mathcal{F}$ and $\mathcal{G}$ are not necessarily logarithmic foliations is consequence of Lemma~\ref{lema:puntos-jac-log-nolog}.

Consider the logarithmic foliations ${\mathcal L}_\lambda^C$ and  ${\mathcal L}_\mu^{D}$  given by $\omega_\lambda=0$ and $\eta_\mu=0$ with
\begin{align*}
\omega_\lambda &  = \prod_{i=1}^r(y-\alpha_i(x)) \sum_{i=1}^{r} \lambda_i \frac{d(y-\alpha_i(x))}{y-\alpha_i(x)} \\
\eta_\mu & = \prod_{i=1}^s(y-\beta_i(x)) \sum_{i=1}^{s} \mu_i \frac{d(y-\beta_i(x))}{y-\beta_i(x)}
\end{align*}
where the curve $C_i$ is given by $y-\alpha_i(x)=0$ with $\alpha_i(x)=\sum_{j=1}^{\infty} a_j^i x^j \in {\mathbb C}\{x\}$ and the curve $D_i$ is given by
$y-\beta_i(x)=0$ with $\beta_i(x)=\sum_{j=1}^{\infty} b_j^i x^j \in {\mathbb C}\{x\}$. Let us denote ${\mathcal J}_{\lambda,\mu}$  the jacobian curve of ${\mathcal L}_\lambda^C$ and  ${\mathcal L}_\mu^{D}$ which is defined by $J_{\lambda,\mu}(x,y)=0$ with
$$ J_{\lambda,\mu}(x,y) =A_{\lambda}(x,y) Q_{\mu}(x,y) - B_{\lambda}(x,y) P_{\mu}(x,y)$$
where we write $\omega_\lambda=A_\lambda(x,y) dx + B_\lambda(x,y) dy$ and $\eta_\mu=P_\mu(x,y) dx + Q_\mu(x,y) dy$. More precisely, we can write

\begin{equation}\label{eq:jacobiano-lambda-mu}
J_{\lambda,\mu}(x,y) = f(x,y) g(x,y) \left|
\begin{array}{cc}
-\sum\limits_{i=1}^{r} \lambda_i \frac{\alpha_i'(x)}{y-\alpha_i(x)} & \sum\limits_{i=1}^{r} \frac{\lambda_i}{y-\alpha_i(x)} \\ & \\
-\sum\limits_{i=1}^{s} \mu_i \frac{\beta'_i(x)}{y-\beta_i(x)} & \sum\limits_{i=1}^{s} \frac{\mu_i}{y-\beta_i(x)}
\end{array}
\right|
\end{equation}
where $f(x,y)=\prod_{i=1}^r(y-\alpha_i(x))$ and $g(x,y)=\prod_{i=1}^s(y-\beta_i(x))$ are equations of the curves $C$ and $D$ respectively.

Let $\pi_Z: X_Z \to ({\mathbb C}^2,0)$ be the minimal reduction of singularities of $Z$. Take $E$ a bifurcation divisor of $G(Z)$ with $v(E)=p$ and consider $\pi_E : X_E \to ({\mathbb C}^2,0)$ the reduction of $\pi_Z$ to $E$. Since the jacobian curve behaves well by a change of coordinates (see Remark~\ref{rk:cambio-coord}),  we can assume that the coordinates $(x,y)$ are adapted to $E$. Take $(x_p,y_p)$ coordinates in the first chart of $E_{red} \subset X_E$ such that $\pi_E(x_p,y_p)=(x_p,x_p^p y_p)$ and $E_{red}=(x_p=0)$.
Let us compute the strict transform of ${\mathcal J}_{\lambda,\mu}$ by $\pi_E$. %We consider two situations: $E$ being the first bifurcation divisor of $G(Z)$ or not.

Let us denote $I=\{1,\ldots, r\}$,  $J=\{1,\ldots, s\}$, $I^E=\{ i \in I \ : \ E \text{ belongs to}  \linebreak \text{ the geodesic of } C_i\}$ and $J^E = \{ j \in J \ : \ E \text{ belongs to the geodesic of } D_i\}$.
We can write
\begin{align*}
\omega_\lambda &  = f(x,y) \left[ \sum_{i \in I \smallsetminus I^E} \lambda_i \frac{ -\alpha'_i(x) dx + dy}{y-\alpha_i(x)} + \sum_{i \in  I^E} \lambda_i \frac{ -\alpha'_i(x) dx + dy}{y-\alpha_i(x)} \right] \\
\eta_\mu &= g(x,y) \left[ \sum_{i \in J \smallsetminus J^E} \mu_i \frac{ -\beta'_i(x) dx + dy}{y-\beta_i(x)} + \sum_{i \in  J^E} \mu_i \frac{ -\beta'_i(x) dx + dy}{y-\beta_i(x)} \right]
\end{align*}
and hence, the jacobian curve  ${\mathcal J}_{\lambda,\mu}$ is given by $J_{\lambda,\mu}(x,y)=0$ with
$$
J_{\lambda,\mu}(x,y) = f(x,y) g(x,y) \left|
\begin{array}{cc}
  \sum\limits_{i \in I \setminus I^E} \lambda_i \frac{-\alpha'_i(x)}{y-\alpha_i(x)} + \sum\limits_{i \in   I^E} \lambda_i \frac{-\alpha'_i(x)}{y-\alpha_i(x)}  &  \sum\limits_{i \in I \setminus I^E} \frac{\lambda_i}{y-\alpha_i(x)} + \sum\limits_{i \in   I^E}  \frac{\lambda_i}{y-\alpha_i(x)}   \\
  \sum\limits_{i \in J \setminus J^E} \mu_i \frac{-\beta'_i(x)}{y-\beta_i(x)} + \sum\limits_{i \in   J^E} \mu_i \frac{-\beta'_i(x)}{y-\beta_i(x)}  &  \sum\limits_{i \in J \setminus J^E} \frac{\mu_i }{y-\beta_i(x)} + \sum\limits_{i \in   J^E}  \frac{\mu_i}{y-\beta_i(x)}
\end{array}
\right|
$$
Since $v(E)=p$ and $(x,y)$ are coordinates adapted to $E$, we have that
\begin{align*}
  \text{ord}_x(\alpha_i(x)) &  \geq p \quad \text{ if } \quad  i \in I^E; \quad &  \text{ord}_x(\alpha_i(x))=n_i & < p  \quad \text{ if } \quad i \in I \smallsetminus I^E;   \\
\text{ord}_x(\beta_i(x)) & \geq p \quad \text{ if } \quad i \in J^E; &  \text{ord}_x(\beta_i(x))=o_i&  <p \quad \text{ if } \quad i \in J \smallsetminus J^E.
\end{align*}
Thus, $\alpha_i(x) = \sum_{j \geq p} a_j^i x^j $ if $i \in I^E$ and
$\beta_i(x) = \sum_{j \geq p} b_j^i x^j$ if $i \in J^E$, but $\alpha_i (x) = \sum_{j \geq n_i} a_j^i x^j$ with $n_i < p$ if $i \in I \smallsetminus I^E$ and  $\beta_i (x) = \sum_{j \geq o_i} b_j^i x^j$ with $o_i < p$ if $i \in J \smallsetminus J^E$. Then, $J_{\lambda,\mu}(x_p,x_p^p y_p)$ is given by
$$
J_{\lambda,\mu}(x_p,x_p^p y_p)  = f(x_p,x_p^p y_p) g(x_p,x_p^p y_p) M_E^*(x_p,y_p)
$$
with
\begin{align*}
 &  M_E^*(x_p,y_p)    =  \\
&  \frac{1}{x_p^{p+1}} \left|
\begin{array}{cc}
  \sum\limits_{i \in I \smallsetminus I^E} \lambda_i \frac{-n_i a_{n_i}^i + x_p(\cdots) }{x_p^{p-n_i}y_p-a_{n_i}^i + x_p(\cdots)} + \sum\limits_{i \in   I^E} \lambda_i \frac{-p a_p^i + x_p(\cdots)}{y_p-a_p^i + x_p(\cdots) }  &  \sum\limits_{i \in I \smallsetminus I^E} \frac{\lambda_i x_p^{p-n_i}}{x_p^{p-n_i}y_p-a_{n_i}^i + x_p(\cdots)} + \sum\limits_{i \in   I^E}  \frac{\lambda_i}{y_p-a_p^i + x_p(\cdots)}   \\
  \sum\limits_{i \in J \smallsetminus J^E} \mu_i \frac{-o_i b_{o_i}^i + x_p(\cdots) }{x_p^{p-o_i} y_p-b_{o_i}^i + x_p(\cdots)} + \sum\limits_{i \in   J^E} \mu_i \frac{-pb_p^i + x_p(\cdots)}{y_p-b_p^i + x_p(\cdots)}  &  \sum\limits_{i \in J \smallsetminus J^E} \frac{\mu_i x_p^{p-o_i} }{x_p^{p-o_i} y_p-b_{o_i}^i + x_p(\cdots)} + \sum\limits_{i \in   J^E}  \frac{\mu_i}{y_p-b_p^i + x_p(\cdots)}
\end{array}
\right| \\
   & =  \frac{1}{x_p^{p+1}} M_E(x_p,y_p)
\end{align*}
If $M_E(0,y_p) \not \equiv 0$, then the points $\pi_E^* {\mathcal J}_{\lambda,\mu} \cap E_{red}$, in the first chart of $E_{red}$, are given by $x_p=0$ and $J_E(y_p)=0$ where
$$J_E(y_p)= \tilde{f}(0,y_p) \tilde{g}(0,y_p) M_E(0,y_p),$$

Let $\{R_1^E,\ldots, R_{b_E}^E\}$ be the union of the singular points of $\pi_E^*{\mathcal L}_{\lambda}^C$ and $\pi_E^* {\mathcal L}_\mu^D$ in the first chart of $E_{red}$ where $R_l^E=(0,c_l^E)$ in  coordinates $(x_p,y_p)$. Note that $\{R_1^E,\ldots, R_{b_E}^E\}=\pi_E^* Z \cap E_{red}$.  Denote  $m_l^C=\nu_{R_l^E}(\pi_E^* C) =\sharp  \{ i \in \{1, \ldots, r\} \ : \ \pi_E^*C_i \cap E_{red} = \{R_l^E\}\}$ and  $m_l^D=\nu_{R_l^E}(\pi_E^*D)=\sharp \{ j \in \{1, \ldots, s\} \ : \ \pi_E^*D_j \cap E_{red} = \{R_l^E\}\}$ for $l=1,2,\ldots, b_E$.

We have that
\begin{align*}
M_E(0,y_p) & =  \left|
\begin{array}{ccc}
  \sum\limits_{i \in I \smallsetminus I^E} \lambda_i n_i   + \sum\limits_{i \in   I^E} \lambda_i \frac{-p a_p^i}{y_p-a_p^i}  &  \quad & \sum\limits_{i \in   I^E}  \frac{\lambda_i}{y_p-a_p^i}   \\ & \\
  \sum\limits_{i \in J \smallsetminus J^E} \mu_i o_i   + \sum\limits_{i \in   J^E} \mu_i \frac{-pb_p^i}{y_p-b_p^i}  &   & \sum\limits_{i \in   J^E}  \frac{\mu_i}{y_p-b_p^i}
\end{array}
\right| \\ & \\
& =  \left|
\begin{array}{ccc}
  \sum\limits_{i \in I \smallsetminus I^E} \lambda_i n_i   + \sum\limits_{i \in   I^E} \lambda_i p    &  \quad & \sum\limits_{i \in   I^E}  \frac{\lambda_i}{y_p-a_p^i}   \\ &  \\
  \sum\limits_{i \in J \smallsetminus J^E} \mu_i o_i   + \sum\limits_{i \in   J^E} \mu_i p &   & \sum\limits_{i \in   J^E}  \frac{\mu_i}{y_p-b_p^i}
\end{array}
\right| \\
& = -\kappa_E({\mathcal L}_\lambda^C) \kappa_E({\mathcal L}_\mu^D)  \left|
\begin{array}{ccc}
  1   &  \quad & \sum\limits_{l=1}^{b_E}  \frac{\mathcal{I}_{R_l^E}(\pi_E^*{\mathcal L}_\lambda^C,E_{red})}{y_p-c_l^E}   \\ & \\
  1   &   & \sum\limits_{l=1}^{b_E}  \frac{\mathcal{I}_{R_l^E}(\pi_E^*{\mathcal L}_\mu^D,E_{red})}{y_p-c_l^E}
\end{array}
\right| \\
\end{align*}
where $\kappa_E({\mathcal L}_\lambda^C) =  \sum\limits_{i \in I \smallsetminus I^E} \lambda_i n_i   + \sum\limits_{i \in   I^E} \lambda_i p $ and $\kappa_E({\mathcal L}_\mu^D) = \sum\limits_{i \in J \smallsetminus J^E} \mu_i o_i   + \sum\limits_{i \in   J^E} \mu_i p$ are the residues of the logarithmic foliations along the divisor $E$ (see equation~\eqref{eq:peso-divisor-log}) and we use the expression of the Camacho-Sad index for a logarithmic foliation given in~\eqref{eq:camacho-sad-log}. Consequently, we obtain that
$$M_E(0,y_p)= -\kappa_E({\mathcal L}_\lambda^C) \kappa_E({\mathcal L}_\mu^D) {\mathcal M}_E(y_p)$$
where $\mathcal{M}_E(z)$ is the rational function associated to the divisor $E$ for the foliations  $\mathcal{L}^C_\lambda$ and $\mathcal{L}^D_\mu$ (see expression~\eqref{eq:funcion-meromorfa})
Then, the points $\pi_E^*{\mathcal J}_{\lambda,\mu} \cap E_{red}$, in the first chart of $E_{red}$, are given by $x_p=0$ and
$$
\prod_{i=1}^{b_E} (y_p -c_l^E)^{m_l^C+m_l^D} {\mathcal M}_E(y_p)=0.
$$
(note that the curve $C \cup D$ has only non-singular irreducible components). Let $\{q_1,\ldots, q_{s_E}\}$ be the set of zeros of ${\mathcal M}_E(z)$. For $l=1,2, \ldots, s_E$, put $Q_l^E=(0,q_l)$ and denote by $t_{Q_l^E}$ the multiplicity of $q_l$ as a zero of ${\mathcal M}_E(z)$. Thus, the points in $\pi_E^* {\mathcal J}_{\lambda,\mu} \cap E_{red}$ belong to $\operatorname{Col}(E) \cup \operatorname{NCol}(E) \cup M(E)$. Moreover, the multiplicity of $\pi_E^* {\mathcal J}_{\lambda,\mu}$ at a point $P \in E_{red}$, in the first chart of $E_{red}$,  is given by
$$\nu_P(\pi_E^* {\mathcal J}_{\lambda,\mu}) = \nu_P(\pi_E^*C) + \nu_P(\pi_E^* D) + \tau_E(P)$$
where $\tau_E(P)$ was defined by the expression~\eqref{eq:tau_E}.
\end{proof}

\begin{Remark}\label{rmk:x-tangente}
With the notations of the  proof above, if the first divisor $E^1$ is non-collinear, we have that the tangent cone of ${\mathcal J}_{\mathcal{L}^C_\lambda,\mathcal{L}^D_\mu}$ is given by $J_1(x,y)=0$
where
\begin{align*}
 J_1(x,y)  =& - \left(\sum_{i=1}^{r} \lambda_i a_1^i \prod_{j\neq i} (y-a_1^j x ) \right) \left( \sum_{i=1}^{s} \mu_i \prod_{j \neq i} (y-b_1^j x) \right) \\
   & + \sum_{i=1}^{r} \lambda_i \prod_{j \neq i} (y-a_1^jx) \sum_{i=1}^{s} \mu_i b_1^i \prod_{j \neq i} (y - b_1^j x)
\end{align*}
Thus, $x=0$ is not tangent to the jacobian curve ${\mathcal J}_{\mathcal{L}^C_\lambda,\mathcal{L}^D_\mu}$ provided that
\begin{equation}\label{eq:x-no-tg}
\kappa_{E^1}({\mathcal L}_\lambda^C) \sum_{i=1}^s \mu_i b_1^i - \kappa_{E^1}({\mathcal L}_\mu^D) \sum_{i=1}^r \lambda_ i a_1^i \neq 0
\end{equation}
where we recall that $\kappa_{E^1}({\mathcal L}_\mu^D)= \sum_{i=1}^{r} \lambda_i$ and $\kappa_{E^1}({\mathcal L}_\mu^D)=\sum_{i=1}^s \mu_i$.
By Lemma~\ref{lema:parte-incial-modelo-log}, the above remarks hold for the jacobian curve ${\mathcal J}_{\mathcal{F},\mathcal{G}}$ for any $\mathcal{F} \in \mathbb{G}_{C,\lambda}$, $\mathcal{G} \in \mathbb{G}_{D,\mu}$.
In the example given in Remark~\ref{rmk:x-cono-tg} we have that $\sum_{i=1}^{3} \lambda_i a_1^i = \sum_{i=1}^{3} \mu_i b_1^i =0$ and hence the condition in \eqref{eq:x-no-tg} does not hold whereas in the example given in Remark~\ref{rmk:CE-ME} condition in \eqref{eq:x-no-tg} holds and hence $x=0$ is not tangent to the jacobian curve.
\end{Remark}

The reader can find in appendix~\ref{ap:equisingularidad} some definitions related to the equisingularity data of curves used in the statements of the following results.

\medskip
Consider now   two consecutive bifurcation divisors $E$ and $E'$ in $G(Z)$ such that  $E'$  arises from $E$ at $P$. As we have explained in section~\ref{sec:collinear}, this means that there is a chain of consecutive divisors
$$E_0=E < E_1 < \cdots < E_{k-1} < E_k=E'$$
with $b_{E_l}=1$ for $l=1, \ldots, k-1$ and the morphism $\pi_{E'}=  \pi_E  \circ \sigma$ where $\sigma: X_{E'} \to (X_E,P)$ is a composition of $k$ punctual blow-ups
%$$
\begin{equation} \label{eq:bifurcacion-consecutivos_2}
(X_E, P) \overset{\sigma_1}{\longleftarrow} (X_{E_1},P_1)  \overset{\sigma_2}{\longleftarrow}  \cdots \overset{\sigma_{k-1}}{\longleftarrow} (X_{E_{k-1}},P_{k-1})  \overset{\sigma_{k}}{\longleftarrow} X_{E'}.
\end{equation}

Now we can explain the behaviour of the branches of the jacobian curve going through a non-collinear point.
Next corollary states that the branches of the jacobian curve going through a non-collinear point $P$ in a bifurcation divisor as above go through the points $P_1, \ldots, P_{k-1}$  given in the sequence \eqref{eq:bifurcacion-consecutivos_2}, that is, the divisor $E'$ is in the geodesic of those branches of $\mathcal{J}_{\mathcal{F},\mathcal{G}}$ going through $P$ in $E_{red}$.
\begin{Corollary}\label{cor:consecutive-divisors}
Let $E$ and $E'$ be two consecutive bifurcation divisors in $G(Z)$ with $E {<_P} E'$. If $P \in \operatorname{NCol}(E)$,  we have that
$$
\nu_P(\pi_E^* {\mathcal J}_{\mathcal{F},\mathcal{G}}) = \sum_{Q \in E'_{red}} \nu_Q(\pi_{E'}^* {\mathcal J}_{\mathcal{F},\mathcal{G}})
$$
In particular, we get that there is no irreducible component $\delta$ of ${\mathcal J}_{\mathcal{F},\mathcal{G}}$ such that $\pi_{E'}^*\delta$ is attached to some intermediate component $E_i$, $1\leq i\leq k-1$, in the chain $E < E_1 < \cdots < E_{k-1} < E'$.
Moreover,
\begin{equation}\label{eq1:consecutive-divisors}
1 +\sum_{Q \in M(E')} t_Q  = \sharp \operatorname{NCol}(E').
\end{equation}
Hence $E'$ is non-collinear.
\end{Corollary}
\begin{Remark}
Note that from the above result, we get that there is no irreducible component $\delta$ of ${\mathcal J}_{\mathcal{F},\mathcal{G}}$ such that
$$v(E) < {\mathcal C} (\delta,\gamma_{E'}) < v(E')$$
where $\gamma_{E'}$ is a $E'$-curvette.
\end{Remark}
\begin{proof}[Proof of Corollary~\ref{cor:consecutive-divisors}]
Let $E$ and $E'$ be two consecutive bifurcation divisors in $G(Z)$ with $E <_P E'$ and assume that $P \in \operatorname{NCol}(E)$, thus $\Delta_E(P)\neq 0$. By Theorem~\ref{th-multiplicidad-jac} we have that
$$\nu_P(\pi_E^*{\mathcal J}_{\mathcal{F},\mathcal{G}})=\nu_P(\pi_E^*C)+\nu_P(\pi_E^*D)-1.$$
Recall that $E <_P E'$ implies the existence of
a chain of consecutive divisors
$$E_0=E < E_1 < \cdots < E_{k-1} < E_k=E'$$
with $b_{E_l}=1$ for $l=1, \ldots, k-1$ and the morphism $\pi_{E'}= \pi_E \circ \sigma  $ where $\sigma: X_{E'} \to (X_E,P)$ is a composition of $k$ punctual blow-ups
\begin{equation*} %\label{eq:bifurcacion-consecutivos_2}
(X_E, P) \overset{\sigma_1}{\longleftarrow} (X_{E_1},P_1)  \overset{\sigma_2}{\longleftarrow}  \cdots \overset{\sigma_{k-1}}{\longleftarrow} (X_{E_{k-1}},P_{k-1})  \overset{\sigma_{k}}{\longleftarrow} X_{E'}.
\end{equation*}
Since $P$ is non-collinear, then each $P_i$ is non-collinear by Lemma~\ref{lema:explosio-colineal}, thus $\Delta_{E_i}(P_i) \neq 0$ and hence $M(E_i)=\emptyset$ for $i=1,\ldots, k-1$. In particular, using again Theorem~\ref{th-multiplicidad-jac}, we have that
$$
\nu_{P_i}(\pi_{E_i}^*{\mathcal J}_{\mathcal{F},\mathcal{G}})=\nu_{P_i}(\pi_{E_i}^*C)+\nu_{P_i}(\pi_{E_i}^*D)-1. % = \nu_P(\pi_E^*C)+m_P(\pi_E^*D)-1.
$$
Since the curves $C$ and $D$ have only non-singular irreducible components, and $P_i$ is the only infinitely near point of both curves in $E_i$, we have that $\nu_{P_i}(\pi_{E_i}^*C)=\nu_P(\pi_E^*C)$ and $\nu_{P_i}(\pi_{E_i}^*D)=\nu_P(\pi_E^*D)$ for all $i=1,\ldots, k-1$. Consequently,
$$
\nu_{P_i}(\pi_{E_i}^*{\mathcal J}_{\mathcal{F},\mathcal{G}})=\nu_P(\pi_E^*{\mathcal J}_{\mathcal{F},\mathcal{G}}), \quad \text{ for } i=1,\ldots, k-1.
$$
Since $E'$ is a bifurcation divisor, we get that
$$
\nu_{P_{k-1}}(\pi_{E_{k-1}}^*{\mathcal J}_{\mathcal{F},\mathcal{G}})=\sum_{Q \in E'_{red}} \nu_Q(\pi_{E'}^* {\mathcal J}_{\mathcal{F},\mathcal{G}}).
$$
Hence, from all the equalities above, we deduce that
$$
\nu_P(\pi_E^* {\mathcal J}_{\mathcal{F},\mathcal{G}}) = \sum_{Q \in E'_{red}} \nu_Q(\pi_{E'}^* {\mathcal J}_{\mathcal{F},\mathcal{G}}),
$$
which proves the first statement of the corollary.
Finally, in order to prove the equality given in~\eqref{eq1:consecutive-divisors}, it is enough to show that
$$\sum_{R \in \operatorname{NCol}(E')} \Delta_{E'}(R) \neq 0$$
 by Lemma~\ref{lema:E-no-colineal}.
Let us assume that $\sum_{R \in \operatorname{NCol}(E')} \Delta_{E'}(R) = 0$, which implies
$$
\sum_{R \in \operatorname{NCol}(E')} \mathcal{I}_{R}(\pi_{E'}^* {\mathcal F},E_{red}')= \sum_{R \in \operatorname{NCol}(E')} \mathcal{I}_{R}(\pi_{E'}^* {\mathcal G},E_{red}')
$$
and, by the properties of the Camacho-Sad indices, we deduce that
$$
 \mathcal{I}_{\widetilde{P}_{k-1}}(\pi_{E'}^* {\mathcal F},E_{red}')=  \mathcal{I}_{\widetilde{P}_{k-1}}(\pi_{E'}^* {\mathcal G},E_{red}')
$$
where we denote $\widetilde{P}_{k-1} = \widetilde{E}_{k-1,red} \cap E_{red}'$ and $\widetilde{E}_{k-1,red}$ is the strict transform of $E_{k-1,red}$ by $\sigma_{k}$. Since $\mathcal F$ and $\mathcal G$ are generalized curve foliations, then $\widetilde{P}_{k-1}$ is a simple singularity for $\pi_{E'}^* {\mathcal F}$ and $\pi_{E'}^*{\mathcal G}$ and hence we have that
\begin{align*}
\mathcal{I}_{\widetilde{P}_{k-1}}(\pi_{E'}^* {\mathcal F},E_{red}') \cdot \mathcal{I}_{\widetilde{P}_{k-1}}(\pi_{E'}^* {\mathcal F},\widetilde{E}_{k-1,red})&=1 \\
\mathcal{I}_{\widetilde{P}_{k-1}}(\pi_{E'}^* {\mathcal G},E_{red}') \cdot \mathcal{I}_{\widetilde{P}_{k-1}}(\pi_{E'}^* {\mathcal G},\widetilde{E}_{k-1,red}) &=1.
\end{align*}
Consequently, given that
\begin{align*}
 \mathcal{I}_{P_{k-1}}(\pi_{E_{k-1}}^* {\mathcal F},E_{k-1,red}) & =\mathcal{I}_{\widetilde{P}_{k-1}}(\pi_{E'}^* {\mathcal F},\widetilde{E}_{k-1,red}) +1 = \frac{1}{\mathcal{I}_{\widetilde{P}_{k-1}}(\pi_{E'}^* {\mathcal F},E_{red}')} +1, \\
\mathcal{I}_{P_{k-1}}(\pi_{E_{k-1}}^* {\mathcal G},E_{k-1,red}) & =\mathcal{I}_{\widetilde{P}_{k-1}}(\pi_{E'}^* {\mathcal G},\widetilde{E}_{k-1,red}) +1 = \frac{1}{\mathcal{I}_{\widetilde{P}_{k-1}}(\pi_{E'}^* {\mathcal G},E_{red}')} +1,
\end{align*}
we obtain that $\mathcal{I}_{P_{k-1}}(\pi_{F_{k-1}}^* {\mathcal F},E_{k-1,red})=\mathcal{I}_{P_{k-1}}(\pi_{F_{k-1}}^* {\mathcal G},E_{k-1,red})$.
This  implies that $\Delta_{E_{k-1}}(P_{k-1})=0$ which is not possible since $P_{k-1}$ is a non-collinear point by Lemma~\ref{lema:explosio-colineal}. This ends the proof.
\end{proof}
In order to explain the behaviour of the branches of the jacobian curve going through a collinear point, we introduce the following definition:
\begin{Definition}\label{def:cover}
Let $E$ be a bifurcation divisor of $G(Z)$ and take $P$ a collinear point of $E$. We say that a set of non-collinear bifurcation divisors $\{E_1,\ldots, E_u\}$ is a {\em (non-collinear) cover\/} of $E$ at $P$ if the following conditions hold:
\begin{itemize}
    \item[(i)] $E$ is in the geodesic of each $E_l$;
    \item[(ii)] if $\{E_1^l,\ldots, E_{r(l)}^l\}$ is the set of all  bifurcation divisors in the geodesic of $E_l$  with
    $$E <_P E_1^l < \ldots <E_{r(l)}^l < E_l$$
    then either $r(l)=0$ or else each $E_j^l$ is collinear;
    \item[(iii)] if $Z_j$ is an irreducible component of $Z$ with $\pi_E^* Z_j \cap E_{red}=\{P\}$, then there exists a divisor $E_l$ in the cover such that $\pi_{E_l}^* Z_j \cap E_{l} \neq \emptyset$, that is, there is a divisor $E_l$ in the cover which is in the geodesic of $Z_j$.
\end{itemize}
\end{Definition}
 Given a collinear point $P$ of $E$, there is a unique cover of $E$ at $P$. We can find it as follows: take an irreducible component $Z_j$ of $Z$ with $\pi_E^* Z_j \cap E_{red}=\{P\}$. Let $E'$ be the consecutive bifurcation divisor to $E$ with $E<_P E'$ belonging to the geodesic of $Z_j$. If $E'$ is non-collinear, then $E'$ is one of the bifurcation divisors in the cover of $E$ at $P$, otherwise we repeat the process above with the following bifurcation divisor in the geodesic of $Z_j$. Since the maximal bifurcation divisors are non-collinear (see Remark~\ref{rk:max-div-non-colinear}), we will always find a non-collinear divisor in the geodesic of $Z_j$ verifying condition $(iii)$  in the above definition.

\begin{Theorem}\label{th:colinear-point}
Consider a non-collinear bifurcation divisor $E$ of $G(Z)$ and a collinear point $P$ of $E$. Take a cover $\{E_1,\ldots, E_u\}$ of $E$ at $P$. Then
$$
\nu_P(\pi_E^* {\mathcal J}_{\mathcal{F},\mathcal{G}})-\sum_{l=1}^{u} \sum_{Q \in E_{l,red}} \nu_Q(\pi_{E_l}^* {\mathcal J}_{\mathcal{F},\mathcal{G}})=
t_P + \sum_{l=1}^u (\sharp \operatorname{NCol}(E_l) - t(E_l)).
$$
Consequently, there is a curve $J_{P}^E$ composed by irreducible components of ${\mathcal J}_{\mathcal{F},\mathcal{G}}$ such that, if $\delta$ is a branch of $J_{P}^E$,
\begin{itemize}
  \item $\pi_E^* \delta \cap E_{red}=\{P\}$
  \item ${\mathcal C}(\delta,\gamma_{E_l}) < v(E_l)$ for $l=1,\ldots, u,$
where  $\gamma_{E_l}$ is any $E_l$-curvette.
\end{itemize}
Moreover, we have that
$$\nu_0(J_{P}^E)= t_P + \sum_{l=1}^u (\sharp \operatorname{NCol}(E_l) - t(E_l)).$$
\end{Theorem}

\begin{proof} %[Proof of Theorem~\ref{th:colinear-point}]
Let $E$ be a non-collinear bifurcation divisor of $G(Z)$ and a point $P \in \operatorname{Col}(E)$. Consider a cover $\{E_1,\ldots, E_u\}$ of $E$ at $P$. By Theorem~\ref{th-multiplicidad-jac}, we have that
\begin{align*}
  \nu_P(\pi_E^* {\mathcal J}_{\mathcal{F},\mathcal{G}}) & = \nu_P(\pi_E^* C) + \nu_P(\pi_E^* D) + t_P  \\
  \sum_{l=1}^{u} \sum_{Q \in E_{l,red}} \nu_Q(\pi_{E_l}^* {\mathcal J}_{\mathcal{F},\mathcal{G}}) & =  \sum_{l=1}^{u} \sum_{Q \in E_{l,red}} (\nu_Q(\pi_{E_l}^*C) + \nu_Q(\pi_{E_l}^*D) + \tau_{E_l}(Q))
\end{align*}
By the properties of a cover given in definition~\ref{def:cover}, we have that
$$\nu_P(\pi_E^* C) + \nu_P(\pi_E^* D) =\sum_{l=1}^{u} \sum_{Q \in E_{l,red}} (\nu_Q(\pi_{E_l}^*C) + \nu_Q(\pi_{E_l}^*D))$$
and the result is straightforward.
\end{proof}

The results above allow to give a decomposition of ${\mathcal J}_{\mathcal{F},\mathcal{G}}$  into bunches of branches in the sense of the decomposition theorem of polar curves.
Recall that given a divisor $E$ of $\pi_Z^{-1}(0)$, we denote by $\pi_E : X_E \to ({\mathbb C}^2,0)$ the morphism reduction of $\pi_Z$ to $E$  and we write $\pi_Z= \pi_E \circ \pi_E'$.
Let $B(Z)$ be the set of bifurcation divisors of $G(Z)$. Given any $E \in B(Z)$ which is a non-collinear divisor for $\mathcal F$ and $\mathcal G$, we define $J^E_{nc}$ as the union of the branches $\xi$ of ${\mathcal J}_{\mathcal{F},\mathcal{G}}$ such that
\begin{itemize}
  \item $\pi_E^* \xi \cap \pi_E^* Z = \emptyset$;
  \item if $E' < E$, then $\pi_E^* \xi \cap \pi'_E(E')=\emptyset$
  \item if $E'> E$, then $\pi_{E'}^* \xi \cap E'_{red}=\emptyset$
\end{itemize}
Moreover, given a non-collinear divisor $E$, we denote $J_c^E=\cup_{P \in \text{Col}(E)} J_P^E$ (with $J_c^E=\emptyset$ if $\text{Col}(E)=\emptyset$).

Thus, the previous results allow us to give a decomposition of $${\mathcal J}_{\mathcal{F},\mathcal{G}}= J^* \cup (\cup_{E \in B_N(Z)} J^E )$$ (see below for the precise statement) such there is a certain control of the topology of the irreducible components of $\mathcal{J}_{\mathcal{F},\mathcal{G}}$ obtained from the data of the foliations $\mathcal{F}$ and $\mathcal{G}$ provided that the component of $\mathcal{J}_{\mathcal{F},\mathcal{G}}$ is attached either to a non-collinear   divisor or to  a chain of  collinear divisors which are in between two non-collinear bifurcation divisors.
The irreducible components corresponding to $J^*$ are the one attached to ``isolated'' collinear divisors for which no control is possible.

Given a non-collinear bifurcation divisor $E$ of $G(Z)$, we  denote
$$t^*(E)=\sum_{Q \in M(E) \smallsetminus \operatorname{Col}(E)} t_Q,$$
that is, the number of zeros of ${\mathcal M}_E(z)$  (counting with multiplicities) which do not correspond to collinear points.
Then, we can state the properties of the decomposition of ${\mathcal J}_{\mathcal{F},\mathcal{G}}$ as follows:
\begin{Theorem}\label{th:descomposicion-no-sing}
Consider ${\mathcal F} \in {\mathbb G}_C$ and ${\mathcal G} \in {\mathbb G}_D$ such that $Z=C \cup D$ is a curve with only non-singular irreducible components. Let
 $B_N(Z)$ be the set of non-collinear bifurcation divisors of $G(Z)$. Then there is a unique decomposition ${\mathcal J}_{\mathcal{F},\mathcal{G}}= J^* \cup (\cup_{E \in B_N(Z)} J^E )$ where $J^E=J^E_{nc} \cup J_{c}^E$  with the following properties
\begin{itemize}
  \item[(1)] $ \nu_0(J^E_{nc}) =t^*(E)$. In particular,  $\nu_0(J^E_{nc}) \leq \sharp \operatorname{NCol}(E)-1 \leq b_E-1$.
  \item[(2)] $\pi_E^* J^E_{nc} \cap \pi_E^* Z=\emptyset$
  \item[(3)] if $E' < E$, then $\pi_E^* J^E_{nc} \cap \pi_E'(E')=\emptyset$
  \item[(4)] if $E' > E$, then $\pi_{E'}^* J^E_{nc} \cap E_{red}' = \emptyset$
  \item[(5)] if $\delta$ is a branch of $J_{c}^E$, then $\pi_E^* \delta \cap E_{red}$ is a point in $\operatorname{Col}(E)$.
   \item[(6)] $\nu_0(J_c^E)=\sum_{P \in C(E)} ( t_P + \sum_{l=1}^{u(P)} (\sharp \operatorname{NCol}(E_l^P) - t(E_l^P)))$ where $\{E_1^P, \ldots, E_{u(P)}^P\}$ is a cover of $E$ at $P$.
\end{itemize}
Moreover, if $E$ is a purely non-collinear divisor with $\sum_{R_l^E \in \operatorname{NCol}(E)} \Delta_{E}(R_l^E) \neq 0$, then
\begin{equation}\label{eq:mult-JE-maxima}
\nu_0(J^E)=\nu_0(J^E_{nc})=b_E-1.
\end{equation}
\end{Theorem}
\begin{proof}
We have that
\begin{align*}
\nu_0(J_{nc}^E) & =\sum_{P \in M(E) \smallsetminus \operatorname{Col}(E)} \nu_P(\pi_E^* {\mathcal J}_{\mathcal{F},\mathcal{G}} ) =\sum_{P \in M(E) \smallsetminus \operatorname{Col}(E)} t_P =t^*(E) \\
&  \leq \sum_{ P \in M(E)} t_P \leq \sharp \operatorname{NCol}(E)-1 \leq b_E -1 \\
\end{align*}
where we have used the inequality given in~\eqref{eq:no-colineal} and the fact that $\sharp \operatorname{NCol}(E) \leq b_E$. This gives the first statement of the theorem.

Moreover, if $E$ is a purely non-collinear divisor, then $\text{Col}(E)=\emptyset$, $J^E=J^E_{nc}$ and $\sharp \operatorname{NCol}(E)=b_E$. In addition, when $\sum_{R_l^E \in \operatorname{NCol}(E)} \Delta_{E}(R_l^E) \neq 0$ we have that $\sum_{P \in M(E)} t_P=\sharp \operatorname{NCol}(E)-1$ by Lemma~\ref{lema:E-no-colineal}.  Consequently, we deduce that
  $$\nu_0(J^E)=\sum_{P \in M(E)} \nu_P(\pi_E^* {\mathcal J}_{\mathcal{F},\mathcal{G}} ) =\sum_{ P \in M(E)} t_P = \sharp \operatorname{NCol}(E)-1 = b_E -1
$$
and we obtain expression~\eqref{eq:mult-JE-maxima}.

Properties (2), (3) and (4) are consequence of the definition of $J_{nc}^E$. Properties (5) and (6) follow directly from the definition of $J_c^E$ and Theorem~\ref{th:colinear-point}.
\end{proof}
Note that the properties of $J_{nc}^E$ can be stated in terms of coincidences as follows: if $\delta$ is an irreducible component of $J_{nc}^E$ (with $E$ a non-collinear bifurcation divisor) and $Z_i$ is an irreducible component of $Z=C \cup D$, then
$$
{\mathcal C}(\delta,Z_i)= \left\{
                            \begin{array}{ll}
                              v(E), & \hbox{ if } E \text{ is in the geodesic of } Z_i, \\
                              {\mathcal C}(\gamma_E,Z_i), & \hbox{ otherwise.}
                            \end{array}
                          \right.
$$
where $\gamma_E$ is any $E$-curvette. Observe that $v(E)={\mathcal C}(\gamma_E,Z_i)$ when $E$ is in the geodesic of $Z_i$ and $\gamma_E$   does not intersect $E$ at the points $\pi_E^* Z \cap E$.

Next result determines the intersection multiplicity of $J^E_{nc}$ with the curves of separatrices $C$ and $D$ of the foliations $\mathcal F$ and $\mathcal G$.
\begin{Corollary}\label{cor:mult-int-J-C-D}
If $E$ is a non-collinear bifurcation divisor, then
$$(J_{nc}^E,C)_0= \nu_E(C) \cdot t^*(E); \qquad (J_{nc}^E,D)_0= \nu_E(D) \cdot t^*(E)
$$
where $\nu_E(C)=(C,\gamma_E)_0$ and $\nu_E(D)=(D,\gamma_E)_0$ with $\gamma_E$ any $E$-curvette.
\end{Corollary}
\begin{proof} %[Proof of Corollary~\ref{cor:mult-int-J-C-D}]
Let $E$ be a non-collinear bifurcation divisor of $G(Z)$ and let $\gamma_E$ be any $E$-curvette which does not intersect $E$ at the points $\pi_E^*Z \cap E$.
By the properties of $J_{nc}^E$ given in Theorem~\ref{th:descomposicion-no-sing}, we have that if $\delta$ is a branch of $J_{nc}^E$ then
$${\mathcal C}(\delta,C_i)=\left\{
                             \begin{array}{ll}
                               v(E), & \hbox{if }  E \text{ is in the geodesic of } C_i; \\
                               {\mathcal C}(\gamma_E,C_i) , & \hbox{othewise.}
                             \end{array}
                           \right.
$$
Note that ${\mathcal C}(\gamma_E,C_i)= v(E)$ if $E$ is in the geodesic of $C_i$ (that is, $i \in I_E$). Moreover, since $\gamma_E$ and $C_i$ are non-singular curves, we have that ${\mathcal C}(\gamma_E,C_i)=(\gamma_E,C_i)_0$. Therefore, using the relationship between the coincidence and the intersection multiplicity of two branches given in Remark~\ref{rmk:coincidencia-mult-int}, we have that
$$(\delta,C_i)_0=\nu_0(\delta) \cdot (\gamma_E,C_i)_0$$
for a branch $\delta$ of $J_{nc}^E$.
Now, if we denote by ${\mathcal B}(J_{nc}^E)$ the set of branches of $J_{nc}^E$, we have that
\begin{align*}
  (J_{nc}^E,C)_0 & = \sum_{i=1}^{r} (J_{nc}^E,C_i)_0 =  \sum_{i=1}^{r} \sum_{\delta \in {\mathcal B}(J_{nc}^{E})} (\delta,C_i)_0 \\
   & = \sum_{i=1}^{r} \sum_{\delta \in {\mathcal B}(J_{nc}^{E})} \nu_0 (\delta) \cdot (\gamma_E,C_i)_0 = \nu_0(J_{nc}^E) \sum_{i=1}^{r}  (\gamma_E,C_i)_0 \\
 & = t^*(E) \cdot \nu_E(C)
\end{align*}
\end{proof}

As a consequence of the result above and Propositions \ref{prop:int-sep} and \ref{prop:milnor-number} we obtain next corollary for non-dicritical generalized curve foliations which relates invariants of the foliations $\mathcal F$ and $\mathcal G$, such as the Milnor numbers o the tangency orders, with data coming from the decomposition of the jacobian curve.
\begin{Corollary}\label{cor:sum-mult-int}
With the hypothesis and notations of Theorem~\ref{th:descomposicion-no-sing}, we get that
$$
\sum_{E \in B_N(Z)} \nu_0(J_{nc}^E) \nu_E(C_i) \leq \mu_0({\mathcal F},C_i) + \tau_0({\mathcal G},C_i)
$$
and
$$
\sum_{E \in B_N(Z)} \nu_0(J_{nc}^E) ( \nu_E(C) - \nu_E(D)) \leq \mu_0({\mathcal F}) -  \mu_0({\mathcal G}).
$$
\end{Corollary}

\begin{proof}%[Proof of Corollary~\ref{cor:sum-mult-int}]
We have just proved that $(J_{nc}^E,C_i)_0=\nu_0(J_{nc}^E) \nu_E(C_i)$. Thus, using Proposition~\ref{prop:int-sep}, we get that
$$
\sum_{E \in B_N(Z)} (J_{nc}^E,C_i)_0 = \sum_{E \in B_N(Z)} \nu_0(J_{nc}^E) \nu_E(C_i) \leq ({\mathcal J}_{\mathcal{F},\mathcal{G}},C_i)_0 = \mu_0({\mathcal F},C_i) + \tau_0({\mathcal G},C_i).
$$
Now, from Corollary~\ref{cor:mult-int-J-C-D} and Proposition~\ref{prop:milnor-number}, we obtain that
\begin{align*}
  \sum_{E \in B_N(Z)} ((J_{nc}^E,C)_0 - (J_{nc}^E,D)_0)) & = \sum_{E \in B_N(Z)} \nu_0(J_{nc}^E) (\nu_{E}(C) - \nu_E(D)) \\
  & \leq ({\mathcal J}_{\mathcal{F},\mathcal{G}},C)_0 -({\mathcal J}_{\mathcal{F},\mathcal{G}},D)_0= \mu_0({\mathcal F}) - \mu_0({\mathcal G})
\end{align*}
which gives the second inequality.
\end{proof}

The general case of foliations with separatrices that can have singular irreducible components will be treated in next section. % section~\ref{sec:caso-general}.
\section{General case}\label{sec:caso-general}
Consider two plane curves $C=\cup_{i=1}^r C_i$ and $D=\cup_{j=1}^sD_j$ which can have singular branches. Assume that $C$ and $D$ have no common irreducible components. Let $\rho: ({\mathbb C}^2,0) \to ({\mathbb C}^2,0)$ be a ramification given in coordinates by $\rho(u,v)=(u^n,v)$ such that the curve $\rho^{-1}Z$ has only non-singular irreducible components where $Z=C \cup D$. In this section we will denote $\widetilde{B}$ the curve $\rho^{-1} B$ for any plane curve $B$. See appendix~\ref{ap:ramificacion} for notations concerning ramifications.

Take ${\mathcal F}$ and ${\mathcal G}$ foliations with $C$ and $D$ as curve of separatrices respectively. Let us study the relationship between the curves $\widetilde{\mathcal J}_{{\mathcal F},\mathcal{G}}=\rho^{-1} {\mathcal J}_{{\mathcal F},\mathcal{G}}$ and ${\mathcal J}_{\rho^* \mathcal{F},\rho^* \mathcal{G}}$.

Assume that the foliations $\mathcal F$ and $\mathcal G$ are given by $\omega=0$ and $\eta=0$ with
$$\omega=A(x,y) dx + B(x,y) dy; \quad \eta=P(x,y) dx + Q(x,y) dy,$$
then $\rho^*{\mathcal F}$ and $\rho^*{\mathcal G}$ are given by $\rho^* \omega=0$ and $\rho^* \eta=0$ where
$$\rho^* \omega= A(u^n,v) n u^{n-1} du+B(u^n,v) dv; \quad \rho^* \eta= P(u^n,v) n u^{n-1} du+Q(u^n,v) dv.$$
Therefore, if we write $J(x,y)=A(x,y) Q(x,y) - B(x,y) P(x,y)$, then the curve $\rho^{-1} {\mathcal J}_{{\mathcal F},\mathcal{G}}$ is given by
$J(u^n,v)=0$ whereas ${\mathcal J}_{\rho^* \mathcal{F},\rho^* \mathcal{G}}$ is given by $nu^{n-1}J(u^n,v)=0$. Let us see (Corollary~\ref{cor:mult-jac-ram}) that $\rho^{-1} {\mathcal J}_{{\mathcal F},\mathcal{G}}=\widetilde{\mathcal J}_{{\mathcal F},\mathcal{G}}$ satisfies the statements of Theorem~\ref{th-multiplicidad-jac} with respect to $\rho^{-1} Z=\widetilde{Z}$.

Let  $\pi_{\widetilde{Z}}: \widetilde{X} \to ({\mathbb C}^2,0)$ be the minimal reduction of singularities of $\widetilde{Z}$. We denote by $\widetilde{E}$ any irreducible component of $\pi_{\widetilde{Z}}^{-1}(0)$ and by $\pi_{\widetilde{E}} :\widetilde{X}_{\tilde{E}} \to ({\mathbb C}^2,0)$ the morphism reduction of $\pi_{\widetilde{Z}}$ to $\widetilde{E}$. Let us state some properties concerning the infinitely near points of $\widetilde{\mathcal J}_{\mathcal{F},\mathcal{G}}$ and ${\mathcal J}_{\rho^* \mathcal{F},\rho^* \mathcal{G}}$:

\begin{Lemma}\label{lemma:ptos-inf-prox-ram}
Let  $\widetilde{E}$ be an irreducible component of $\pi_{\widetilde{Z}}^{-1}(0)$.  %such that $\widetilde{E}$ is a non-collinear divisor.
We have that
$$
\pi_{\widetilde{E}}^* \widetilde{\mathcal J}_{{\mathcal F},{\mathcal G}} \cap \widetilde{E}_{red}^* = \pi_{\widetilde{E}}^* {\mathcal J}_{\rho^* \mathcal{F},\rho^* \mathcal{G}}  \cap \widetilde{E}_{red}^*,
$$
where $\widetilde{E}_{red}^*$ denote the points in the first chart of $\widetilde{E}_{red}$.
Moreover, $$\nu_P ( \pi_{\widetilde{E}}^* \widetilde{\mathcal J}_{{\mathcal F},{\mathcal G}}) = \nu_P( \pi_{\widetilde{E}}^* {\mathcal J}_{\rho^* \mathcal{F},\rho^* \mathcal{G}})$$ for each $P  \in \pi_{\widetilde{E}}^* \widetilde{\mathcal J}_{{\mathcal F},{\mathcal G}} \cap \widetilde{E}_{red}^*$.
\end{Lemma}
\begin{proof}
Take $\widetilde{E}$ an irreducible component of $\pi_{\widetilde{Z}}^{-1}(0)$ with $v(\widetilde{E})=p$ and assume that $(u,v)$ are coordinates adapted to $\widetilde{E}$. If we denote $\widetilde{J}(u,v)=J(u^n,v)$, we have that
$${\rm In}_p (nu^{n-1} \widetilde{J};u,v)=nu^{n-1} {\rm In}_p (\widetilde{J};u,v).$$
%Let us denote $\widetilde{J}_I(u,v)={\rm In}_p (\widetilde{J};u,v)$.
Write ${\rm In}_p (\widetilde{J};u,v)=\sum_{i+pj=k} h_{ij}u^iv^j$.
Hence, if $(u_p,v_p)$ are coordinates in the first chart of $\widetilde{E}_{red} \subset \widetilde{X}_{\widetilde{E}}$ such that $\pi_{\widetilde{E}} (u_p,v_p)=(u_p,u_p^p v_p)$ and $\widetilde{E}_{red}=(u_p=0)$, then the points $\pi_{\widetilde{E}}^* \widetilde{\mathcal J}_{{\mathcal F},{\mathcal G}} \cap \widetilde{E}_{red}$, in the first chart of $\widetilde{E}_{red}$, are given by $u_p=0$ and
$\sum_{i+pj=k} h_{ij}v_p^j=0$. %$\widetilde{J}_I(1,v_p)=0$.
This proves that  $\pi_{\widetilde{E}}^* \widetilde{\mathcal J}_{{\mathcal F},{\mathcal G}} \cap \widetilde{E}_{red}^* = \pi_{\widetilde{E}}^* {\mathcal J}_{\rho^* \mathcal{F},\rho^* \mathcal{G}}  \cap \widetilde{E}_{red}^*
$ and that $\nu_P ( \pi_{\widetilde{E}}^* \widetilde{\mathcal J}_{{\mathcal F},{\mathcal G}}) = \nu_P( \pi_{\widetilde{E}}^* {\mathcal J}_{\rho^* \mathcal{F},\rho^* \mathcal{G}})$ for each $P  \in \pi_{\widetilde{E}}^* \widetilde{\mathcal J}_{{\mathcal F},{\mathcal G}} \cap \widetilde{E}_{red}^*$.
\end{proof}
Hence, when $\widetilde{E}$ is a non-collinear divisor for the foliations $\rho^* {\mathcal F}$ and $\rho^* {\mathcal G}$, the curve ${\mathcal J}_{\rho^* \mathcal{F},\rho^* \mathcal{G}}$ satisfies Theorem~\ref{th-multiplicidad-jac} with respect to $\rho^{-1} Z=\widetilde{Z}$, and thanks to the previous lemma, we get the following result for $\rho^{-1} {\mathcal J}_{{\mathcal F},\mathcal{G}}=\widetilde{\mathcal J}_{{\mathcal F},\mathcal{G}}$:
\begin{Corollary}\label{cor:mult-jac-ram}
Take $\widetilde{E}$ an irreducible component  of $\pi^{-1}_{\widetilde{Z}}(0)$ which   is a non-collinear divisor for the foliations $\rho^* {\mathcal F}$ and $\rho^* {\mathcal G}$. Given any $P \in \widetilde{E}_{red}^*$, we have that
$$
\nu_P(\pi_{\widetilde{E}}^*\widetilde{\mathcal J}_{\mathcal{F},\mathcal{G}})= \nu_P(\pi_{\widetilde{E}}^* \widetilde{C}) + \nu_P(\pi_{\widetilde{E}}^*\widetilde{D}) + \tau_{\widetilde{E}}(P).
$$
In particular, if $P \in \widetilde{E}_{red}^*$ with $\nu_P(\pi_{\widetilde{E}}^*{\widetilde{\mathcal J}}_{\mathcal{F},\mathcal{G}})>0$, then  $P$ is an infinitely near point of $\widetilde{Z}$ or a point in $M(\widetilde{E})$.
\end{Corollary}
Let $E$ be a bifurcation divisor of $G(Z)$ and consider $\widetilde{E}_l, \widetilde{E}_{k}$ two bifurcation divisors of $G(\widetilde{Z})$ associated to $E$. Recall that
there is a bijection between the sets of points $\pi_{\widetilde{E}^l}^* \widetilde{Z} \cap \widetilde{E}^l_{red}$ and $\pi_{\widetilde{E}^k}^* \widetilde{Z} \cap \widetilde{E}^k_{red}$ given by the map $\rho_{l,k}: \widetilde{E}^l_{red} \to \widetilde{E}^k_{red}$ (see appendix~\ref{ap:ramificacion}). Thus we will denote  $\{R_1^{\widetilde{E}^l}, R_2^{\widetilde{E}^l}, \cdots, R_{b_{\widetilde{E}^l}}^{\widetilde{E}^l}\}$ and $\{R_1^{\widetilde{E}^k}, R_2^{\widetilde{E}^k}, \cdots, R_{b_{\widetilde{E}^k}}^{\widetilde{E}^k}\}$ the sets of points $\pi_{\widetilde{E}^l}^* \widetilde{Z} \cap \widetilde{E}^l_{red}$ and $\pi_{\widetilde{E}^k}^* \widetilde{Z} \cap \widetilde{E}^k_{red}$ respectively, with $R_t^{\widetilde{E}^k}= \rho_{l,k}(R_t^{\widetilde{E}^l})$ for $t=1,2,\ldots, b_{\widetilde{E}^k}$.
 By the results in appendix~\ref{subsec:ram-log} (see Proposition~\ref{prop:ram-modelo-log} and equation~\eqref{eq:indices-ramificados}), we get that
\begin{align*}
  \mathcal{I}_{R_t^{\widetilde{E}^l}}(\pi_{\widetilde{E}^l}^* \rho^*{\mathcal F},\widetilde{E}^l_{red}) & =  \mathcal{I}_{R_t^{\widetilde{E}^k}}(\pi_{\widetilde{E}^k}^* \rho^*{\mathcal F},\widetilde{E}^k_{red}) \\
   \mathcal{I}_{R_t^{\widetilde{E}^l}}(\pi_{\widetilde{E}^l}^* \rho^*{\mathcal G},\widetilde{E}^l_{red}) & =  \mathcal{I}_{R_t^{\widetilde{E}^k}}(\pi_{\widetilde{E}^k}^* \rho^*{\mathcal G},\widetilde{E}^k_{red})
\end{align*}
which implies
$$
\Delta_{\widetilde{E}^l}(R_t^{\widetilde{E}^l})=\Delta_{\widetilde{E}^k}(R_t^{\widetilde{E}^k}) \quad \text{ for } \quad t=1,2,\ldots, b_{\widetilde{E}^l},
$$
with  $\Delta_{\widetilde{E}^l}(R_t^{\widetilde{E}^l})=\Delta_{\widetilde{E}^l}^{\rho^* \mathcal{F},\rho^*{\mathcal G}}(R_t^{\widetilde{E}^l})$. Thus, $\widetilde{E}^l$ is collinear (resp. non-collinear) if and only if $\widetilde{E}^k$ is also collinear (resp. non-collinear). So we can introduce the following definition
\begin{Definition}
We say that a bifurcation divisor $E$ of $G(Z)$ is {\em collinear (resp. non-collinear)\/} for the foliations $\mathcal F$ and $\mathcal G$ when any of its associated divisors $\widetilde{E}^l$ is collinear (resp. non-collinear) for the foliations $\rho^{\ast} {\mathcal F}$ and $\rho^{\ast} {\mathcal G}$.
\end{Definition}
Moreover, if $R_t^{\widetilde{E}^l}, R_s^{\widetilde{E}^l}$ are two points  $\pi_{\widetilde{E}^l}^* \widetilde{Z} \cap \widetilde{E}^l_{red}$ with $\rho_{\widetilde{E}^l,E}(R_t^{\widetilde{E}^l}) =\rho_{\widetilde{E}^l,E}(R_s^{\widetilde{E}^l})$ where $\rho_{\widetilde{E}^l,E}: \widetilde{E}^l_{red} \to E_{red}$ is the ramification defined in appendix~\ref{ap:ramificacion}, then
$$
\Delta_{\widetilde{E}^l}(R_t^{\widetilde{E}^l})=\Delta_{\widetilde{E}^k}(R_s^{\widetilde{E}^l})
$$
by equation~\eqref{eq:indices-ramificados-mismo-divisor} in appendix~\ref{subsec:ram-log}. Thus, we say that an infinitely near point $R^E$ of $Z$ in $E_{red}$ is a {\em collinear point} (resp. {\em non-collinear point}) for the foliations $\mathcal F$ and $\mathcal G$ if for any associated divisor $\widetilde{E}^l$ and any infinitely near point $R_t^{\widetilde{E}^l}$ of $\rho^{-1} Z$ in $\widetilde{E}^l_{red}$  with $\rho_{\widetilde{E}^l,E}(R_t^{\widetilde{E}^l}) =R^E$, the point $R_t^{\widetilde{E}^l}$ is collinear (resp. non-collinear) for the foliations $\rho^* {\mathcal F}$ and $\rho^*{\mathcal G}$. Given a bifurcation divisor $E$ of $G(Z)$, we denote by $\operatorname{Col}(E)$ the set of collinear points of $E$ and $\operatorname{NCol}(E)$ the set of non-collinear points of $E$.

Corollary~\ref{cor:mult-jac-ram} and the results in section~\ref{sec:jacobian} allow to give a decomposition of ${\mathcal J}_{\mathcal{F},\mathcal{G}}$. By Theorem~\ref{th:descomposicion-no-sing}, we have a decomposition
$$
\widetilde{\mathcal J}_{\mathcal{F},\mathcal{G}} = \widetilde{J}^* \cup (\cup_{\widetilde{E} \in B_{N}(\widetilde{Z})} J^{\widetilde{E}})
$$
with $J^{\widetilde{E}}=J^{\widetilde{E}}_{nc} \cup J^{\widetilde{E}}_{c}$. Given a non-collinear bifurcation divisor $E$ of $G(Z)$, we define $J^E=J_{nc}^E \cup J_c^E$ to be such that
$$
\rho^{-1}{J^E_{nc}}=\bigcup_{l=1}^{\underline{n}_E}  J^{\widetilde{E}^l}_{nc}; \qquad \rho^{-1}{J^E_{c}}=\bigcup_{l=1}^{\underline{n}_E}  J^{\widetilde{E}^l}_{c}
$$
where $\{\widetilde{E}^l\}_{l=1}^{\underline{n}_E}$ are the divisors of $G(\widetilde{Z})$ associated to $E$ and $J^*$ to be such that $\rho^{-1} J^* = \widetilde{J}^*$.
Hence, we can state the main result of this paper
\begin{Theorem}\label{th:desc-general}
Let us write $Z=\cup_{i=1}^{r+s} Z_i$  with $Z_i$ irreducible and denote by $B_N(Z)$  the set of non-collinear bifurcation divisors of $G(Z)$. Then there is a decomposition
$$
{\mathcal J}_{\mathcal{F},\mathcal{G}} = J^* \cup (\cup_{E \in B_{N}(Z)} J^{E})
$$
with $J^E=J^E_{nc} \cup J^E_c$ such that
\begin{itemize}
  \item[(i)] $\nu_0(J^E_{nc}) \leq \left\{
                                   \begin{array}{ll}
                                     \underline{n}_E n_E (b_E-1), & \hbox{if $E$ does not belong to a dead arc;} \\
                                     \underline{n}_E n_E (b_E-1) - \underline{n}_E, & \hbox{otherwise.}
                                   \end{array}
                                 \right.
$
  \item[(ii)] For each irreducible component $\delta$ of $J_{nc}^E$ we have that
\begin{itemize}
  \item ${\mathcal C}(\delta,Z_i)=v(E)$ if $E$ belongs to the geodesic of $Z_i$;
  \item ${\mathcal C}(\delta,Z_j)={\mathcal C}(Z_i,Z_j)$ if $E$ belongs to the geodesic of $Z_i$ but not to the one of $Z_j$.
\end{itemize}
    \item[(iii)] For each irreducible component $\delta$ of $J_{c}^E$, there exists an irreducible component $Z_i$ of $Z$ such that $E$ belongs to its geodesic and
$${\mathcal C}(\delta,Z_i)>v(E).$$
Moreover, if $E'$ is the first non-collinear bifurcation divisor in the geodesic of $Z_i$ after $E$, then $${\mathcal C}(\delta,Z_i)<v(E').$$
\end{itemize}
\end{Theorem}

\section{Jacobian curves of hamiltonian foliations and Polar curves of foliations}\label{sec:polares}
In this section we will explain how our results imply previous results concerning jacobian curves of two plane curves or polar curves of foliations.
\subsection{Jacobian of two curves} \label{sec:curvas}
In \cite{Kuo-P-2002,Kuo-P-2004}, T.-C. Kuo and A. Parusi\'{n}ki consider the Jacobian $f_x g_y-f_y g_x=0$ of a pair of germs of holomorphic functions $f, g$ without common branches and give properties of its Puiseux series which they called {\em polar roots\/} of the Jacobian. They define a tree-model, noted $T(f,g)$, which represents the Puiseux series of the curves $C=(f=0)$ and $D=(g=0)$ and the contact orders among these series.
The tree-model $T(f,g)$ is constructed as follows:
it starts with an horizontal bar $B_*$ called {\em ground bar\/} and a vertical segment on $B_*$ called the {\em main trunk} of the tree. This trunk is marked with $[p,q]$ where $p=\nu_0(C)$ and $q=\nu_0(D)$.  %$C=(f=0)$ and $D=(g=0)$.
Let $\{y_i^C(x)\}_{i=1}^{\nu_0(C)}$ and $\{y_i^D(x)\}_{i=1}^{\nu_0(D)}$ be the Puiseux series of $C$ and $D$ respectively, and denote $\{z_j(x)\}_{j=1}^{N}$ the set $\{y_i^C(x)\}_{i=1}^{\nu_0(C)} \cup \{y_i^D(x)\}_{i=1}^{\nu_0(D)}$ with $N=\nu_0(C)+\nu_0(D)$. Denote by  $h_0=\min \{ \text{ord}_x(z_i(x)-z_j(x)) \ : \ 1 \leq i, j \leq N\}$. A bar $B_0$ is drawn on top of the main trunk with $h(B_0)=h_0$ being the {\em height} of $B_0$. The Puiseux series $\{z_j(x)\}$ are divided into equivalence classes (mod $B_0$) by the following relation: $z_j(x) \sim_{B_0} z_k(x)$ if $\text{ord}_x(z_j(x)-z_k(x))> h_0$. Each equivalence class is represented by a vertical line, called {\em trunk}, drawn on the top of $B_0$. Each trunk is marked by a {\em bimultiplicity} $[s,t]$ where $s$ (resp. $t$) denote the number of Puiseux series   of $C$ (resp. of $D$) in the equivalence class. The same construction is repeated recursively on each trunk. The construction finishes with trunks which have bimultiplicity $[1,0]$ or $[0,1]$ representing each Puiseux series of the curve $Z=C \cup D$.

Let us now consider the curve $\widetilde{Z}=\rho^{-1}Z$ where $\rho:({\mathbb C}^2,0) \to ({\mathbb C}^2,0)$ is any $Z$-ramification given by $\rho(u,v)=(u^n,v)$  (see appendix~\ref{ap:ramificacion}).
Since the branches of $\widetilde{Z}$ are in bijection with the Puiseux series of $Z$ and the valuations of the  bifurcation divisors of $G(\widetilde{Z})$ represent the contact orders among these series, then the tree model above $T(f,g)$ can be recovered from the dual graph of $G(\widetilde{Z})$:
there is a bijection between the set of bars of $T(f,g)$ which are not the ground bar and the bifurcation divisors $B(\widetilde{Z})$ of $G(\widetilde{Z})$. For instance, the first bar $B_0$ corresponds to the first bifurcation divisor $\widetilde{E}_1$ of $G(\widetilde{Z})$ and $h(B_0)=\frac{v(\widetilde{E}_1)}{n}$. The number of trunks on $B_0$ is equal to $b_{\widetilde{E}_1}$, that is, each trunk on $B_0$ correspond to an infinitely near point of $\widetilde{Z}$ on $\widetilde{E}_{1,red}$. In particular, given a trunk with bimultiplicity $[s,t]$ corresponding to a point $R_{i}^{\widetilde{E}_1}$, then $s=\nu_{R_i^{\widetilde{E}_1}} (\pi_{\widetilde{E}_1}^* C)$ and $t=\nu_{R_i^{\widetilde{E}_1}} (\pi_{\widetilde{E}_1}^* D)$.

We shall illustrate with an example the relationship between $T(f,g)$ and $G(\widetilde{Z})$. The following example  corresponds to the Example 1.1 in \cite{Kuo-P-2004}.
\begin{Example} Take positive integers $d,f$ with $d < f$ and non-zero constants $A,B$. Consider
\begin{align*}
f(x,y) & =(y+x) (y-x^{d+1}+Ax^{f+1}) (y+x^{d+1}+Bx^{f+1}) \\
g(x,y) & =(y-x) (y-x^{d+1}-Ax^{f+1}) (y+x^{d+1}-Bx^{f+1}) \\
\end{align*}
and put $C=(f=0)$ and $D=(g=0)$. Since $Z=C \cup D$ has only non-singular irreducible components, we do not need to consider a ramification. The dual graph $G(Z)$ is given by
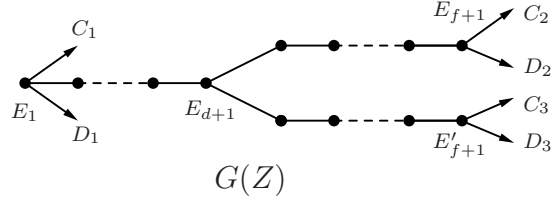
\begin{center}
\begin{texdraw}
\arrowheadtype t:F \arrowheadsize l:0.08 w:0.04 \drawdim mm \setgray
0

\move(10 20)
%\htext{$G(C)$}
\fcir f:0 r:0.8 \avec (17 15) \move (10 20)  \avec (17 25)  \move (10 20)   \rlvec (7 0)  \fcir f:0 r:0.8 \lpatt(1 1)  \rlvec (10 0)  \fcir f:0 r:0.8 \lpatt () \rlvec (7 0) \fcir
f:0 r:0.8  \rlvec (10 5)  \fcir f:0 r:0.8 \rlvec (7 0) \fcir f:0 r:0.8   \lpatt(1 1) \rlvec (10 0)  \fcir f:0 r:0.8 \lpatt ()   \rlvec (7 0)  \fcir
f:0 r:0.8 \avec (75 30) \move(68 25) \avec (75 22)

\move (34 20)  \rlvec (10 -5) \fcir f:0 r:0.8   \lpatt () \fcir
f:0 r:0.8 \rlvec (7 0) \fcir f:0 r:0.8  \lpatt(1 1)  \rlvec (10 0) \lpatt() \fcir f:0 r:0.8  \rlvec (7 0)  \fcir
f:0 r:0.8   \avec (75 18) \move(68 15) \avec (75 12)

%\fcir f:0 r:0.8 \rlvec (10 -5)  \fcir f:0 r:0.8 \rlvec (10 0)
%\avec (75 23) \move (70 20) \avec (75 17)

\move(8 15) \htext{\tiny{$E_1$}} \move(31 15)
\htext{\tiny{$E_{d+1}$}}
\move(64 28) \htext{\tiny{$E_{f+1}$}}
\move(64 10) \htext{\tiny{$E_{f+1}'$}}

\move (16 26)  \htext{\tiny{$C_1$}}
\move (16 12)  \htext{\tiny{$D_1$}}
\move (76 28)  \htext{\tiny{$C_2$}}
\move (76 21)  \htext{\tiny{$D_2$}}
\move (76 16)  \htext{\tiny{$C_3$}}
\move (76 11)  \htext{\tiny{$D_3$}}

\move (35 5) \htext{$G(Z)$} %\move (60 5) \htext{$G(\rho^{-1}C)$}
\end{texdraw}
\end{center}
while the tree-model is given by
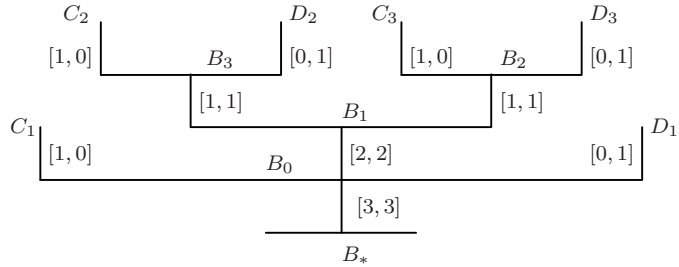
\vspace{\baselineskip}
\begin{center}
\begin{texdraw}
\arrowheadtype t:F \arrowheadsize l:0.08 w:0.04 \drawdim mm \setgray
0
\move (40 20)
\rlvec (20 0) \move (50 20) \rlvec (0 7)
\move (10 27) \rlvec (80 0) \rlvec (0 7)
\move (10 27) \rlvec (0 7) \move (50 27) \rlvec (0 7)
\move (30 34) \rlvec (40 0) \rlvec (0 7) \move (30 34) \rlvec (0 7)
\move (18 41) \rlvec (24 0) \rlvec (0 7) \move (18 41) \rlvec (0 7)
\move (58 41) \rlvec (24 0) \rlvec (0 7) \move (58 41) \rlvec (0 7)

\move (50 16) \htext{\tiny{$B_*$}}
\move (52 22)  \htext{\tiny{$[3,3]$}}
\move (40 28) \htext{\tiny{$B_0$}}
\move (51 29)  \htext{\tiny{$[2,2]$}}  \move (11 29)  \htext{\tiny{$[1,0]$}} \move (83 29)  \htext{\tiny{$[0,1]$}}
\move (50 35) \htext{\tiny{$B_1$}}
\move (31 36)  \htext{\tiny{$[1,1]$}} \move (71 36)  \htext{\tiny{$[1,1]$}}
\move (32 42) \htext{\tiny{$B_3$}} \move (11 42)  \htext{\tiny{$[1,0]$}}  \move (43 42)  \htext{\tiny{$[0,1]$}}
\move (71 42) \htext{\tiny{$B_2$}} \move (59 42)  \htext{\tiny{$[1,0]$}}  \move (83 42)  \htext{\tiny{$[0,1]$}}

\move (6 33)  \htext{\tiny{$C_1$}} \move (91 33)  \htext{\tiny{$D_1$}}
\move (13 48)  \htext{\tiny{$C_2$}}  \move (43 48)  \htext{\tiny{$D_2$}}
\move (54 48)  \htext{\tiny{$C_3$}}  \move (83 48)  \htext{\tiny{$D_3$}}
\end{texdraw}
\end{center}
where we have indicated the branches of $C$ and $D$ corresponding to the terminal trunks. Thus,  the bijection among the bars in $T(f,g)$ and the bifurcation divisors in $G(Z)$ is given by
$$
B_0 \longleftrightarrow E_1; \qquad  B_1 \longleftrightarrow E_{d+1}; \qquad  B_2 \longleftrightarrow E'_{f+1}; \qquad  B_3 \longleftrightarrow E_{f+1}
$$
with $h(B_0)=v(E_1)=1$, $h(B_1)=v(E_{d+1})=d+1$, $h(B_2)=h(B_3)=v(E_{f+1})=v(E_{f+1}')=f+1$.
\end{Example}

Let us show that the notion of collinear point and collinear divisor given in section~\ref{sec:collinear} correspond to the ones given in \cite{Kuo-P-2002,Kuo-P-2004} thanks to the bijections explained above. Let $B$ be a bar of $T(f,g)$  and consider a Puiseux series $z_k(x)$ of $Z$ which goes through $B$, this means, that
$$
z_k(x)= z_B(x) + c x^{h(B)} + \cdots
$$
where $z_B(x)$ depends only on the bar $B$ and $c$ is uniquely determined by $z_k(x)$. If $T$ is a trunk which contains $z_k(x)$, then it is said that the trunk $T$ {\em grows on} $B$ {\em at} $c$. Let $\widetilde{E}$ be the bifurcation divisor of $G(\widetilde{Z})$ corresponding to $B$ and consider
$$
v_k(u)=z_k(u^n)= z_B(u^n) + c u^{nh(B)} + \cdots
$$
Note that the curve given by $v-v_k(u)=0$ is a branch of $\widetilde{Z}$ such that $\widetilde{E}$ belongs to its geodesic. This curve determines a unique point $R$ in $\widetilde{E}_{red}$; in this way we can establish the  bijection among the trunks on $B$ and the infinitely near points of $\widetilde{Z}$ on $\widetilde{E}_{red}$.

Let now $B$ be a bar of $T(f,g)$ that corresponds to a bifurcation divisor $\widetilde{E}$ of $G(\widetilde{Z})$ and $T_k$, $1\leq k \leq b_{\widetilde{E}}$, be the set of trunks on $B$ with bimultiplicity $[p_k,q_k]$ where the trunk $T_i$ grows on $B$ at $c_i$. Let us denote $\{R_1^{\widetilde{E}}, \ldots, R_{b_{\widetilde{E}}}^{\widetilde{E}}\}$ the set of infinitely near points of $\widetilde{Z}$ in $\widetilde{E}_{red}$ with $R_i^{\widetilde{E}}$ corresponding to the trunk $T_i$.

In \cite{Kuo-P-2004}, the authors define %it is introduced the following notation
$$\Delta_{B}(c_k)=\left| \begin{array}{cc}
\nu_f(B) & p_k \\
\nu_g(B) & q_k
\end{array} \right|, \qquad 1 \leq k \leq b_{\widetilde{E}}$$
where $\nu_f(B)=\text{ord}_x(f(x,z_B(x)+cx^{h(B)}))$ for $c \in {\mathbb C}$ generic (resp. $\nu_g(B)$), and   the {\em rational function associated to} $B$ as
$$
{\mathcal M}_B(z)= \sum_{k=1}^{b_{\widetilde{E}}} \frac{\Delta_B(c_k)}{z-c_k}.
$$
Note that, if $E$ is the bifurcation divisor of $G(Z)$ such that $\widetilde{E}$ is associated to $E$, then the curve given by $y=z_B(x)+cx^{h(B)}$ is an $E$-curvette. Thus, taking into account Proposition 2.5.3 of \cite{Cas-00} for instance, we get that
$$
\nu_f(B)= \frac{(C,\gamma_E)_0}{m(E)}
$$
with $\gamma_E$ any $E$-curvette. Moreover, it is easy to verify that
$$
\nu_f(B)= \frac{1}{n}  \sum_{i=1}^{\nu_0(C)} (\gamma_i^C, \gamma_{\widetilde{E}})_0
$$
where $\gamma_i^C$ is the curve given by $v-y_i^C(u^n)=0$ and $\gamma_{\widetilde{E}}$ is an $\widetilde{E}$-curvette. Note that we can compute the intersection multiplicity $(\gamma_i^C,\gamma_{\widetilde{E}})_0=\sum_{\widetilde{E}' \leq \widetilde{E}} \varepsilon_{\widetilde{E}'}^{\gamma_i^C}$ where the sum runs over all the divisors $\widetilde{E}'$ in $G(\widetilde{C})$ in the geodesic of $\widetilde{E}$ and $\varepsilon_{\widetilde{E}'}^{\gamma_i^C}=1$ if the geodesic of $\gamma_i^C$ contains the divisor $\widetilde{E}'$ and $\varepsilon_{\widetilde{E}'}^{\gamma_i^C}=0$ otherwise. Thus, with the notations given in section~\ref{sec:logarithmic}, we have that
$$\nu_f(B)=\frac{1}{n} \kappa_{\widetilde{E}}({\mathcal L}^{\widetilde{C}})
$$
where ${\mathcal L}^{\widetilde{C}}={\mathcal G}_{\tilde{f}}$ is the logarithmic foliation in ${\mathbb G}_{\widetilde{C}}$ with $\lambda=(1,1,\ldots,1)$, that is, the hamiltonian foliation defined by $d \tilde{f}=0$ with $\tilde{f}(u,v)=f(u^n,v)$. Moreover,
$$
p_k=\nu_{R_k^{\widetilde{E}}}(\pi_{\widetilde{E}}^* \widetilde{C}), \qquad k=1,\ldots, b_{\widetilde{E}}
$$
and thus
$$
\mathcal{I}_{R_k^{\widetilde{E}}}(\pi_{\widetilde{E}}^*{\mathcal L}^{\widetilde{C}},\widetilde{E}_{red})= - \frac{p_k}{n \nu_f(B)}, \qquad k=1,\ldots, b_{\widetilde{E}}.
$$
Consequently, with the notations introduced in section~\ref{sec:collinear}, we have that
$$\Delta_B(c_k) = -\frac{1}{n} \Delta_{\widetilde{E}}(R_k^{\widetilde{E}}) \quad
\text{ and } \quad
{\mathcal M}_{\widetilde{E}}(z) = - n {\mathcal M}_B(z).
$$
Thus the  notions of collinear divisor and collinear point  given in section~\ref{sec:collinear} correspond to the ones given in \cite{Kuo-P-2004} for bars and points on them, and  the results given in section~\ref{sec:jacobian} imply  some of the Theorems proved in \cite{Kuo-P-2004}.
%Theorem~\ref{th-multiplicidad-jac} implies Theorem T and Corollary 2.4 in \cite{Kuo-P-2004}, Corollary~\ref{cor:consecutive-divisors} implies Corollary 2.6 in \cite{Kuo-P-2004} and Theorem~\ref{th:colinear-point} implies Theorem C in \cite{Kuo-P-2004}. %correspond to the results....

%%%%%%%%%%%%%%%%%%%%%%%%%%%%%%%%%%%%%%%%%%%%%%%%%%%%%%%%
\subsection{Semiroots and Approximate roots}\label{sec:aprox-roots}
The notion of approximated root was introduced by Abhyankar and Moh in \cite{Abh-M} where they proved the following result:
\begin{Proposition}
  Let $A$ be an integral domain and $P(y) \in A[y]$ be a monic polynomial of degree $d$. If $p$ is invertible in $A$ and $p$ divides $d$, then there exists a unique monic polynomial $Q(y) \in A[y]$ such that the degree of $P-Q^p$ is less than $d-d/p$.
\end{Proposition}
The unique polynomial $Q$ given by the previous proposition is called the {\em $p$-th approximate root\/} of $P$. Let us consider $f(x,y) \in {\mathbb C}\{x\}[y]$ an irreducible Weierstrass polynomial with characteristic exponents $\{\beta_0,\beta_1,\ldots,\beta_g\}$ and denote $e_k=\gcd(\beta_0,\beta_1,\ldots, \beta_k)$ for $k=1,\ldots, g$. Thus $e_k$ divides $\beta_0=deg_y f$. We will denote $f^{(k)}$ the $e_k$-approximate root of $f$ and we call them the {\em characteristic approximate roots\/} of $f$.
Next result (Theorem 7.1 in \cite{Abh-M}) gives the  main properties of the characteristic approximate roots of $f$ (see also \cite{Gwo-P-1995,Pop-2003}):
\begin{Proposition}\label{prop:char-root}
Let $f(x,y) \in {\mathbb C}\{x\}[y]$ be  an irreducible Weierstrass polynomial with characteristic exponents $\{\beta_0,\beta_1,\ldots,\beta_g\}$. Then the characteristic approximate roots $f^{(k)}$ for $k=0,1,\ldots, g-1$ verify:
\begin{itemize}
  \item[(1)] The degree in $y$ of $f^{(k)}$ is equal to $\beta_0/e_k$ and ${\mathcal C}(f,f^{(k)})=\beta_{k+1}/\beta_0$.
  \item[(2)] The polynomial $f^{(k)}$ is irreducible with characteristic exponents \linebreak $\{\beta_0/e_k,\beta_1/e_k,\ldots, \beta_k/e_k\}$.
\end{itemize}
\end{Proposition}
In \cite{Gar-G}, E. García Barroso and J. Gwo\'{z}dziewicz   studied the jacobian curve of $f$ and $f^{(k)}$ and they give a result concerning its factorization (see Theorem 1 of \cite{Gar-G}). In this section, we will prove that this result of factorization can be obtained as a consequence of Theorem~\ref{th:desc-general}.

\begin{Remark}
In \cite{Pop-2003}, P. Popescu-Pampu proved that all polynomials in ${\mathbb C}\{x\}[y]$ satisfying condition (1) in the proposition above also verify condition (2). Hence, given  an irreducible Weierstrass polynomial $f(x,y) \in {\mathbb C}\{x\}[y]$ with characteristic exponents $\{\beta_0,\beta_1,\ldots, \beta_g\}$, we can consider the monic polynomials in ${\mathbb C}\{x\}[y]$  satisfying condition (1) above which are called {\em $k$-semiroots\/} of $f$ (see Definition 6.4 in \cite{Pop-2003}). Since we only need the properties of characteristic approximate roots given in Proposition~\ref{prop:char-root}, in the rest of the section, we will denote by $f^{(k)}$ a $k$-semiroot of $f$, $0 \leq k \leq g-1$.
\end{Remark}

Let $C$ be the curve defined by $f=0$ and denote $C^{(k)}$ the curve given by $f^{(k)}=0$ with $0 \leq k \leq g-1$. Consider ${\mathcal F} \in {\mathbb G}_{C}$ and ${\mathcal F}^{(k)} \in {\mathbb G}_{C^{(k)}}$. Note that the minimal reduction of singularities $\pi_C: X_C \to ({\mathbb C}^2,0)$ of the curve $C$ gives also a reduction of singularities of $C \cup C^{(k)}$. There are $g$ bifurcation divisors in $G(C)$. The set of bifurcation divisors of $G(C)$ will be denote by $\{E_1,\ldots, E_g\}$ with $v(E_i)=\frac{\beta_i}{\beta_0}$.
Remark that the dual graph  $G(C \cup C^{(k)})$ is given by (see \cite{Pop-2003} for instance):

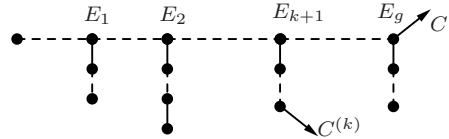
\vspace*{0.5\baselineskip}
\begin{center}
\begin{texdraw}
\arrowheadtype t:F \arrowheadsize l:0.08 w:0.04 \drawdim mm \setgray
0
%\move(10 20) \fcir f:0 r:0.8 \rlvec (5 0) \fcir f:0 r:0.8 \lpatt(1 1)  \rlvec (5 0)  \fcir f:0 r:0.8 \lpatt( ) \rlvec (5 0) \fcir f:0 r:0.8  \rlvec (0 -4)
\move(15 20) \fcir f:0 r:0.8  \lpatt(1 1)  \rlvec (10 0)  \fcir f:0 r:0.8 \lpatt( )   \rlvec (0 -4)
\fcir f:0 r:0.8 \lpatt(1 1)  \rlvec (0 -4)  \fcir f:0 r:0.8
\move(25 20) \fcir f:0 r:0.8 \rlvec (10 0) \fcir f:0 r:0.8  \lpatt( ) \rlvec (0 -4) %\fcir f:0 r:0.8
\fcir f:0 r:0.8 \lpatt(1 1)  \rlvec (0 -4)  \fcir f:0 r:0.8 \lpatt( ) \rlvec (0 -4)  \fcir f:0 r:0.8
\move(35 20) \fcir f:0 r:0.8 \lpatt(1 1) \rlvec (15 0) \fcir f:0 r:0.8  \lpatt() \rlvec (0 -4)
\fcir f:0 r:0.8 \lpatt(1 1)  \rlvec (0 -5)  \fcir f:0 r:0.8
\move(50 20) \fcir f:0 r:0.8 \rlvec (15 0) \fcir f:0 r:0.8 \lpatt() \rlvec (0 -4)
\fcir f:0 r:0.8 \lpatt(1 1)  \rlvec (0 -5)  \fcir f:0 r:0.8
\lpatt( ) \move(50 11) \avec (55 7) \htext{\tiny$C^{(k)}$}
\move(65 20) \avec(70 24)
\move (70 21) \htext{\tiny$C$}
% \rlvec
%(0 -7) \fcir f:0 r:0.8 \move (20 20) \avec (25 23) \move (26 21)

%\htext{\tiny{$C_1$}}

%\move (20 13)  \avec (25 16)

\move (24 22) \htext{\tiny$E_1$}
\move (34 22) \htext{\tiny{$E_2$}}
\move (49 22) \htext{\tiny{$E_{k+1}$}}
\move (63 22) \htext{\tiny{$E_{g}$}}
%\move (16 10) \htext{\tiny{$E_2$}} \move (26 14) \htext{\tiny{$C_2$}}
\end{texdraw}
\end{center}
Thus the sets of bifurcation divisors of $G(C)$ and $G(C \cup C^{(k)})$ coincide. All bifurcation divisors of $G(C \cup C^{(k)})$ are Puiseux divisors for $C$ while only $E_1, \ldots, E_k$ are Puiseux divisors for $C^{(k)}$.
Then we have
\begin{Lemma}\label{lema:colinear-irreducible}
The set of non-collinear bifurcation divisors of $G(C \cup C^{(k)})$ for the foliations $\mathcal F$ and ${\mathcal F}^{(k)}$ is
$\{E_{k+1}, \ldots, E_g\}.$
\end{Lemma}
\begin{proof}
Let $\{(m_1,n_1), \ldots, (m_g,n_g)\}$ be the Puiseux pairs of $C$, then we remind that
$\beta_0=\nu_0(C)= n_1 \cdots n_g$, $e_k=n_{k+1} \cdots n_g$ and $\beta_k/\beta_0=m_k/n_1 \cdots n_k$ for $k=1,\ldots, g$.
Given a bifurcation divisor $E_l$ of $G(C \cup C^{(k)})$, we have that $n_{E_l}=n_l$,
$\underline{n}_{E_l}=n_1 \cdots n_{l-1}=\beta_0/e_{l-1}$ and $m(E_l)=\underline{n}_{E_l}n_E = n_1 \cdots n_l$.

Consider now the ramification  $\rho:({\mathbb C}^2,0) \to ({\mathbb C}^2,0) $   given by $\rho(u,v)=(u^n,v)$ with $n=\beta_0$ and denote $\widetilde{C}=\rho^{-1}C$, $\widetilde{C}^{(k)}=\rho^{-1} C^{(k)}$.
Take a bifurcation divisor $E_l$  and let $\{\widetilde{E}_l^t\}_{t=1}^{\underline{n}_{E_l}}$ be the set of bifurcation divisors of $G(\widetilde{C}\cup \widetilde{C}^{(k)})$ associated to $E_l$.

In the case $l < k+1$, we have that $\pi_{\widetilde{E}_l^t}^* \widetilde{C} \cap \widetilde{E}_{l,red}^t=\pi_{\widetilde{E}_l^t}^* \widetilde{C}^{(k)} \cap \widetilde{E}_{l,red}^t$
with $b_{\widetilde{E}_l^t}=n_l$ in $G(\widetilde{C} \cup \widetilde{C}^{(k)})$. Let us denote $\pi_{\widetilde{E}_l^t}^* \widetilde{C} \cap \widetilde{E}_{l,red}^t=\{R_1^{\widetilde{E}_l^t}, \ldots, R_{b_{\widetilde{E}_l^t}}^{\widetilde{E}_l^t}\}$.
Using the equations given in section~\ref{sec:logarithmic} and appendix~\ref{ap:ramificacion}, the computation of the Camacho-Sad indices for the foliations $\rho^* {\mathcal F}$ and $\rho^* {\mathcal F}^{(k)}$ gives
\begin{align}
\mathcal{I}_{R_s^{\widetilde{E}_l^t}}(\pi_{\widetilde{E}_l^t}^* \rho^* {\mathcal F}, \widetilde{E}_{l,red}^t) & =- \frac{n_{l+1} \cdots n_g}{\sum_{s=1}^{n} \sum_{\widetilde{E} \leq \widetilde{E}_{l}^t } \varepsilon_{\widetilde{E}}^{\sigma_s}} \label{eq:indice-C} \\
\mathcal{I}_{R_s^{\widetilde{E}_l^t}}(\pi_{\widetilde{E}_l^t}^* \rho^* {\mathcal F^{(k)}}, \widetilde{E}_{l,red}^t) & =- \frac{n_{l+1} \cdots n_k}{\sum_{s=1}^{n_1 \cdots n_k} \sum_{\widetilde{E} \leq \widetilde{E}_{l}^t } \varepsilon_{\widetilde{E}}^{\sigma_s^{(k)}}} \notag
\end{align}
where $\widetilde{C}=\cup_{s=1}^n \sigma_s$ and $\widetilde{C}^{(k)}= \cup_{s=1}^{n_1 \cdots n_k} \sigma_s^{(k)}$. Hence, taking into account the results of appendix~\ref{ap:ramificacion}, we have
$$\mathcal{I}_{R_s^{\widetilde{E}_l^t}}(\pi_{\widetilde{E}_l^t}^* \rho^* {\mathcal F}, \widetilde{E}_{l,red}^t) =\mathcal{I}_{R_s^{\widetilde{E}_l^t}}(\pi_{\widetilde{E}_l^t}^* \rho^* {\mathcal F^{(k)}}, \widetilde{E}_{l,red}^t), \qquad s=1, \ldots, b_{\widetilde{E}_l^t}$$
and consequently, $\Delta_{\widetilde{E}_l^t}(R_s^{\widetilde{E}_l^t}) =0$, $s=1, \ldots, b_{\widetilde{E}_l^t}$, for the foliations $\rho^* {\mathcal F}$ and $\rho^* {\mathcal F}^{(k)}$. This proves that the bifurcation divisors $E_l$ of $G(C \cup C^{(k)})$ with $l < k+1$ are collinear for ${\mathcal F}$ and ${\mathcal F}^{(k)}$.

Consider now the bifurcation divisor $E_{k+1}$ of $G(C \cup C^{(k)})$ and let
$\{\widetilde{E}_{k+1}^t\}_{t=1}^{\underline{n}_{E_{k+1}}}$ be the set of bifurcation divisors of $G(\widetilde{C}\cup \widetilde{C}^{(k)})$ associated to $E_{k+1}$. Although the curve $\pi_{E_{k+1}}^* C^{(k)}$ does not intersect $E_{k+1,red}$, the curve $\pi_{\widetilde{E}_{k+1}^t}^* \widetilde{C}^{(k)}$ intersects $\widetilde{E}_{k+1,red}^t$ in one point for each $t=1,\ldots, \underline{n}_{E_{k+1}}$ which is different from the $n_{k+1}$ points where $\pi_{\widetilde{E}_{k+1}^t}^* \widetilde{C}$ intersects $\widetilde{E}_{k+1,red}^t$. Note that $b_{E_{k+1}}=2$ in $G(C \cup C^{(k)})$ and, by equation~\eqref{eq:b-tilde}, $b_{\widetilde{E}_{k+1}^t}=n_{k+1}+1$ in $G(\widetilde{C} \cup \widetilde{C}^{(k)})$. Let $\{R_1^{\widetilde{E}_{k+1}^t}, \ldots, R_{b_{\widetilde{E}_{k+1}^t}}^{\widetilde{E}_{k+1}^t}\}$ be the set of points $(\pi_{\widetilde{E}_{k+1}^t}^* \widetilde{C} \cap \widetilde{E}_{k+1,red}^t) \cup (\pi_{\widetilde{E}_{k+1}^t}^* \widetilde{C}^{(k)} \cap \widetilde{E}_{k+1,red}^t)$ with $R_{b_{\widetilde{E}_{k+1}^t}}^{\widetilde{E}_{k+1}^t}=\pi_{\widetilde{E}_{k+1}^t}^* \widetilde{C}^{(k)} \cap \widetilde{E}_{k+1,red}^t$. Thus, we can compute the Camacho-Sad of $\rho^*\mathcal{F}$ and $\rho^*\mathcal{F}^{(k)}$ at these points as in the previous case, and prove that
$\Delta_{\widetilde{E}_{k+1}^t}(R_s^{\widetilde{E}_{k+1}^t}) \neq 0$, $s=1,\ldots, b_{\widetilde{E}_{k+1}^t}$, for the foliations $\rho^* {\mathcal F}$ and $\rho^* {\mathcal F}^{(k)}$. Consequently $E_{k+1}$ is a non-collinear divisor for $\mathcal F$ and ${\mathcal F}^{(k)}$. However, we have that
\begin{equation}\label{eq:Delta-E-k+1}
\sum_{s=1}^{b_{\widetilde{E}_{k+1}^t}} \Delta_{\widetilde{E}_{k+1}^t}(R_s^{\widetilde{E}_{k+1}^t}) = 0.
\end{equation}
In fact, the divisor $\widetilde{E}_{k+1}^t$ arises from one of the divisors $\widetilde{E}_{k}^r$ at one of the points $R_t^{\widetilde{E}_{k}^r}$ of the set $\pi_{\widetilde{E}_k^r}^* \widetilde{C} \cap \widetilde{E}_{k,red}^r$. Since $\widetilde{E}_{k}^r$ is a collinear divisor, then $R_t^{\widetilde{E}_{k}^r}$ is a collinear point and equation~\eqref{eq:Delta-E-k+1}
follows from Corollary~\ref{cor:colinear}.

Consider now a bifurcation divisor $E_l$ of $G(C \cup C^{(k)})$ with $l > k+1$. In this case, we have that the curve $\pi_{E_l}^* C^{(k)}$ does not intersect $E_{l,red}$ , the curve $\pi_{\widetilde{E}_{l}^t}^* \widetilde{C}^{(k)}$ does not intersect $\widetilde{E}_{l,red}^t$ and $b_{\widetilde{E}_{l}^t}=n_{l}$ in $G(\widetilde{C} \cup \widetilde{C}^{(k)})$ (see equation~\eqref{eq:b-tilde}).
Let us denote  $\{R_1^{\widetilde{E}_l^t}, \ldots, R_{b_{\widetilde{E}_l^t}}^{\widetilde{E}_l^t}\}$ the set of points $\pi_{\widetilde{E}_{l}^t}^* \widetilde{C} \cap \widetilde{E}_{l,red}^t =  \pi_{\widetilde{E}_{l}^t}^* (\widetilde{C} \cup \widetilde{C}^{(k)}) \cap \widetilde{E}_{l,red}^t$. With the notations above, for $s \in \{1,\ldots, b_{\widetilde{E}_l^t}\}$, we have that $\mathcal{I}_{R_s^{\widetilde{E}_l^t}}(\pi_{\widetilde{E}_l^t}^* \rho^* {\mathcal F}, \widetilde{E}_{l,red}^t)$ is given by equation \eqref{eq:indice-C} while $\mathcal{I}_{R_s^{\widetilde{E}_l^t}}(\pi_{\widetilde{E}_l^t}^* \rho^* {\mathcal F^{(k)}}, \widetilde{E}_{l,red}^t)=0$ since the points $R_s^{\widetilde{E}_l^t}$ are non-singular points for $\rho^*{\mathcal F}^{(k)}$. This implies $\Delta_{\widetilde{E}_l^t}(R_s^{\widetilde{E}_l^t}) \neq 0$, $1 \leq s \leq b_{\widetilde{E}_l^t}$, and hence $E_l$ is a non-collinear divisor for the foliations $\rho^* {\mathcal F}$ and $\rho^* {\mathcal F}^{(k)}$. Moreover, we have that
$$
\sum_{s=1}^{b_{\widetilde{E}_{l}^t}} \Delta_{\widetilde{E}_{l}^t}(R_s^{\widetilde{E}_{l}^t}) = - \sum_{s=1}^{b_{\widetilde{E}_{l}^t}} \mathcal{I}_{R_s^{\widetilde{E}_l^t}}(\pi_{\widetilde{E}_l^t}^* \rho^* {\mathcal F}, \widetilde{E}_{l,red}^t) = 1 + \mathcal{I}_Q(\pi_{\widetilde{E}_l^t}^* \rho^* {\mathcal F}, \widetilde{E}_{l,red}^t)
$$
where $Q$ is the only singular point of $\pi_{\widetilde{E}_l^t}^* \rho^* {\mathcal F}$ in $\widetilde{E}_l^t$ different from the points $R_s^{\widetilde{E}_l^t}$. By Proposition 4.4 in \cite{Cor-2003}, we know that $\mathcal{I}_Q(\pi_{\widetilde{E}_l^t}^* \rho^* {\mathcal F}, \widetilde{E}_{l,red}^t) \neq -1$ and hence
$$
\sum_{s=1}^{b_{\widetilde{E}_{l}^t}} \Delta_{\widetilde{E}_{l}^t}(R_s^{\widetilde{E}_{l}^t}) \neq 0.
$$
\end{proof}

Thus, by Theorem~\ref{th:desc-general} there is a decomposition
$${\mathcal J}_{\mathcal{F},{\mathcal F}^{(k)}}= J^* \cup \left( \bigcup_{i=k+1}^g J^{i} \right)$$
where $J^i=J^{E_i}_{nc}$, such that
\begin{itemize}
  \item[(i)] $\nu_0(J^{k+1}) < n_1 \cdots n_{k+1}$.
  \item[(ii)] $\nu_0(J^i)= n_1 \cdots n_{i-1} (n_{i} -1)$ for $k+2 \leq i \leq g$.
  \item[(iii)] if $\gamma$ is a branch of $J^i$, $k+1 \leq i \leq g$, we have that ${\mathcal C}(\gamma,C)=\frac{\beta_i}{\beta_0}$.
  \item[(iv)] if $\gamma$ is a branch of $J^*$, then  ${\mathcal C}(\gamma,C)< \frac{\beta_{k+1}}{\beta_0}$.
\end{itemize}
Note that $\delta$ is a branch of $J^*$ if it is not a branch of any of the curves $J^{i}$, $k+1 \leq i \leq g$, and hence, $\pi_C^*\gamma$ intersects a component $E$ of the exceptional divisor $\pi^{-1}_C(0)$ which appears in the reduction of singularities of $C$ before than $E_{k+1}$. Consequently we have that ${\mathcal C}(\gamma,C)< v(E_{k+1})=\frac{\beta_{k+1}}{\beta_0}$ which gives property (iv).

 Let us prove that $J^{k+1}= \emptyset$. From section~\ref{sec:caso-general} we have that $\rho^{-1} J^{k+1}= \bigcup_{t=1}^{n_1 \cdots n_{k}} J_{nc}^{\widetilde{E}^t_{k+1}}$ with the notations of the proof of Lemma~\ref{lema:colinear-irreducible}. Let us compute ${\mathcal M}_{\widetilde{E}^t_{k+1}}(z)$ for any $t \in \{1,\ldots, n_1 \cdots n_k\}$. To simplify notations, let us denote $\tilde{b}=b_{\widetilde{E}^t_{k+1}}=n_{k+1}+1$,  $\widetilde{R}_s= R_s^{\widetilde{E}^t_{k+1}}$ and $\Delta (\widetilde{R}_l) =\Delta_{\widetilde{E}^t_{k+1}}(\widetilde{R}_l)$. Thus, we have that
$$
{\mathcal M}_{\widetilde{E}^t_{k+1}}(z) = \sum_{s=1}^{\tilde{b}-1}\frac{\Delta (\widetilde{R}_s)}{z-\xi^s} + \frac{\Delta (\widetilde{R}_{\tilde{b}})}{z}
$$
where $\xi$ is a primitive $n_{k+1}$-root of a value $a=a^{\widetilde{E}^t_{k+1}}$ determined by the Puiseux parametrizations of $C$. From the proof of Lemma~\ref{lema:colinear-irreducible}, we obtain that $\Delta(\widetilde{R}_s)=\Delta(\widetilde{R}_t)$ for any $s,t \in \{1,\ldots, \tilde{b}-1\}$. Thus, taking into account equation~\eqref{eq:Delta-E-k+1}, we get that
\begin{align*}
{\mathcal M}_{\widetilde{E}^t_{k+1}}(z)  & = \Delta (\widetilde{R}_1)\frac{\sum_{s=1}^{\tilde{b}-1} \prod_{t=1 \atop t \neq s}^{n_{k+1}}(z-\xi^t)}{z^{n_{k+1}}-a} + \frac{\Delta (\widetilde{R}_{\tilde{b}})}{z} \\
 & = \Delta(\widetilde{R}_1) \frac{n_{k+1} z^{n_{k+1} -1}}{z^{n_{k+1}}-a} + \frac{\Delta (\widetilde{R}_{\tilde{b}})}{z} =\frac{ - a \Delta(\widetilde{R}_{\tilde{b}})}{z(z^{n_{k+1}}-a)}
\end{align*}
where $a \Delta(\widetilde{R}_{\tilde{b}}) \neq 0$. Hence, $J_{nc}^{\widetilde{E}^t_{k+1}}=\emptyset$ for all $t=1, \ldots,n_1 \cdots n_{k}$ and consequently $J^{k+1}=\emptyset$.

\begin{Corollary}\label{Cor-raiz-aprox}
Let $C$ be an irreducible curve and $C^{(k)}$ the curve given by the $k$-characteristic approximate root (or by a $k$-semiroot) with $0 \leq k \leq g-1$. Consider ${\mathcal F} \in {\mathbb G}_{C}$ and ${\mathcal F}^{(k)} \in {\mathbb G}_{C^{(k)}}$. Thus, the jacobian curve ${\mathcal J}_{\mathcal{F},{\mathcal F}^{(k)}}$ has a decomposition

$${\mathcal J}_{\mathcal{F},{\mathcal F}^{(k)}}= J^* \cup \left( \bigcup_{i=k+2}^g J^{i} \right)$$
 such that
\begin{itemize}
  \item[(1)] $\nu_0(J^i)= n_1 \cdots n_{i-1} (n_{i} -1)$ for $k+2 \leq i \leq g$.
  \item[(2)] if $\gamma$ is a branch of $J^i$, $k+2 \leq i \leq g$, we have that ${\mathcal C}(\gamma,C)=\frac{\beta_i}{\beta_0}$.
  \item[(3)] if $\gamma$ is a branch of $J^*$, then  ${\mathcal C}(\gamma,C)< \frac{\beta_{k+1}}{\beta_0}$.
\end{itemize}
\end{Corollary}
In particular, the result above implies the result of E. García Barroso and J. Gwo\'{z}dziewicz  (Theorem 1 in \cite{Gar-G}) concerning the jacobian curve of a plane curve and its characteristic approximate roots.

Moreover, in \cite{Sar}, it is considered the jacobian curve ${\mathcal J}_{{\mathcal F},{\mathcal G}^{(k)}}$ of a foliation $\mathcal F$ with an irreducible separatrix $f=0$ and the hamiltonian foliation ${\mathcal G}^{(k)}$ defined by $df^{(k)}=0$ with $f^{(k)}$ a characteristic approximate root of $f$. Corollary \ref{Cor-raiz-aprox} also implies the main result of N. E. Saravia in \cite{Sar} (Theorem 3.4.1) concerning factorization  of ${\mathcal J}_{{\mathcal F},{\mathcal G}^{(k)}}$ .

\begin{Remark}
Note that $E^1$ is always a collinear divisor for the foliations $\mathcal F$ and $\mathcal{F}^{(k)}$,  hence the hypothesis of Lemma~\ref{lema:multiplicidad-E1} are not satisfied and  the multiplicity of ${\mathcal J}_{\mathcal{F},\mathcal{F}^{(k)}}$ can be greater than $\nu_0(\mathcal{F}) + \nu_0(\mathcal{F}^{(k)})=\nu_0(C)+\nu_0(C^{(k)})-2$ as showed in the examples given in  \cite{Gar-G} or \cite{Sar}.
\end{Remark}
\subsection{Polar curves of foliations}\label{sec:polars}
Given a germ of foliation $\mathcal F$ in $({\mathbb C}^2,0)$, a polar curve of $\mathcal F$ corresponds to the jacobian curve of $\mathcal F$ and a non-singular foliation $\mathcal G$. If we are interested in the topological properties of a generic polar curve of $\mathcal F$, it is enough to consider a generic curve ${\mathcal P}_{[a:b]}^\mathcal{F}$ in the family of curves given by
$$\omega \wedge (a dy - b dx) =0$$
where $\omega=0$ is a 1-form defining $\mathcal F$ and $[a:b] \in {\mathbb P}^1_{\mathbb C}$ (see \cite{Cor-2003} section 2). When $\mathcal F$ is a hamiltonian foliation given by $df=0$ we recover the notion of polar curve of a plane curve. As we mention in the introduction, polar curves play an important role in the study of singularities of plane curves and also of foliations. There is a result, known as ``decomposition theorem'', which describes the minimal topological properties of the generic polar curve of a plane curve $C$ in terms of the topological type of the curve $C$ (see \cite{Mer} for the case of $C$ irreducible; \cite{Gar} for $C$ with several branches). In the case of foliations, the decomposition theorem also holds for the generic polar curve of a generalized curve foliation $\mathcal F$ with an irreducible separatrix (see \cite{Rou}). In the general case of a generalized curve foliation $\mathcal F$ whose curve of separatrices is not irreducible, the decomposition theorem for its generic polar curve only holds under some conditions on the foliation $\mathcal F$ (see \cite{Cor-2003}). Let us see that all these results  can be recovered from the results in this paper. In particular, we show that we can prove Theorems 5.1 and 6.1 in \cite{Cor-2003} which give the decomposition theorem for the polar curve of a generalized curve foliation $\mathcal F$ and hence we get all the other results concerning decompositions theorems.

Let $\mathcal F$ be a generalized curve foliation in $({\mathbb C}^2,0)$ with $C$ as curve of separatrices and denote by $\mathcal{P}^\mathcal{F}$ a generic polar curve of $\mathcal F$. We can assume that $\mathcal{P}^\mathcal{F}={\mathcal J}_{\mathcal{F},\mathcal{G}}$ where $\mathcal{G}$ is a non-singular foliation. Note that the curve  of separatrices $D$ of $\mathcal G$ is a non-singular irreducible  plane curve. Let us assume  first that  $C$ has  only non-singular irreducible components, all of them different from $D$, and take the notations of section \ref{sec:collinear}. Thus the minimal reduction of singularities $\pi_C: X_C \to ({\mathbb C}^2,0)$ is also the minimal reduction of singularities of $Z=C \cup D$. Note that the dual graph $G(Z)$ is obtained from $G(C)$ adding an arrow to the first divisor $E^1$ which represents the curve $D$. Hence, if we denote $b_E^Z$, $b_E^C$ the number associated to a divisor $E$ in $G(Z)$ or $G(C)$ respectively, as defined in subsection \ref{Apendice:grafo-dual}, then $b_{E^1}^Z=b_{E^1}^C+1$ and $b_{E}^Z=b_{E}^C$ otherwise.

Consider $E$ an irreducible component  of the exceptional divisor $\pi_C^{-1}(0)$. If $E=E^1$ is the divisor which appears after the blow-up of the origin, then $\pi_{E^1}^* D \cap E_{red}^1=\{Q\}$ and $Q \not \in \pi_{E^1}^* C \cap E_{red}^1$. Thus, for $R \in  E_{red}^1$ we have that

$$\Delta_{E^{1}}^{\mathcal{F},\mathcal{G}} (R) = \left\{
                                             \begin{array}{ll}
                                               \mathcal{I}_{Q} (\pi_{E^1}^*\mathcal{G},E_{red}^1), & \hbox{ if } R=Q \\
                                               -\mathcal{I}_{R} (\pi_{E^1}^*\mathcal{F},E_{red}^1), & \hbox{ otherwise.}
                                             \end{array}
                                           \right.
$$
with  $\mathcal{I}_{Q} (\pi_{E^1}^*\mathcal{G},E_{red}^1)=-1$.
If $E \neq E_1$, then $\pi^{*}_E Z \cap E_{red}= \pi^*_E C \cap E_{red}$ and then $\Delta_E^{\mathcal{F},\mathcal{G}} (R)=-\mathcal{I}_{R} (\pi_E^*\mathcal{F},E_{red})$ for any $R \in E_{red}$.

With the hypothesis above and the notations of section \ref{sec:collinear}, we have that
\begin{Lemma} The following conditions are equivalent:
\begin{itemize}
  \item[(i)] There is no corner in $\pi^{-1}_C(0)$ such that $\pi_C^* \mathcal{F}$ has  Camacho-Sad index equal to $-1$.
  \item[(ii)] All the components of the exceptional divisor $\pi^{-1}_C(0)$ are purely non-collinear.
\end{itemize}
\end{Lemma}
\begin{proof}
Assume that (i) holds and that there is a component $E$ of the exceptional divisor which is not purely non-collinear, that is, there is a singular point $R \in E_{red}$ of $\pi_E^* \mathcal{F}$ with $\mathcal{I}_{R} (\pi_E^*\mathcal{F},E_{red})=0$. Then $R$ is not a simple singularity for $\pi_E^*\mathcal{F}$ and hence, if $\sigma: X_{E'} \to  X_E$ is the blow-up with center in $R$, and we denote by $\tilde{E}_{red}$ the strict transform of $E_{red}$ by $\sigma$, then
we have that  $\mathcal{I}_{\tilde{R}} (\pi_{\tilde{E}}^*\mathcal{F},\tilde{E}_{red})=-1$ where $\tilde{R}=E'_{red} \cap  \tilde{E}_{red}$. Thus we get a corner in $\pi_C^{-1}(0)$ with Camacho-Sad index equal to $-1$.

Conversely, assume now that all the components of the exceptional divisor $\pi^{-1}_C(0)$ are purely non-collinear and there is a corner $\tilde{R} = E_{k-1} \cap E_{k}$ with $I_{\tilde{R}} (\pi_{C}^*\mathcal{F},E_{k-1})=-1$. Consider the morphism $\pi_{E_{k-1}}: X_{E_{k-1}} \to (\mathbb{C}^2,0)$ and take the point $R \in E_{k-1,red}$ that we have to blow-up to obtain the divisor $E_{k}$. Thus, by the properties of the Camacho-Sad, we have that $\mathcal{I}_{R} (\pi_{E_{k-1,red}}^*\mathcal{F},E_{k-1,red})=0$ but this contradicts that $E_{k-1}$ is purely non-collinear.
\end{proof}
In particular, let us see that Theorem 6.1 in \cite{Cor-2003} is consequence of Theorem \ref{th-multiplicidad-jac}. Assume that the logarithmic model of $\mathcal F$ is non resonant, this implies condition (i) in the lemma above (by Proposition 4.4 in \cite{Cor-2003}) and hence all the divisors in $G(Z)$ are non-collinear. Note that $E^1$ is always a bifurcation divisor in $G(Z)$ and we have that $\sum_{R \in E^1} \Delta_{E^1}(R)=0$ by Remark \ref{rmk:suma-delta}.
If we write $\pi_{E^1}^* C \cap E_{red}^1= \{R_1^{E^1},\ldots, R_{b_{E^1}^C}^{E^1}\}$ and $\pi_{E^1}^* D \cap E_{red}^1 = \{Q\}$ with $R_{l}^{E^1}=(0,c_l^{E^1})$, $l=1,\ldots, b_{E^1}^C$, and $Q=(0,d)$ in coordinates in the first chart of $E_{red}^1$, then the set of zeros of ${\mathcal M}_{E^1}(z)$ are given by the roots of the polynomial
$$
\prod_{j =1}^{b_{E^1}^C} (z-c_j^{E^1}) + (z-d) \sum_{l =1}^{b_{E^1}^C} \mathcal{I}_{R_l^{E^1}}(\pi_{E^1}^* \mathcal{F},E_{red}^1) \prod_{j \neq l} (z -c_{j}^{E^1}) $$
which has multiplicity equal to $b_{E^1}^C-1$ provided that $$\sum_{j =1}^{b_{E^1}^C} c_j^{E^1} + d \sum_{l =1}^{b_{E^1}^C} \mathcal{I}_{R_l^{E^1}}(\pi_{E^1}^* \mathcal{F},E_{red}^1) \neq 0.$$
Note that we can assume that this condition holds since we are consider a generic polar curve.

Consider a component $E$ of the exceptional divisor $\pi_{C}^{-1}(0)$. We have that $\operatorname{NCol}(E)=\pi_{E}^* C \cap E_{red}$ if $E \neq E^1$ and $\operatorname{NCol}(E^1)=(\pi_{E^1}^* C \cap E_{red}^1) \cup (\pi_{E^1}^* D \cap E_{red}^1)$. Given  a point $P$ in $\pi_{E}^* C \cap E_{red}$, by Theorem \ref{th-multiplicidad-jac}, we have that
$$\nu_P( \pi_E^* {\mathcal P}^\mathcal{F})=\nu_P(\pi_E^* C)-1$$
and hence we have Theorem 6.1 in \cite{Cor-2003}. If we take the point $Q$ given by $\pi_{E^1}^* D \cap E_{red}^1$ we have that $\nu_Q( \pi_{E^1}^* {\mathcal P}^\mathcal{F})=\nu_Q(\pi_{E^1}^* D)-1=0$. Thus,
 by Theorem \ref{th:descomposicion-no-sing} we obtain that ${\mathcal P}^\mathcal{F}=\cup_{E \in B(C)} J^E$ with the following properties
\begin{itemize}
  \item[(1)] $\nu_0(J^E)=b_E^C-1$
  \item[(2)] $\pi_E^*J^E \cap \pi_E^* C=\emptyset$
  \item[(3)] if $E' < E$, then $\pi_E^* J^E \cap \pi_E'(E')=\emptyset$
  \item[(4)] if $E' > E$, then $\pi_{E'}^* J^E \cap E_{red}'=\emptyset$
\end{itemize}
which in particular implies the decomposition of the generic polar curve given in Corollary 6.2 of \cite{Cor-2003} for $C$ with non-irreducible components. Thus the decomposition in the general case (Theorem 5.1 in \cite{Cor-2003}) follows from Theorem \ref{th:desc-general}.

\appendix

\section{Equisingularity data and Ramification}\label{ap:ramificacion}
The aim of this appendix is to explain the behaviour of plane curves and their invariants under the action of a ramification. Although some of these results can be found in \cite{Cor-2009} (Appendix B), we include them here for completeness.

\subsection{Equisingularity data}\label{ap:equisingularidad}
In subsection~\ref{Apendice:grafo-dual} we have introduced some notations concerning equisingularity of plane curves that will be used in the sequel. This appendix completes subsection~\ref{Apendice:grafo-dual} with more notations related with equisingularity data that have already been  used to prove some results or that will be useful in order to describe the effect of ramification over a plane curve.

Recall that $\pi_C : X_C \to (\mathbb{C}^2,0)$ is the minimal reduction of singularities of a curve $C=\cup_{i=1}^r C_i$. Given an irreducible component $E$ of $\pi_C^{-1}(0)$,
a {\em curvette\/} $\tilde{\gamma}$ of the divisor $E$ is a
non-singular curve transversal to $E$ at a non-singular point of
$\pi_{C}^{-1}(0)$. The projection $\gamma=\pi_{C}(\tilde{\gamma})$
is a germ of plane curve in $({\mathbb C}^2,0)$ and we say that
$\gamma$ is an {\em $E$-curvette}. We denote by $m(E)$ the multiplicity at
the origin of any $E$-curvette and by $v(E)$ the coincidence ${\mathcal C}(\gamma_E,\gamma_E')$ of two $E$-curvettes $\gamma_E, \gamma_E'$ which cut $E$ in different points.
Note that $v(E) <
v(E')$ if $E < E'$. Recall  that the {\em coin\-cidence\/}
${\mathcal C}(\gamma,\delta)$ between two irreducible curves
$\gamma$ and $\delta$ is defined as
\begin{equation}\label{eq:coincidencia}
{\mathcal C}(\gamma,\delta) = \sup_{\ 1 \leq i \leq \nu_0(\gamma)
\atop 1 \leq j \leq \nu_0(\delta)} \{
\text{ord}_{x}(y^{\gamma}_i(x)-y^{\delta}_j(x))\ \}
\end{equation}
where $\{y_{i}^{\gamma}(x)\}_{i=1}^{\nu_0(\gamma)}$,
$\{y_{j}^{\delta}(x)\}_{j=1}^{\nu_0(\delta)}$ are the Puiseux series
of $\gamma$ and $\delta$ respectively.

\begin{Remark}\label{rmk:coincidencia-mult-int}
Note that the coincidence ${\mathcal C}(\gamma,\delta)$ between two irreducible curves $\gamma$ and $\delta$ and the intersection multiplicity $(\gamma,\delta)_0$ of both curves at the origin are related as follows (see Merle \cite{Mer}, prop. 2.4): if $\{\beta_0,\beta_1,\ldots,\beta_g\}$ are the characteristic exponents of $\gamma$ and $\alpha$ is a rational number such that $\beta_q \leq \alpha < \beta_{q+1}$ ($\beta_{g+1}=\infty$), then the following statements are equivalent:
\begin{align*}
    \text{(i)}& \ \ {\mathcal C}(\gamma,\delta) =
    \frac{\alpha}{\nu_0(\gamma)} \hspace{8cm}\\
    \text{(ii)} & \ \ \frac{(\gamma,\delta)_0}{\nu_0(\delta)}
    =\frac{\bar{\beta}_q}{n_1 \cdots n_{q-1}} +
    \frac{\alpha-\beta_q}{n_1 \cdots n_q}
\end{align*}
where $\{(m_i,n_i)\}_{i=1}^g$ are the Puiseux pairs of $\gamma$ ($n_0=1$) and $\{\bar{\beta}_0,\bar{\beta}_1,\ldots, \bar{\beta}_g\}$ is a minimal system of generators of the semigroup $S(\gamma)$ of $\gamma$.
\end{Remark}

Consider any curvette $\tilde{\gamma}_E$ of $E$, then $\pi_{E}'(\tilde{\gamma}_E)$ is also
a curvette of $E_{red} \subset X_E$ and it is clear  that $m(E)=m(E_{red})$ and $v(E)=v(E_{red})$. Let $\{ \beta_0^E,\beta_1^E, \ldots, \beta_{g(E)}^{E}\}$ be the characteristic exponents of
$\gamma_E = \pi_{C}(\tilde{\gamma}_E)$. Then we have that $m(E)=\beta_0^E=\nu_0(\gamma_{E})$. There are two possibilities for the value of $v(E)$:
\begin{itemize}
    \item[1.] either $\pi_E$ is the minimal reduction of
    singularities of $\gamma_E$ and then
    $v(E)=\beta_{g(E)}^E/\beta_0^E$. We say that $E$ is a {\em
    Puiseux divisor\/} for $\pi_C$ (or $C$);
    \item[2.] or $\pi_E$ is obtained by blowing-up $q \geq 1$ times
    after the minimal reduction of singularities of $\gamma_E$ and
    in this situation $v(E) = (\beta_{g(E)}^{E}+q )
    /\beta_0^E$. In this situation, if $E$ is a bifurcation divisor, we say that $E$ is a {\em contact divisor\/} for
    $\pi_C$ (or $C$).
\end{itemize}
 Moreover, a bifurcation divisor $E$
can belong to a dead arc only if it is a Puiseux divisor.

Take $E$ a bifurcation divisor of $G(C)$ and let $\{(m_1^E,n_1^E), (m_2^E,n_2^E), \ldots,
 (m_{g(E)}^E,n_{g(E)}^E)\}$ be the Puiseux pairs of an $E$-curvette $\gamma_E$, we denote
$$n_E =
\left\{
  \begin{array}{ll}
    n_{g(E)}, & \hbox{ if $E$ is a Puiseux divisor;} \\
    1, & \hbox{otherwise,}
  \end{array}
\right.
$$
and $\underline{n}_E=m(E)/n_E$. Observe that, if $E$ belongs to a
dead arc with terminal divisor $F$, then $m(F)=\underline{n}_E$. We
define $k_E$ to be
$$
k_E =\left\{
        \begin{array}{ll}
          g(E)-1, & \hbox{if $E$ is a Puiseux divisor;} \\
          g(E), & \hbox{if $E$ is a contact divisor.}
        \end{array}
      \right.
$$
Hence we have that $\underline{n}_E=n_1^E \cdots n_{k_E}^E$.
\begin{Remark}\label{rk:bif-div}
Let $E$ be a bifurcation divisor  of $G(C)$ which is a Puiseux divisor for $C$ and take any $i \in I_E$ (that is, $E$ belongs to the geodesic of the curve $C_i$). We have two possibilities concerning $v(E)$:
\begin{itemize}
  \item either $v(E)=\beta_{k_E+1}^i/\beta_0^i$, then we say that $E$ is a {\em Puiseux divisor for\/} $C_i$;
  \item or $v(E)$ corresponds to the coincidence of $C_i$ with another branch of $C$, in this situation we say that $E$ is a {\em contact divisor for\/} $C_i$.
\end{itemize}
Note that, if $E$ is a Puiseux divisor for $C$, then it is a Puiseux divisor for at least one irreducible component $C_i$ with $i \in I_E$, but it can be a contact divisor for other branches $C_j$ with $j \in I_E$, $j\neq i$. Consider for instance the curve $C=C_1 \cup C_2$ with $C_1=(y^2-x^3=0)$ and $C_2=(y-x^2=0)$. The dual graph $G(C)$ is given by

\vspace*{0.5\baselineskip}
\begin{center}
\begin{texdraw}
\arrowheadtype t:F \arrowheadsize l:0.08 w:0.04 \drawdim mm \setgray
0 \move(10 20) \fcir f:0 r:0.8 \rlvec (10 0) \fcir f:0 r:0.8 \rlvec
(0 -7) \fcir f:0 r:0.8 \move (20 20) \avec (25 23) \move (26 21)
\htext{\tiny{$C_1$}}

\move (20 13)  \avec (25 16)

\move (9 16) \htext{\tiny$E_1$} \move (21 17) \htext{\tiny{$E_3$}}
\move (16 10) \htext{\tiny{$E_2$}} \move (26 14) \htext{\tiny{$C_2$}}
\end{texdraw}
\end{center}
Thus the divisor $E_3$ is a Puiseux divisor for $C$ and $C_1$ but it is a contact divisor for $C_2$ since $v(E_3)=3/2={\mathcal C}(C_1,C_2)$.
\end{Remark}
\begin{Example}
  If we consider a semiroot $C^{(k)}$ of a curve $C$ as in subsection 7.2, then all the bifurcation divisors $E_1,E_2, \ldots, E_g$ are Puiseux divisors for $C$; the divisors $E_1,E_2, \ldots, E_k$ are also Puiseux divisors for $C^{(k)}$ but $E_{k+1}$ is a contact divisor for $C^{(k)}$.
\end{Example}

\subsection{Ramification}\label{subsec:ramification}
Consider a plane curve $C=\cup_{i=1}^r C_i$ in  $({\mathbb
C}^2,0)$. Let $\rho: ({\mathbb C}^2,0) \rightarrow ({\mathbb
C}^2,0)$ be any $C$-ramification, that is, $\rho$ is transversal to
C and  $\widetilde{C} = \rho^{-1}C$ has only non-singular irreducible
components. Consequently, if the ramification is given by $x=u^n, y=v$,  then it is
required that $n \equiv 0 \mod(n^1,n^2,\ldots,n^r)$  to
ensure that $\widetilde{C}$ has only non-singular irreducible components where $n^{i}=\nu_0(C_i)$.
Each curve $\widetilde{C}_i=\rho^{-1}C_i$ has exactly $n^i$ irreducible components and
the number of irreducible components of $\widetilde{C}$ is
equal to $\nu_0(C)= n^1+\cdots+n^r$.

More precisely, let $y^i(x) = \sum_{l \geq n^i} a_{l}^{i} x^{l/n^{i}}$ be a Puiseux
series of $C_i$, thus all its Puiseux series are given by
$$
y_{j}^{i}(x) = \sum_{l \geq n^{i}} a_{l}^{i} \varepsilon_i^{lj}
x^{l/n^{i}} \ \ \ \text{ for }j=1,2,\ldots, n^i,
$$
where $\varepsilon_i$ is a primitive $n^i$-root of the unity. Then
$f_i(x,y)=\prod_{l=1}^{n^{i}} (y-y_{l}^{i}(x))$ is a reduced
equation of $C_i$. If we put $v_j^{i}(u)=y_{j}^{i}(u^n)$, then
$v_{j}^{i}(u) \in {\mathbb C}\{u\}$ since $n/n^{i} \in {\mathbb N}$.
Hence the curve $\sigma_{j}^{i}=(v-v_j^{i}(u)=0)$ is
non-singular and it is one of the irreducible components of
$\rho^{-1}C_i$. Thus an equation of $\rho^{-1}C_i$ is given by
$$
g_i(u,v)= f_i(u^n,v)=\prod_{l=1}^{n^{i}} (v-v_{l}^{i}(u)).
$$
In particular we have that the irreducible
components $\{\sigma_{j}^{i}\}_{j=1}^{n^{i}}$ of $\rho^{-1}C_i$ are
in bijection with the Puiseux series of $C_i$.

It is well-known that the equisingularity type of a curve $C$ is
determined by the characteristic exponents
$\{\beta_0^{i},\beta_1^{i},\ldots,\beta_{g_i}^{i}\}_{i=1}^{r}$ of
its irreducible components and the intersection multiplicities
$\{(C_i,C_j)_0\}_{i\neq j}$. In \cite{Cor-2009} it is proved that
the equisingularity data of $C$ can be recovered from the curve $\rho^{-1}C$.

\vspace{\baselineskip}

Let us explain now  the relationship between the dual graphs $G(C)$ and $G(\widetilde{C})$ of the minimal reduction of singularities of $C$ and $\widetilde{C}$ respectively.
Observe   that, if
$\widetilde{E}$ and $\widetilde{E}'$ are two consecutive vertices of
$G(\widetilde{C})$ with $\widetilde{E} < \widetilde{E}'$, then
$v(\widetilde{E}')=v(\widetilde{E})+1$. Thus, $G(\widetilde{C})$ is completely
determined once we know  the bifurcation divisors, the order
relations  among them and the number of edges   which leave from
each bifurcation divisor.

Let $K_i$ be the geodesic in $G(C)$ of a branch $C_i$ of $C$ and let
$\widetilde{K}_i$ be the sub-graph of $G(\widetilde{C})$ corresponding to
the geodesics of the irreducible components
$\{\sigma_{l}^{i}\}_{l=1}^{n^i}$ of $\rho^{-1}C_i$. Let us explain how
to construct $\widetilde{K}_i$ from $K_i$.  Denote by $B(\widetilde{K}_i)$ and $B(K_i)$
the bifurcation vertices of $\widetilde{K}_i$ and $K_i$ respectively. We
say that a vertex $\widetilde{E}$ of $B(\widetilde{K}_i)$ is {\em associated\/} to a
vertex $E$ of $B(K_i)$ if $v(\widetilde{E})=n v(E)$. Note that there can be other bifurcation vertices in $G(\widetilde{C}) \smallsetminus B(\widetilde{K}_i)$ with valuation equal to $nv(E)$ but they are not associated to $E$.

Take a vertex $E$ of $B(K_i)$. Consider first the case of $E$ being the first
bifurcation divisor of $B(K_i)$ and take $E'$ its consecutive vertex
in $B(K_i)$. Then $E$ has only one associated vertex $\widetilde{E}$ in
$B(\widetilde{K}_i)$ and there are two possibilities for the number of
edges which leave from it:
\begin{itemize}
    \item If $E$ is a Puiseux divisor for $C_i$, %$v(E)={\beta_1^i}/{\beta_0^i}$ (we say that $E$ is a Puiseux divisor for $C_i$),
    then there are $n_1^{i}$ edges which leave from $\widetilde{E}$ in
    $\widetilde{K}_i$; then  $E'$ has $n_1^i$ associated
    vertices in $B(\widetilde{K}_i)$.
    \item If $E$ is a contact divisor for $C_i$, then
    there is only one edge which leave from $\widetilde{E}$ in $\widetilde{K}_i$ and thus
    $E'$ has only one associated  vertex in $B(\widetilde{K}_i)$.
\end{itemize}
Note that, if $E$ is a Puiseux divisor for $C$, then $E$ is a Puiseux divisor for at least one irreducible component $C_i$ but it can be a contact divisor for all the other irreducible components (see Remark~\ref{rk:bif-div}).
Let $E$ now be any vertex   of $B(K_i)$ and assume that we know the part
of $\widetilde{K}_i$ corresponding to the vertices of $K_i$ with
valuation $\leq v(E)$. Then there are $\underline{n}_E= n_1^i \cdots
n_{k_E}^i$ vertices $\{ \widetilde{E}^{l} \}_{l=1}^{\underline{n}_E}$
associated to $E$ and
\begin{itemize}
    \item If $E$ is a Puiseux divisor for $C_i$,
    then there are $n_{k_E+1}^i$ edges which leave from each
    vertex $\widetilde{E}_l$ in $\widetilde{K}_i$.
    \item If $E$ is a contact divisor for $C_i$,
then there is only one edge which leaves from each vertex $\widetilde{E}_l$
in $\widetilde{K}_i$.
\end{itemize}
The dual graph $G(\widetilde{C})$ is constructed   by
gluing the graphs $\widetilde{K}_i$. Thus we deduced that, if $\widetilde{E}$ is a divisor of $G(\widetilde{C})$ associated to a divisor $E$ of $G(C)$, then
\begin{equation}\label{eq:b-tilde}
b_{\widetilde{E}} = \left\{
                  \begin{array}{ll}
                    b_E, & \hbox{if $E$ is a contact divisor for $C$;} \\
                    (b_E-1)n_E, & \hbox{if $E$ is a Puiseux divisor for $C$ which
                                belongs} \\
                                & \hbox{to a dead arc;} \\
                    (b_E-1)n_E +1 , & \hbox{if $E$ is a Puiseux divisor for $C$ which does
                                        not} \\
                                 & \hbox{belong to a dead arc.}
                  \end{array}
                \right.
\end{equation}
Thus all vertices in $G(\widetilde{C})$ associated to a divisor $E$ of $G(C)$ have the same valence.
Moreover, if $\gamma_E$ is an $E$-curvette of a
bifurcation divisor $E$ of $G(C)$, the curve $\rho^{-1}\gamma_E$ has
$m(E)=\underline{n}_E n_E$ irreducible components which are all
non-singular and each divisor $\widetilde{E}^l$ belongs to the geodesic
of exactly $n_E$ branches of $\rho^{-1}\gamma_E$ which are curvettes
of $\widetilde{E}^l$ in different points.

Observe that there are non-bifurcation divisors  of $G(C)$ without
associated divisors in $G(\widetilde{C})$.

Due to the bijection   between the Puiseux series of
$C_i$ and the irreducible components of $\rho^{-1}C_i$, we have that
 the choice of a vertex $\widetilde{E}^l
\in B(\widetilde{K}_i)$ associated to a bifurcation divisor $E$ is
equivalent to the choice of a $\underline{n}_{E}$-th root $\xi_l$ of the
unity. This implies that the vertex $\widetilde{E}^l$ belongs to the geodesic of $e_{E}^{i}= n^{i}/\underline{n}_E$
irreducible components $\{\sigma_{lt}^{i}\}_{t=1}^{e_{E}^{i}}$ of
$\rho^{-1}C_i$.
Moreover, the curve $\sigma_{lt}^{i}$ is given by
$(v-\eta_{lt}^{i}(u)=0)$ where
$$
\eta_{lt}^{i}(u)=\sum_{s \geq n^i} a_{s}^{i} (\zeta_{ilt})^s
u^{sn/n^i}, \text{ for } t=1, \ldots, e_{E}^{i},
$$
and $\{\zeta_{ilt}\}_{t=1}^{e_{E}^{i}}$ are the $e_{E}^{i}$-th roots
of $\xi_l$. The cardinal of the set $\pi_{\widetilde{E}^l}^* \widetilde{C}_i \cap \widetilde{E}^l_{red}$ is given by
\begin{equation*}
\sharp(\pi_{\widetilde{E}^l}^* \widetilde{C}_i \cap \widetilde{E}^l_{red})= \left\{
  \begin{array}{ll}
    n_E, & \hbox{ if $E$ is a Puiseux divisor for $C_i$;} \\
    1, & \hbox{ if $E$ is a contact divisor for $C_i$.}
  \end{array}
\right.
\end{equation*}
Furthermore, if $E$ is a Puiseux divisor for $C_i$ and $P, Q$ are two different points in $\pi_{\widetilde{E}^l}^* \widetilde{C}_i \cap \widetilde{E}^l_{red}$, we have that
\begin{equation}\label{eq:mult-C_i}
  \nu_P(\pi_{\widetilde{E}^l}^* \widetilde{C}_i) = \nu_Q(\pi_{\widetilde{E}^l}^* \widetilde{C}_i)=\frac{e_E^i}{n_E};
\end{equation}
and if $E$ is a contact divisor for $C_i$ and we denote $P$ the only point in the set $\pi_{\widetilde{E}^l}^* \widetilde{C}_i \cap \widetilde{E}^l_{red}$, then
\begin{equation}\label{eq:mult-C_i-2}
\nu_P(\pi_{\widetilde{E}^l}^* \widetilde{C}_i)=e_E^i.
\end{equation}

Consider now two divisors $\widetilde{E}^l$ and $\widetilde{E}^k$   associated to the same bifurcation divisor $E$ of $G(C)$, and let  $\xi_l$ and $\xi_k$ be the $\underline{n}_E$-th roots of the unity corresponding to the divisors $\widetilde{E}^l$ and $\widetilde{E}^k$.
We can define a bijection
$\rho_{l,k}: \widetilde{E}^l_{red} \to \widetilde{E}^k_{red}$ as follows:
 the map $\rho_{l,k}$ sends the
``infinity point" of $\widetilde{E}^l_{red}$ (that is, the origin of the
second chart of $\widetilde{E}^l_{red}$) into  the ``infinity point" of
$\widetilde{E}^k_{red}$.
Given another point $P$ of $\widetilde{E}^l_{red}$, we
take an $\widetilde{E}^l$-curvette
$\gamma_{\widetilde{E}^l}^P=(v-\psi_{\widetilde{E}^l}^{P}(u)=0)$ with
\begin{equation}\label{eq:curvette-tildeEl}
\psi_{\widetilde{E}^l}^{P}(u) = \sum_{i=1}^{v(\widetilde{E}^l)-1}
a_i^{\widetilde{E}^l} u^{i} + a_{v(\widetilde{E}^l)}^{P} u^{v(\widetilde{E}^l)},
\end{equation}
and such that $\pi_{\widetilde{E}^l}^* \gamma_{\widetilde{E}^l}^P \cap
\widetilde{E}^l_{red} = \{P\}$.
We define $\rho_{l,k}(P)$ to be the point $\pi_{\widetilde{E}^k}^* \gamma_{\widetilde{E}^k}^{\rho_{l,k}(P)} \cap
\widetilde{E}^k_{red}$, where the curve
$\gamma_{\widetilde{E}^k}^{\rho_{l,k}(P)}= (v-\psi_{\widetilde{E}^k}^{\rho_{l,k}(P)}(u)=0)$ is given by
$$\psi_{\widetilde{E}^k}^{\rho_{l,k}(P)}(u)=\sum_{i=1}^{v(\widetilde{E}^l)-1}
a_i^{\widetilde{E}^l} \left(\frac{\xi_{k}}{\xi_l}\right)^i  u^{i} + a_{v(\widetilde{E}^l)}^{P} \left(\frac{\xi_{k}}{\xi_l}\right)^{v(\widetilde{E}^l)} u^{v(\widetilde{E}^l)}.$$
Note that $\gamma_{\widetilde{E}^k}^{\rho_{l,k}(P)}$ is an $\widetilde{E}^k$-curvette.

Recall also that given any bifurcation divisor $E$ of $G(C)$ or $E=E^1$ and any of its associated divisors $\widetilde{E}^l$ in $G(\widetilde{C})$, there is a morphism  $\rho_{\widetilde{E}^l,E}:\widetilde{E}^l_{red} \to E_{red}$ which is a ramification of order $n_E$ (see \cite{Cor-2009}, Lemma 8). The map $\rho_{\widetilde{E}^l,E}$ is defined as follows:   $\rho_{\widetilde{E}^l,E}$ sends the
``infinity point" of $\widetilde{E}^l_{red}$
%(that is, the origin of the second chart of $\widetilde{E}^l_{red}$)
into  the ``infinity point" of
$E_{red}$ and the origin of the first chart of $\widetilde{E}^l_{red}$ is sent to the origin of the first chart of $E_{red}$. For any other point $P$ of $\widetilde{E}^l_{red}$, we
can take an $\widetilde{E}^l$-curvette
$\gamma_{\widetilde{E}^l}^P=(v-\psi_{\widetilde{E}^l}^{P}(u)=0)$ with $\pi_{\widetilde{E}^l}^* \gamma_{\widetilde{E}^l}^P \cap
\widetilde{E}^l_{red} = \{P\}$. Thus if
$\psi_{\widetilde{E}^l}^{P}(u)$ is given by equation~\eqref{eq:curvette-tildeEl}, we can consider the $E$-curvette $\gamma_E^P$ given by the Puiseux series
$$
y^{P}(x) = \sum_{i=1}^{v(\widetilde{E}^l)-1}
a_i^{\widetilde{E}^l} x^{i/m(E)} + a_{v(\widetilde{E}^l)}^{P} x^{v(\widetilde{E}^l)/m(E)},
$$
and $\rho_{\widetilde{E}^l,E}(P)$ is the only point $\pi^*_E \gamma_E^P \cap E_{red}$. Observe that, if $E$ is a bifurcation divisor in the geodesic of a branch $C_i$ of $C$ and $\widetilde{E}^l$ is any associated divisor to $E$ in $G(\widetilde{C})$, then the morphism $\rho_{\widetilde{E}^l,E}$ maps all the points in $\pi_{\widetilde{E}^l}^* \widetilde{C}_i \cap \widetilde{E}^l_{red}$ to the only point in $\pi_E^*C_i \cap E_{red}$.

Moreover, note that, if $\widetilde{E}^l$ and $\widetilde{E}^k$ are two divisors of $G(\widetilde{C})$ associated to a bifurcation divisor $E$ of $G(C)$, then the following diagram
$$
\xymatrix{
\widetilde{E}^l_{red} \ar[rr]^{\rho_{l,k}} \ar[rd]_{\rho_{\widetilde{E}^l,E}}  &  &  \widetilde{E}^{k}_{red}  \ar[ld]^{\rho_{\widetilde{E}^k,E}} \\
  & E_{red}  &
}
$$
is commutative.

Finally remark that, if $\gamma_{E_t}$ is a curvette of a terminal
divisor $E_t$ of a dead arc with bifurcation divisor $E$, then
$\rho^{-1}\gamma_{E_t}$ is composed by $m(E_t)=\underline{n}_{E}$
non-singular irreducible components and each divisor $\widetilde{E}^l$
belongs to the geodesic of exactly  one branch of
$\rho^{-1}\gamma_{E_t}$, where
$\{\widetilde{E}^l\}_{l=1}^{\underline{n}_E}$ are the divisors
associated to $E$ in $G(\widetilde{C})$.

\subsection{Logarithmic foliations and ramification}\label{subsec:ram-log}
Consider the logarithmic foliation $\mathcal{L}^C_\lambda$ defined by
$$
f_1 \cdots f_r \sum_{i=1}^{r} \lambda_i \frac{df_i}{f_i}=0
$$
with $\lambda=(\lambda_1,\ldots,\lambda_r) \in \mathbb{C}^r$  and $f_i \in \mathbb{C}\{x,y\}$ (see section~\ref{sec:logarithmic} for notations concerning logarithmic foliations). Let us see the behaviour of $\mathcal{L}^C_\lambda$ after a ramification.

Consider the curve $C=\cup_{i=1}^r C_i$ with $C_i=(f_i=0)$ and let $\rho: ({\mathbb C}^2,0) \to ({\mathbb C}^2,0)$ be any $C$-ramification, that is, $\rho$ is transversal to $C$ and the curve $\widetilde{C}=\rho^{-1}C$ has only non-singular irreducible components. We refer to subsections~\ref{Apendice:grafo-dual}, \ref{ap:equisingularidad} and \ref{subsec:ramification} for notations concerning equisingularity data of curves and ramifications.

We have that $\rho^*{\mathcal L}_{\lambda}^C = {\mathcal L}_{\lambda^*}^{\widetilde{C}}$ with
$$
\lambda^* = (\overbrace{\lambda_{1}, \ldots,\lambda_{1}}^{n^1},
\ldots, \overbrace{\lambda_{r}, \ldots, \lambda_{r}}^{n^{r}})
$$
where $n^i=\nu_0(C_i)$ for $i=1,\ldots,r$. We put $\rho^{-1}C_i= \widetilde{C}_i=\{ \sigma_t^i\}_{t=1}^{n^i}$ where each $\sigma_t^i$ is an irreducible curve. Moreover, we have that logarithmic models behave well under ramification. More precisely,
\begin{Proposition}\label{prop:ram-modelo-log}(see \cite{Cor-2003})
Let $\mathcal{F} \in \mathbb{G}_C$ and $\mathcal{L}_\lambda^C$ a logarithmic model of $\mathcal{F}$. If $\rho: ({\mathbb C}^2,0) \to ({\mathbb C}^2,0)$ is a $C$-ramification, then $\rho^*\mathcal{L}_\lambda^C$ is a logarithmic model of $\rho^*\mathcal{F}$.
\end{Proposition}
Let $\pi_C : X_C \to ({\mathbb C}^2,0)$ be the minimal reduction of singularities of $C$ and  $\pi_{\widetilde{C}} : X_{\widetilde{C}} \to ({\mathbb C}^2,0)$ be the minimal reduction of singularities of $\widetilde{C}$. Take $E$ a bifurcation divisor of $G(C)$ and let $\widetilde{E}^l$ be any bifurcation divisor of $G(\widetilde{C})$ associated to $E$.
Given any $i \in I_E$ (that is, $E$ is in the geodesic of $C_i$), there are $e_E^i$ branches of $\rho^{-1}C_i$ such that $\widetilde{E}^l$ belongs to their geodesics where $e_E^i=n^i/\underline{n}_E$  and we have that the residue of the logarithmic foliation ${\mathcal L}_{\lambda^*}^{\widetilde{C}}$ along the divisor $\widetilde{E}^l$ is given by
$$
\kappa_{\widetilde{E}^l}({\mathcal L}_{\lambda^*}^{\widetilde{C}})=\sum_{i=1}^{r} \lambda_i \sum_{t=1}^{n^i} \sum_{\widetilde{E} \leq \widetilde{E}^l} \varepsilon_{\widetilde{E}}^{\sigma_t^i}.
$$
(see equation~\eqref{eq:peso-divisor-log} for its definition).
Let $\{R_1^{\widetilde{E}^l}, R_2^{\widetilde{E}^l}, \cdots, R_{b_{\widetilde{E}^l}}^{\widetilde{E}^l}\}$ be the set of points $\pi_{\widetilde{E}^l}^* \widetilde{C} \cap \widetilde{E}^l_{red}$ and put $m_{R_t^{\widetilde{E}^l}}^i=\nu_{R_t^{\widetilde{E}^l}} ( \pi_{\widetilde{E}^l}^* \widetilde{C}_i)$ for $t=1,2, \ldots, b_{\widetilde{E}^l}$. Note that $m_{R_t^{\widetilde{E}^l}}^i= \sharp \{ s \in \{1,\ldots, n^i\} \ : \ \pi_{\widetilde{E}^l}^* \sigma_s^i \cap \widetilde{E}^l_{red}= \{ R_t^{\widetilde{E}^l} \} \} =\frac{e_E^i}{n_E}$ (the last equality follows from equations~\eqref{eq:mult-C_i} and \eqref{eq:mult-C_i-2} in appendix~\ref{subsec:ramification} where $m_{R_t^{\widetilde{E}^l}}^i$ is also computed). With these notations we have that
\begin{equation}\label{eq:camacho-sad-log-ram}
\mathcal{I}_{R_t^{\widetilde{E}^l}}(\pi_{\widetilde{E}^l}^* {\mathcal L}_{\lambda^*}^{\widetilde{C}},\widetilde{E}^l_{red})= -\frac{\ \  \sum_{i \in I_E} \lambda_i m_{R_t^{\widetilde{E}^l}}^i\ \ }{\kappa_{\widetilde{E}^l}({\mathcal L}_{\lambda^*}^{\widetilde{C}})}.
\end{equation}
Observe that if $\widetilde{E}^l$ and $\widetilde{E}^k$  are two bifurcation divisors of $G(\widetilde{C})$ associated to the same divisor $E$ of $G(C)$, we have that $\kappa_{\widetilde{E}^l}({\mathcal L}_{\lambda^*}^{\widetilde{C}}) = \kappa_{\widetilde{E}^k}({\mathcal L}_{\lambda^*}^{\widetilde{C}})$. Moreover,
there is a bijection between the sets of points $\pi_{\widetilde{E}^l}^* \widetilde{C} \cap \widetilde{E}^l_{red}$ and $\pi_{\widetilde{E}^k}^* \widetilde{C} \cap \widetilde{E}^k_{red}$ induced by the map $\rho_{l,k}: \widetilde{E}^l_{red} \to \widetilde{E}^k_{red}$ (see appendix~\ref{subsec:ramification}). Hence, if  $\{R_1^{\widetilde{E}^k}, R_2^{\widetilde{E}^k}, \cdots, R_{b_{\widetilde{E}^k}}^{\widetilde{E}^k}\}$ is the set of points $\pi_{\widetilde{E}^k}^* \widetilde{C} \cap \widetilde{E}^k_{red}$ with $R_t^{\widetilde{E}^k}= \rho_{l,k}(R_t^{\widetilde{E}^l})$ for $t=1,2,\ldots, b_{\widetilde{E}^k}$, we have that
\begin{equation}\label{eq:indices-ramificados}
\mathcal{I}_{R_t^{\widetilde{E}^k}}(\pi_{\widetilde{E}^k}^* {\mathcal L}_{\lambda^*}^{\widetilde{C}},\widetilde{E}^k_{red}) = \mathcal{I}_{R_t^{\widetilde{E}^l}}(\pi_{\widetilde{E}^l}^* {\mathcal L}_{\lambda^*}^{\widetilde{C}},\widetilde{E}^l_{red}), \qquad t=1,2,\ldots, b_{\widetilde{E}^k}.
\end{equation}
Moreover, if $R_t^{\widetilde{E}^l}, R_s^{\widetilde{E}^l}$ are two points in $\pi_{\widetilde{E}^l}^* \widetilde{C} \cap \widetilde{E}^l_{red}$ with $\rho_{\widetilde{E}^l,E}(R_t^{\widetilde{E}^l}) =\rho_{\widetilde{E}^l,E}(R_s^{\widetilde{E}^l})$ where $\rho_{\widetilde{E}^l,E}: \widetilde{E}^l_{red} \to E_{red}$ is the ramification defined in appendix~\ref{subsec:ramification}, then
\begin{equation}\label{eq:indices-ramificados-mismo-divisor}
\mathcal{I}_{R_t^{\widetilde{E}^l}}(\pi_{\widetilde{E}^l}^* {\mathcal L}_{\lambda^*}^{\widetilde{C}},\widetilde{E}^l_{red}) =\mathcal{I}_{R_s^{\widetilde{E}^l}}(\pi_{\widetilde{E}^l}^* {\mathcal L}_{\lambda^*}^{\widetilde{C}},\widetilde{E}^l_{red})
\end{equation}
since $m_{R_t^{\widetilde{E}^l}}^i=m_{R_s^{\widetilde{E}^l}}^i=\frac{e_E^i}{n_E}$ for $i \in I_E$ by equations~\eqref{eq:mult-C_i} and \eqref{eq:mult-C_i-2}. % in Appendix~\ref{ap:ramificacion}.

\section{Intersection multiplicities}

\medskip
We state now two results concerning the intersection multiplicity of the jacobian curve of two foliations either with a single separatrix of one of the foliations and with the curve of all separatrices. These intersection multiplicities are computed in terms of local invariants of the foliations (see subsection~\ref{subsec:inv-foliations} for notations).
Consider two foliations $\mathcal F$ and $\mathcal G$    in $({\mathbb C}^2,0)$ and denote by $\mathcal{J}_{\mathcal{F},\mathcal{G}}$ the jacobian curve of $\mathcal F$ and $\mathcal G$.

\begin{Proposition}\label{prop:int-sep}
Assume that $\mathcal F$ and $\mathcal G$ have no common separatrix.
If $S$ is an irreducible separatrix of $\mathcal F$, we have that
$$(\mathcal{J}_{\mathcal{F},\mathcal{G}},S)_0 = \mu_0({\mathcal F},S)+\tau_0({\mathcal G},S).
$$
\end{Proposition}
\begin{proof}
Let us write $\omega=A(x,y) dx + B(x,y)dy$ and $\eta=P(x,y) dx + Q(x,y) dy$ the 1-forms defining $\mathcal F$ and $\mathcal G$ respectively.
Let $\gamma(t)=(x(t),y(t))$ be a parametrization of the curve $S$. We can assume, without loss of generality, that $x(t) \neq 0$  and thus $\dot{x}(t) \neq 0$. Since $S$ is a separatrix of $\mathcal F$, then
$A(\gamma(t)) \dot{x}(t) + B(\gamma(t)) \dot{y}(t)=0$. Thus, we have that
\begin{align*}
  (\mathcal{J}_{\mathcal{F},\mathcal{G}},S)_0  & = \text{ord}_t\{A(\gamma(t)) Q(\gamma(t)) - B(\gamma(t)) P(\gamma(t))\} \\
   & = \text{ord}_t \left\{ \frac{-B(\gamma(t)) \dot{y}(t)}{\dot{x}(t)} Q(\gamma(t)) - B(\gamma(t)) P(\gamma(t)) \right\} \\
   & = \text{ord}_t (B(\gamma(t)))- (\text{ord}_t (x(t))-1) + \text{ord}_t \{P(\gamma(t)) \dot{x}(t) + Q(\gamma(t)) \dot{y}(t) \} \\
   &= \mu_0({\mathcal F},S)+\tau_0({\mathcal G},S)
\end{align*}
where the last equality comes from the expression of $\mu_0({\mathcal F},S)$ given in \eqref{multiplicidaderelativa1} and the definition of
$\tau_0({\mathcal G},S)$ given in \eqref{def-tangencia}.
\end{proof}
When  $\mathcal G$ is a non-singular foliation, we obtain Proposition 1 in \cite{Can-C-M} for the polar intersection number with respect to a branch of the curve of separatrices of $\mathcal F$. Note that, although in \cite{Can-C-M} it is assumed that the foliation $\mathcal F$ is non-dicritical, the result is also true when $\mathcal F$ is a dicritical foliation.
Using property (iv) in Theorem~\ref{th-curva-gen}, we get following consequence of the above result:
\begin{Corollary}
If $\mathcal G$ is a non-dicritical second type foliation and $S_{\mathcal G}$ is the curve of separatrices of $\mathcal G$,  we have that
$$(\mathcal{J}_{\mathcal{F},\mathcal{G}},S)_0 = \mu_0({\mathcal F},S)+(S_{\mathcal G},S)_0 -1.$$
\end{Corollary}
Next result gives a relationship among the intersection multiplicities of the jacobian curve with the curves of separatrices and the Milnor number of the foliations.
\begin{Proposition}\label{prop:milnor-number}
Consider two non-dicritical second type foliations ${\mathcal F}$ and ${\mathcal G}$ without common separatrices. Thus
$$({\mathcal J}_{{\mathcal F},{\mathcal G}},S_{\mathcal F})_0-({\mathcal J}_{{\mathcal F},{\mathcal G}},S_{\mathcal G})_0=\mu_0({\mathcal F})-\mu_0({\mathcal G}),$$
where $S_{\mathcal F}$, $S_{\mathcal G}$ are the curves of separatrices of $\mathcal F$ and $\mathcal G$ respectively.
\end{Proposition}
\begin{proof}
Let ${\mathcal B}({\mathcal J}_{{\mathcal F},{\mathcal G}})$ be the set of irreducible components of ${\mathcal J}_{{\mathcal F},{\mathcal G}}$. Given any branch $\Gamma \in {\mathcal B}({\mathcal J}_{{\mathcal F},{\mathcal G}})$, we denote by $\gamma_\Gamma(t)=(x_\Gamma(t),y_\Gamma(t))$ any primitive parametrization of $\Gamma$. Assume that the foliations $\mathcal F$ and $\mathcal G$ are defined by the 1-forms $\omega=A dx + Bdy$ and $\eta=P dx + Qdy$ respectively. Thus we have that $A(\gamma_\Gamma(t)) Q(\gamma_\Gamma(t)) - B(\gamma_\Gamma(t)) P(\gamma_\Gamma(t))=0$. Since $\Gamma$ is not a separatrix of $\mathcal G$, then either $Q(\gamma_\Gamma(t)) \not \equiv 0$ or $P(\gamma_\Gamma(t)) \not \equiv 0$. We will assume that $Q(\gamma_\Gamma(t)) \not \equiv 0$. Let us compute the intersection multiplicity of $({\mathcal J}_{{\mathcal F},{\mathcal G}},S_{\mathcal F})_0$ taking into account property (iv) in Theorem~\ref{th-curva-gen}:
\begin{align*}
 ({\mathcal J}_{{\mathcal F},{\mathcal G}},S_{\mathcal F})_0 %& =({\mathcal J}_{{\mathcal F},{\mathcal G}},C)_0 + ({\mathcal J}_{{\mathcal F},{\mathcal G}},D)_0  \\
&   =  \sum_{\Gamma \in {\mathcal B}({\mathcal J}_{{\mathcal F},{\mathcal G}})} (\Gamma,S_{\mathcal F})_0 = \sum_{\Gamma \in {\mathcal B}({\mathcal J}_{{\mathcal F},{\mathcal G}})} ( \tau_0({\mathcal F},\Gamma) + 1)\\
 & = \sum_{\Gamma \in {\mathcal B}({\mathcal J}_{{\mathcal F},{\mathcal G}})} (\text{ord}_t \{A(\gamma_\Gamma(t)) \dot{x}_\Gamma(t) + B(\gamma_\Gamma(t)) \dot{y}_\Gamma(t) \} +1) \\
%& + \sum_{\Gamma \in {\mathcal B}({\mathcal J}_{{\mathcal F},{\mathcal G}})} \text{ord}_t \{P(\gamma_\Gamma(t)) \dot{x}_B(t) + Q(\gamma_\Gamma(t)) \dot{y}_B(t)  +1\}  \\
& = \sum_{\Gamma \in {\mathcal B}({\mathcal J}_{{\mathcal F},{\mathcal G}})} \left(\text{ord}_t \left\{ \frac{B(\gamma_\Gamma(t)) P(\gamma_\Gamma(t))}{Q(\gamma_\Gamma(t))} \dot{x}_\Gamma(t) + B(\gamma_\Gamma(t)) \dot{y}_\Gamma(t) \right\}+1  \right) \\
& = \sum_{\Gamma \in {\mathcal B}({\mathcal J}_{{\mathcal F},{\mathcal G}})} \text{ord}_t\{B(\gamma_\Gamma(t))\} -  \sum_{\Gamma \in {\mathcal B}({\mathcal J}_{{\mathcal F},{\mathcal G}})} \text{ord}_t\{Q(\gamma_\Gamma(t))\}
\\ &+   \sum_{\Gamma \in {\mathcal B}({\mathcal J}_{{\mathcal F},{\mathcal G}})} (\text{ord}_t \{P(\gamma_\Gamma(t)) \dot{x}_\Gamma(t) + Q(\gamma_\Gamma(t)) \dot{y}_\Gamma(t)\}  +1) \\
& = \mu_0({\mathcal F}) - \mu_0({\mathcal G}) + \sum_{\Gamma \in {\mathcal B}({\mathcal J}_{{\mathcal F},{\mathcal G}})} ( \tau_0({\mathcal G},\Gamma) + 1)\\
%\\ &+ 2  \sum_{\Gamma \in {\mathcal B}({\mathcal J}_{{\mathcal F},{\mathcal G}})} \text{ord}_t \{P(\gamma_\Gamma(t)) \dot{x}_B(t) + Q(\gamma_\Gamma(t)) \dot{y}_B(t)  +1\} \\
& = \mu_0({\mathcal F}) - \mu_0({\mathcal G}) + ({\mathcal J}_{{\mathcal F},{\mathcal G}},S_{\mathcal G})_0.
\end{align*}
\end{proof}
\bibliographystyle{abbrv}

%\bibliography{biblio-jacobian}
%\end{document}

\end{document}